\documentclass[11pt]{article}

\usepackage{amssymb,amsthm,amsmath,amsfonts}
\usepackage{graphicx}
\usepackage{algorithm,algorithmic}
\usepackage[usenames]{color}
\graphicspath{{./pics/}}
\usepackage[colorlinks,linktocpage,linkcolor=blue]{hyperref}

\newtheorem{exam}{Example}

\title{A Two-stage Method for Inverse Medium Scattering}
\author{Kazufumi Ito\footnote{Department of Mathematics and Center for
Research in Scientific Computation, North Carolina State University,
Raleigh, North Carolina (kito@unity.ncsu.edu).} \and Bangti
Jin\footnote{Department of Mathematics and Institute for Applied Mathematics
and Computational Science, Texas A\&M University, College Station, Texas
77843-3368, USA (btjin@math.tamu.edu).} \and Jun Zou\footnote{Department
of Mathematics, Chinese University of Hong Kong, Shatin, N.T., Hong Kong
(zou@math.cuhk.edu.hk).}}
\date{\today}
\begin{document}
\maketitle

\begin{abstract}
We present a novel numerical method to the time-harmonic inverse medium scattering
problem of recovering the refractive index from near-field scattered data. The approach
consists of two stages, one pruning step of detecting the scatterer support, and one
resolution enhancing step with mixed regularization. The first step is strictly direct
and of sampling type, and faithfully detects the scatterer support. The second step is an
innovative application of nonsmooth mixed regularization, and it accurately resolves the
scatterer sizes as well as intensities. The model is efficiently solved by a semi-smooth
Newton-type method. Numerical results for two- and three-dimensional examples indicate
that the approach is accurate, computationally efficient, and robust with respect to data
noise.
\end{abstract}

\section{Introduction}
In this work we study the inverse medium scattering problem (IMSP) of determining the
refractive index from near-field measurements for time-harmonic wave propagation
\cite{ColtonKress:1998}. Consider a homogeneous background space $\mathbb{R}^\mathrm{d}$
$(\mathrm{d}=2,3)$ that contains some inhomogeneous media occupying a bounded domain
$\Omega$. {Let $u^{inc}=e^{ik\,x\cdot d}$ be} an incident plane wave, with the incident
direction $d\in \mathbb{S}^{\mathrm{d}-1}$ and the wave number $k$. Then the total field
$u$ induced by the inhomogeneous medium scatterers satisfies the Helmholtz
equation~\cite{ColtonKress:1998}
\begin{equation}\label{eq:acoustic}
    \Delta u+k^2\,n^2(x) u=0,
\end{equation}
where the function $n(x)$ is the refractive index, i.e., the ratio of the wave speed in
the homogeneous background to that in the concerned medium at location $x$. The model
describes not only time-harmonic acoustic wave propagation, but also electromagnetic wave
propagation in either the transverse magnetic or transverse electric modes
\cite[Appendix]{ItoJinZou:2011a}.

Next we let $\eta=(n^2-1)k^2$, which combines the relative refractive index $n^2-1$ with
the wave number $k$. Clearly the coefficient $\eta$ characterizes the material properties
of the inhomogeneity and is supported in the scatterer $\Omega \subset
\mathbb{R}^\mathrm{d}$. We denote by $I=\eta u$ the induced current, by $G(x,y)$ the
fundamental solution to the Helmholtz equation in the homogeneous background, i.e.,
\begin{equation*}\displaystyle
   G(x,y)=\left\{
   \begin{aligned}\frac{i}{4}H^1_0(k\,|x-y|), &\quad\mathrm{d} =2,\\
   \frac{1}{4\pi}\frac{e^{ik|x-y|}}{|x-y|},&\quad \mathrm{d}=3,\\
   \end{aligned}\right.
\end{equation*}
where the function $H_0^1$ refers to Hankel function of the first kind and zeroth-order.
Then we can express the total field $u$ as follows \cite{ColtonKress:1998}
\begin{equation} \label{int}
  u=u^{inc}+\int_\Omega G(x,y)I(y)\,dy\,.
\end{equation}
By multiplying both sides of equation \eqref{int} by $\eta$, we arrive at the following
integral equation of the second kind for the induced current $I$
\begin{equation}\label{eqn:indcur}
   I(x) = \eta u^{inc} + \eta \int_\Omega G(x,y) I (y)dy.
\end{equation}
The reformulation \eqref{eqn:indcur} is numerically amenable since all computation is
restricted to the scatterer support $\Omega$, which is much smaller than the whole space
$\mathbb{R}^\mathrm{d}$. Hence, the complexity is also very low. We will approximate the
solution to \eqref{eqn:indcur} by the mid-point rule (cf. Appendix \ref{app:int}).

Next we let $u^s= u-u^{inc}$ be the scattered field, which is measured on a closed
curve/surface $\Gamma$ enclosing the scatterers $\Omega$. Then the IMSP is to retrieve
the refractive index $n^2$ or equivalently the coefficient $\eta$ from (possibly very
noisy) measurements of the scattered field $u^s$, corresponding to one or several
incident fields. In the literature, a number of reconstruction techniques for the IMSP
have been developed. These methods can be roughly divided into two groups: support
detection and coefficient estimate. The former group (including MUSIC
\cite{Devaney:2006,Cheney:2001}, linear sampling method \cite{ColtonKirsch:1996,
CakoniColtonMonk:2011} and factorization method \cite{Kirsch:1998} etc.) usually is of
sampling type, and aims at detecting the scatterer support efficiently. The latter group
generally relies on the idea of regularization (including Tikhonov regularization
\cite{BaoLi:2005,RachowiczZdunek:2011}, iterative regularization method
\cite{Hohage:2001, BakushinskyKokurinKozlov:2005, Hohage:2006}, contrast source inversion
\cite{BergBroekhovenAbubakar:1999}, subspace regularization \cite{ChenZhong:2009} and
propagation-backpropagation method \cite{Vogeler:2003}), and aims at retrieving a
distributed estimate of the index function. These approaches generally are more
expensive, but their results may profile the inhomogeneities more precisely.

In this paper, we shall develop a novel two-stage numerical method for solving the IMSP.
The first step employs a direct sampling method, recently developed in
\cite{ItoJinZou:2011a}, to detect the scatterer support $\Omega$ stably and accurately.
It is based on the following index function
\begin{equation} \label{indx}
  \Phi(x_p)=\frac{|\langle u^s, \,G(\cdot,x_p)\rangle_{L^2(\Gamma)}|}{
    \|u^s\|_{L^2(\Gamma)}\|G(\cdot,x_p)\|_{L^2(\Gamma)}}\quad \forall \,x_p\in\widetilde{\Omega},
\end{equation}
where $\widetilde{\Omega}\supset\Omega$ is a sampling domain. Numerically, the method is
strictly direct and does not incur any linear matrix operations. The method can detect
reliably the scatterer support $\Omega$ even in the presence of a large amount of data
noises \cite{ItoJinZou:2011a}. In particular, a (much smaller) computational domain
$D\subset\widetilde{\Omega}$ can be determined from the index $\Phi$, and furthermore the
restriction $\Phi|_D$ may serve as a first approximation to the coefficient $\eta$.

The second step enhances the image resolution by a novel application of (nonsmooth) mixed
regularization. With the approximation $\Phi|_D$ from the sampling step at hand, equation
\eqref{eqn:indcur} gives an approximate induced current $\widehat{I}$ as well as an
approximate total field $\widehat{u}$. Then we seek a regularized solution $\eta$ to the
linearized scattering equation
\begin{equation} \label{scat}
   \int_D G(x,y)\,\widehat{u}(y)\,\eta(y)\,dy=u^s(x),
\end{equation}
by an innovative regularization incorporating both $L^1$ and $H^1$ penalties. The $L^1$
penalty promotes the sparsity of the solution \cite{Tibshirani:1996,
CandesRombergTao:2006, Donoho:2006}. However, the estimate tends to be very spiky with
the $L^1$ penalty used alone. Meanwhile, the conventional $H^1$ penalty can only yield
globally smooth but often overly diffusive profiles. In this work we shall propose a
novel mixed model that consists of a suitable combination of the  $L^1$ and $H^1$
penalties. As we will see, this mixed model produces well clustered and yet distributed
solutions, thereby overcoming the aforementioned drawbacks. It is the mixed model that
enables us to obtain a clear and accurate reconstruction of the inclusions: The
homogeneous background is vividly separated from the scatterers and both support and
intensity of the inclusions are accurately resolved.

Numerically, the $L^1$ penalty term gives rise to a nonsmooth optimality condition, which
renders its direct numerical treatment inconvenient. Fortunately, by using
complementarity functions, the optimality condition reduces to a coupled nonlinear system
for the sought-for coefficient $\eta$ and the Lagrangian multiplier, which is amenable to
efficient numerical solution. We shall develop an efficient and stable semi-smooth Newton
solver for the model via a primal-dual active-set strategy \cite{ItoKunisch:2008}.
Overall, the direct sampling method is very cheap and reduces the computational domain
$D$ in the mixed model (cf. \eqref{scat}) significantly, which in turn makes the
semi-smooth Newton method for the mixed model very fast. Hence, the proposed inverse
scattering method is computationally very efficient.

The rest of the paper is structured as follows. In Section \ref{sec:omega}, we recall a
novel direct sampling method for screening the scatterer support $\Omega$, and derive
thereby an initial guess to the coefficient $\eta$. Then in Section \ref{sec:tikh} we
develop an enhancement technique based on the idea of mixed regularization, and an
efficient semi-smooth Newton solver. Finally, we present numerical results for two- and
three-dimensional examples to demonstrate the accuracy and efficiency of the proposed
inverse scattering method.

\section{A direct sampling method}\label{sec:omega}
In this section, we describe a direct sampling method to determine the shape of the
scatterers, recently derived in \cite{ItoJinZou:2011a}. We only briefly recall the
derivation, but refer the readers to \cite{ItoJinZou:2011a} for more details. Consider a
circular curve $\Gamma$ ($d=2$) or a spherical surface $\Gamma$ ($d=3$). Let $G(x,x_p)$
be the fundamental solution in the homogeneous background. Then using the definitions of
the fundamental solutions $G(x,x_p)$ and $G(x,x_q)$ and Green's second identity, we
deduce
\begin{equation} \label{eqn:green}
  2i\Im(G(x_p,x_q))=\int_{\Gamma}\left[\overline{G}(x,x_q)\frac{\partial G(x,x_p)}{\partial n}
    -G(x,x_p)\frac{\partial \overline{G}(x,x_q)}{\partial n}\right]ds,
\end{equation}
where the points $x_p,x_q\in\Omega_\Gamma$, the domain enclosed by the boundary $\Gamma$.

Next we approximate the right hand side of identity \eqref{eqn:green} by means of the
Sommerfeld radiation condition for the Helmholtz equation, i.e.,
\begin{equation*}
\frac{\partial G(x,x_p)}{\partial n}=ikG(x,x_p) + \mathrm{h.o.t.}.
\end{equation*}
Consequently, we arrive at the following approximate relation
\begin{equation*}
\int_{\Gamma}G(x,x_p)\overline{G}(x,x_q)ds \approx k^{-1}\Im(G(x_p,x_q)),
\end{equation*}
which is valid if the points $x_p$ and $x_q$ are not close to the boundary $\Gamma$.

Now, we consider a sampling domain $\widetilde{\Omega}$ enclosing the scatterer support
$\Omega$. Upon dividing $\widetilde{\Omega}$ into small elements $\{\tau_j\}$, we may
approximate the integral in the scattering relation \eqref{int} by
\begin{equation} \label{eqn:appsc}
   u^s(x)=\int_{\widetilde{\Omega}} G(x,y)I(y)dy\approx \sum_{j}w_j\,G(x,y_j),
\end{equation}
where the weight $w_j$ is given by $|\tau_j|I_j$ with $|\tau_j|$ being the volume of the
$j$th element $\tau_j$. The relation \eqref{eqn:appsc} is plausible if the induced
current $I$ is regular in each element and the elements $\{\tau_j\}$ are sufficiently
fine. It also admits a nice physical interpretation: the scattered field $u^s$ at any
fixed point $x\in\Gamma$ is a weighted average of that due to point scatterers located at
$\{y_j\}$.

Combining the preceding two relations yields
\begin{equation}\label{eqn:appind}
   \int_\Gamma u^s(x)\overline{G}(x,x_p)ds\approx k^{-1}\sum_j w_j\Im(G(y_j,x_p)).
\end{equation}
Hence, if the sampling point $x_p$ is close to some point scatterer $y_j$, i.e.,
$y_j\in\Omega$, then $G(y_j,x_p)$ is nearly singular and takes a very large value,
contributing significantly to the sum in \eqref{eqn:appind}. Conversely, if the point
$x_p$ is far away from all the physical point scatterers,  then the sum will be very
small due to the decay property of $G(y_j,x_p)$.

These facts lead us to the index function  $\Phi(x_p)$ in \eqref{indx} for any $x_p$ in
the sampling region $\widetilde{\Omega}$. In practice, if a point $x_p$ satisfies
$\Phi(x_p)\approx 1$, then it likely lies within the support; whereas if
$\Phi(x_p)\approx 0$, then the point $x_p$ most probably lies outside the support. Hence
it serves as an indicator of the scatterer support. Consequently, we can determine a
domain $D\subset\widetilde{\Omega}$ as one approximate scatterer support, and moreover,
the restriction $\Phi|_D$ of the index $\Phi$ to the domain $D$ may be regarded as a
first approximation to the sought-for coefficient $\eta$. The subdomain $D$ may be chosen
as $D=\{x\in \widetilde{\Omega}: \Phi(x)\ge \mu \max_{x\in \widetilde{\Omega}} \Phi(x)\}$
with $\mu$ being a cut-off value, i.e., the union of elements whose index values are not
less than a specified fraction of the largest index value over the sampling region
$\widetilde{\Omega}$. This determination of subdomain $D$ will be adopted in our
numerical experiments.

This method is of sampling type (cf. \cite{Potthast:2006} for an overview of existing
sampling methods), and its flavor closely resembles \textit{multiple signal
classification} \cite{Schmidt:1986,Cheney:2001,Devaney:2006} and \textit{linear sampling
method} \cite{ColtonKirsch:1996,Kirsch:1998}. However, unlike these existing techniques,
it works with a few (e.g., one or two) incident waves, is highly tolerant with respect to
noise, and involves only computing inner products with fundamental solutions rather than
expensive matrix operations as in the other two techniques. The robustness of $\Phi$ is
attributed to the fact that the (high-frequency) noise is roughly orthogonal to the
(smooth) fundamental solutions.

\section{Mixed regularization}\label{sec:tikh}

The direct sampling method in Section \ref{sec:omega} extracts an accurate estimate $D$
to the scatterer support $\Omega$ as well as a reasonable initial guess to the medium
coefficient $\eta$, i.e., $\Phi|_D$. In this part we refine the approximation $\Phi|_D$
by exploiting the idea of nonsmooth mixed regularization. Given the approximation
$\Phi|_D$, we can compute the induced current $\widehat{I}$ via \eqref{eqn:indcur} for
each incident wave and the corresponding total field $\widehat{u}$ on the domain $D$ from
\eqref{int}. By substituting the approximation $\widehat{u}$ into equation \eqref{int},
we arrive at the following linearized problem
\begin{equation*}
   \int_DG(x,y)\widehat{u}(y)\eta(y)dy = u^s(x)\,, \quad x\in\Gamma.
\end{equation*}
It is convenient to introduce a linear integral operator $K:L^2(D)\mapsto L^2(\Gamma)$
defined by
\begin{equation}\label{eqn:linint}
  (K\eta)(x)=\int_DG(x,y)\widehat{u}(y)\,\eta(y)\,dy.
\end{equation}

We observe that the kernel $G(x,y)\widehat{u}(y)$ is smooth due to the analyticity of
fundamental functions $G(x,y)$ away from the singularities and standard Sobolev
smoothness of the total field $\widehat{u}(y)$ (following from elliptic regularity theory
\cite{ColtonKress:1998}). Hence, the linear operator $K:L^2(D)\mapsto L^2(\Gamma)$ is
compact. As a consequence, the linearized problem \eqref{eqn:linint} is ill-posed in the
sense that small perturbations in the data can lead to huge deviations in the solution,
according to the classical inverse theory \cite{TikhonovArsenin:1977}, which is
reminiscent of the severe ill-posedness of the IMSP, and its stable and accurate
numerical solution calls for regularization techniques.

We determine an enhanced estimate of the coefficient $\eta$ from the linearized problem
\eqref{eqn:linint} by solving the following variational problem:
\begin{equation} \label{mixed}
  \min\quad \frac{1}{2}\int_\Gamma|K\eta -u^s|^2ds +
    \alpha\int_D|\eta|dx+\frac{\beta}{2}\int_D|\nabla \eta|^2dx.
\end{equation}
In comparison with more conventional regularization techniques, the most salient feature
of the model \eqref{mixed} lies in two penalty terms: it contains both the $L^1$ penalty
and the $H^1$ penalty, which exert drastically different a priori knowledge of the
sought-for solution. The scalars $\alpha$ and $\beta$ are regularization parameters
controlling the strength of respective regularization.

This variational problem \eqref{mixed} allows us to determine a coefficient $\eta$ which
is distributed yet clustered, i.e., exhibiting a clear groupwise sparsity structure in
the canonical pixel basis. This a priori knowledge is plausible for localized
inclusions/inhomogeneities in a homogeneous background. The model \eqref{mixed} is
derived from the following widely accepted observations. The $L^1$ penalty promotes the
sparsity of the solution \cite{Tibshirani:1996,CandesRombergTao:2006,Donoho:2006}, i.e.,
the solution is very much localized. Hence, the estimated background is homogeneous.
However, if the $L^1$ penalty is used alone, the solution tends to be very spiky and may
miss numerous physically relevant pixels in the sought-for groups. That is, the desirable
groupwise structure is not preserved. Meanwhile, the more conventional $H^1$ penalty
\cite{TikhonovArsenin:1977} yields a globally smooth profile, but the solution is often
overly diffusive. Consequently, the overall structure stands out clearly, but the
retrieved background is very blurry, which may lead to erroneous diagnosis of the number
of the inclusions and their sizes. In order to preserve simultaneously these distinct
features of the sought-for coefficient, i.e., sparsely distributed groupwise structures
in a homogeneous background, a natural idea would be to combine the $L^1$ penalty with
the $H^1$ penalty, in the hope of retaining the strengths of both models. As we shall see
below, the idea does work very well, and the model is very effective for enhancing the
resolution of the estimate to the coefficient $\eta$.

The general idea of mixed regularization, i.e., using multiple penalties, has proved very
effective in promoting several distinct features simultaneously. This general idea has
been pursued in the imaging community \cite{ItoKunisch:2000,LuShenXu:2007}. However, to
the best of our knowledge, the model \eqref{mixed} has not been explored in the
literature, let alone its efficient and accurate numerical treatment. A detailed
mathematical analysis of the model \eqref{mixed} is beyond the scope of the present
paper. We refer interested readers to \cite{ItoJinTakeuchi:2011} for some preliminary
results on mixed regularization and to \cite{JinLorenzSchiffler:2009} for a related model
(elastic-net).

To fully explore the potentials of the model \eqref{mixed}, an efficient and accurate
solver is required. We shall develop a semi-smooth Newton type method, which allows extracting very
detailed features of the solutions to the model \eqref{mixed}. The starting point of the
algorithm is the necessary optimality condition of the variational problem \eqref{mixed},
which reads
\begin{equation} \label{fixed}
   K^*K\eta-\beta\Delta\eta-K^\ast u^s \in -\alpha\partial \psi(\eta),
\end{equation}
where $\psi(\eta)=\|\eta\|_{L^1}$ and the subdifferential $\partial\psi(\eta)$
\cite{ItoKunisch:2008} is the set-valued signum function, which is defined pointwise as
\begin{equation*}
\partial\psi(\eta)(x)=\left\{\begin{array}{ll}
       1, &\mbox{if  }\eta(x)>0, \\
       \,[-1,1\,], &\mbox{if  }\eta(x)=0,\\
       -1, &\mbox{if  }\eta(x)<0.
       \end{array}\right.
\end{equation*}

Due to the convexity of the functional, the relation \eqref{fixed} is also a sufficient
condition. Hence it suffices to solve the inclusion \eqref{fixed}, for which there are
several different ways, e.g., iterative soft shrinkage
\cite{DaubechiesDefrise:2004,WrightNowakFigueiredo:2009}, augmented Lagrangian
method/alternating direction method or semi-smooth Newton method \cite{ItoKunisch:2008}.
We shall develop a (new) semi-smooth Newton method to efficiently solve the inclusion
\eqref{fixed}. To this end, we first recall the complementarity condition
\cite{ItoKunisch:2008}
\begin{equation} \label{eqn:complrel}
   \lambda=\frac{\lambda+c\eta}{\max(1,|\lambda+c\,\eta|)}
\end{equation}
for any $c>0$, which will be fixed at a constant in the implementation, and $\lambda$
serves as a Lagrange multiplier. It can be directly verified by pointwise inspection that
the complementarity condition \eqref{eqn:complrel} is equivalent to the inclusion
$\lambda\in\partial\psi(\eta)$ (cf.\,\cite{ItoKunisch:2008}). With the help of the
complementarity condition \eqref{eqn:complrel}, we arrive at the following equivalent
nonlinear system in the primal variable $\eta$ and dual variable $\lambda$:
\begin{equation*}
  \left\{\begin{aligned}
      K^*K\eta-\beta\Delta\eta-K^\ast u^s +\alpha\lambda &= 0,\\
      \lambda-\frac{\lambda+c\eta}{\max(1,|\lambda+c\,\eta|)}&=0.
  \end{aligned}\right.
\end{equation*}

Then we apply the semi-smooth Newton algorithm using a primal-dual active set strategy.
The complete implementation is listed in Algorithm \ref{alg:ssn}. The technical details
for deriving the crucial Newton step (Step 5) are deferred to Appendix \ref{app:ssn},
which involves damping and regularization. A natural choice of the stopping criterion at
Step 6 is based on monitoring the change of the active set $\mathcal{A}=\{x\in
D:|\lambda+c\,\eta|\le 1\}$: if the active sets for two consecutive iterations coincide,
then we can terminate the algorithm, cf. \cite{ItoKunisch:2008}.

\begin{algorithm}[h] \caption{Primal-dual active set method}\label{alg:ssn}
   \begin{algorithmic}[1]
      \STATE Initialize $\eta^0$ and $\lambda^0$, and set $c > 0$ and $k=0$.
      \FOR {$k=0,\dots,K$}
          \STATE Set the active set ${\cal A}^k$ and inactive set ${\cal I}^k$ respectively by
             \begin{equation*}
                \begin{aligned}
                   {\cal A}^k&=\{x\in D:|\lambda^k+c\,\eta^k|\le 1\},\\
                   {\cal I}^k&=\{x\in D: |\lambda^k+c\,\eta^k|>1\}.
                \end{aligned}
             \end{equation*}
          \STATE Compute $a$ and $b$ by
             \begin{equation*}
                  a=\frac{\lambda^k}{\max(|\lambda^k|,1)} \quad\mbox{and}\quad b=\frac{\lambda^k+c\,\eta^k}{|\lambda^k+c\,\eta^k|},
             \end{equation*}
             and set $d^k=|\lambda^k+c\eta^k|$ and $F^k=ab^t$.
         \STATE Solve for $(\eta^{k+1},\lambda^{k+1})$ from the system
            \begin{equation*}
               \begin{aligned}
                  K^\ast K\eta^{k+1}-\beta\,\Delta\eta^{k+1}-K^\ast u^s+\alpha\,\lambda^{k+1}&=0\mbox{ on } {\cal I}^k,\\
                  \lambda^{k+1}-c\frac{1}{d^k-1}(I-F^k) \eta^{k+1}-\frac{\lambda^k}{\max(|\lambda^k|,1)}&=0,\\
                  \eta^{k+1}&=0 \mbox{ on }{\cal A}^k.
               \end{aligned}
            \end{equation*}
         \STATE Check the stopping criterion.
      \ENDFOR
      \STATE {\textbf{output}} approximation $\eta^K$.
   \end{algorithmic}
\end{algorithm}

The main computational effort of the algorithm lies in the Newton update at Step 5: each
iteration requires solving a (dense) linear system. We note that the dual variable
$\lambda$ can be expressed in terms of the primal variable $\eta$ and on the active set
$\mathcal{A}$, the coefficient $\eta$ vanishes identically. Thus in practice, we solve
only a linear system for $\eta$ on the inactive set $\mathcal{I}=D\setminus\mathcal{A}$.
An important feature of the algorithm is that the linear system becomes smaller and
smaller and also less and less ill-conditioned as the iteration goes on, while the
iterate captures more and more refined details of the nonhomogeneous medium regions. If
the exact solution is indeed sparse (many zero entries), then the system size, i.e.,
$|\mathcal{I}|$, usually shrinks quickly as the iteration proceeds. The numerical
experiments indicate that the convergence of the algorithm is rather steady and fast.

\section{Numerical experiments}

In this part, we present numerical results for several two- and three-dimensional
examples to showcase the proposed two-stage inverse scattering method, for both exact and
noisy data. The wave number $k$ is fixed at $2\pi$, and the wavelength is set to
$\lambda=1$. The exact scattered field $u^s$ is obtained by first solving the integral
equation \eqref{eqn:indcur} for the induced current $I$ and then substituting the current
$I$ into the integral representation \eqref{int}. Here the integral equation
\eqref{eqn:indcur} is discretized by a mid-point rule; see Appendix \ref{app:int} for
details. The noisy scattered data $u_\delta^s$ are generated pointwise by the formula
\begin{equation*}
   u_\delta^s(x) = u^s(x) + \epsilon \zeta \max_{x\in\Gamma}|u^s(x)|,
\end{equation*}
where $\epsilon$ refers to the relative noise level, and both the real and imaginary
parts of the noise $\zeta=\zeta(x)$ follow the standard Gaussian distribution. The index
$\Phi$, its restriction $\Phi|_D$ and the enhanced approximation $\eta$ by the mixed
model will be displayed. As is mentioned in Section 2, we choose the subdomain $D$
(approximate scatterer support) based on the formula $D=\{x\in \widetilde{\Omega}:
\Phi(x)\ge \mu \max_{x\in \widetilde{\Omega}} \Phi(x)\}$, where the cut-off value $\mu$
is taken in the range $(0.5,0.7)$. The choice of the cutoff value $\mu$ affects directly
the size of the domain $D$, but does not cause much effects on the reconstructions.

Like in any regularization technique, an appropriate choice of regularization parameters
$(\alpha,\beta)$ in the mixed model \eqref{mixed} is crucial for the success of the
proposed imaging algorithm. There have been a number of choice rules
\cite{ItoJinTakeuchi:2011b} for one single parameter, but very little is known about the
mixed model. We shall choose the pair $(\alpha,\beta)$ in a trial and error manner, which
suffices our goal of illustrating the significant potentials of the mixed model for
inverse scattering. In Algorithm \ref{alg:ssn}, the parameter $c$ is set to $50$, and
both $\eta^0$ and $\lambda^0$ are initialized to zero. The maximum number $K$ of Newton
iterations is $50$. In all the experiments, the convergence of the algorithm is achieved
within about $10$ iterations. All the computations were performed on \textsc{MATLAB}
7.12.0 (R2011a) on a dual-core desktop computer with 2GB RAM.

\subsection{Two-dimensional examples}
Unless otherwise specified, one incident direction $d$ is employed for two-dimensional
problems, and it is fixed at $\frac{1}{\sqrt{2}} (1,1)^\mathrm{T}$. The scattered field
$u^s$ is measured at $30$ points uniformly distributed on a circle of radius $5$. The
sampling domain $\widetilde{\Omega}$ is fixed at $[-2,2]^2$, which is divided into a
uniform mesh consisting of small squares of width $h=0.01$. The subdomain $D$ for the
integral equation \eqref{eqn:linint} is divided into a coarser uniform mesh consisting of
small squares of width $0.02$.

\begin{figure}
  \centering
  \begin{tabular}{cccc}
     \includegraphics[width=.25\textwidth]{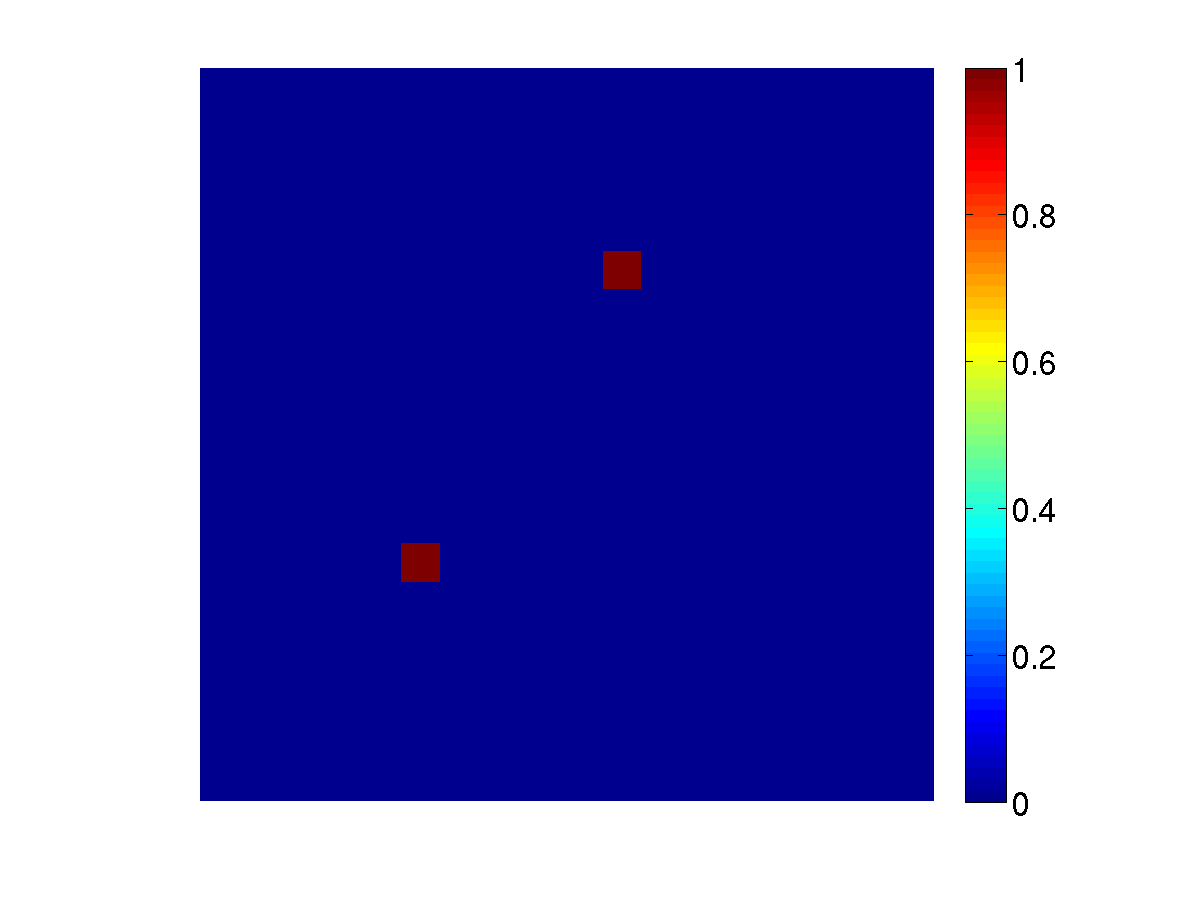} & \includegraphics[width=.25\textwidth]{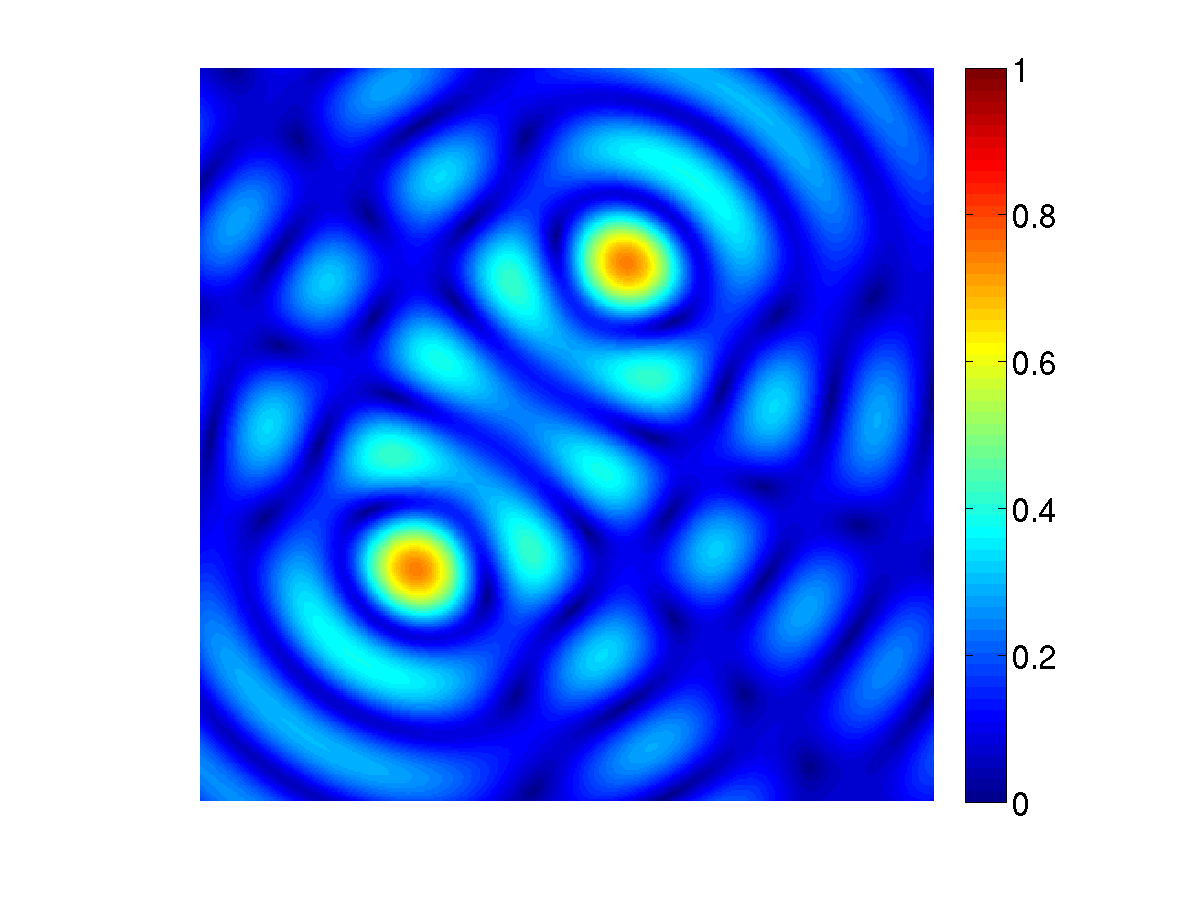}
    &\includegraphics[width=.25\textwidth]{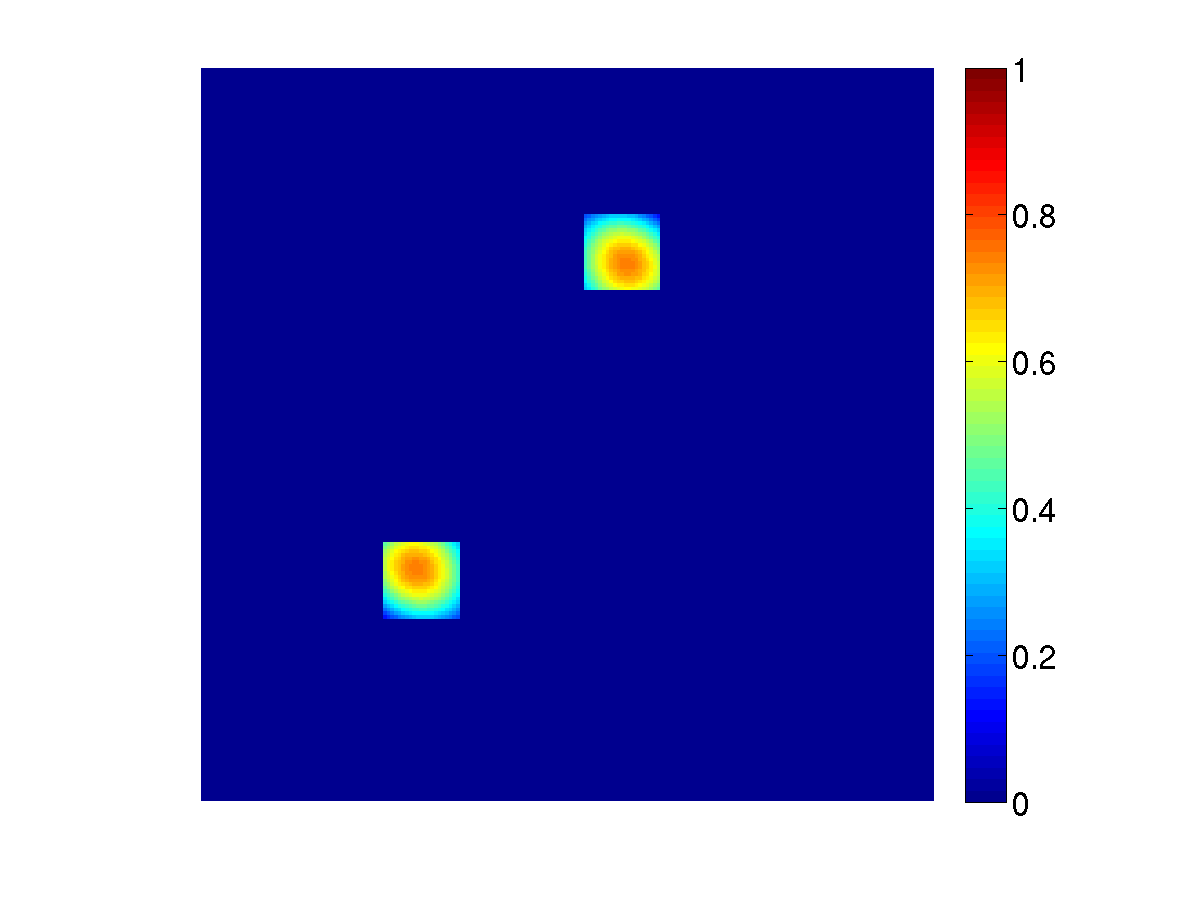} & \includegraphics[width=.25\textwidth]{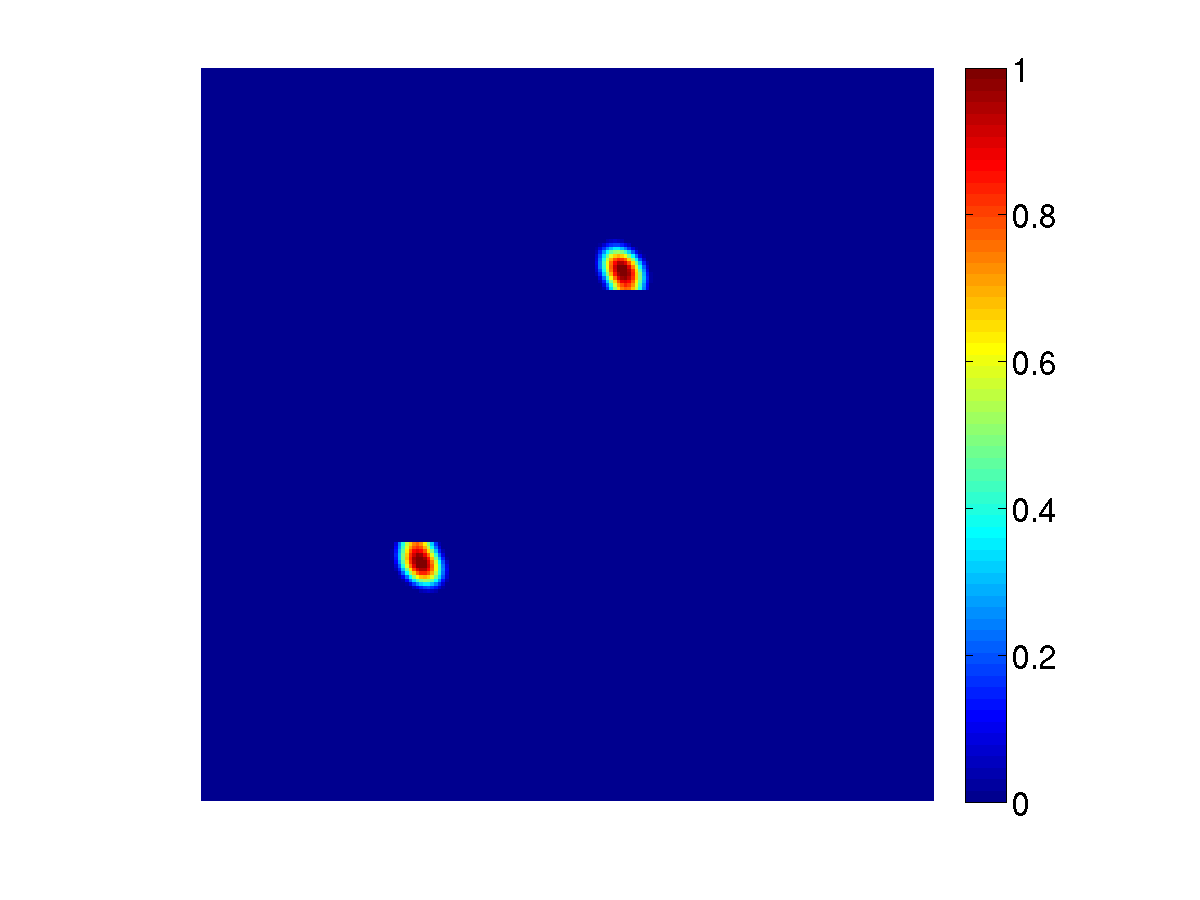}\\
     \includegraphics[width=.25\textwidth]{exam3_recon_true.png} & \includegraphics[width=.25\textwidth]{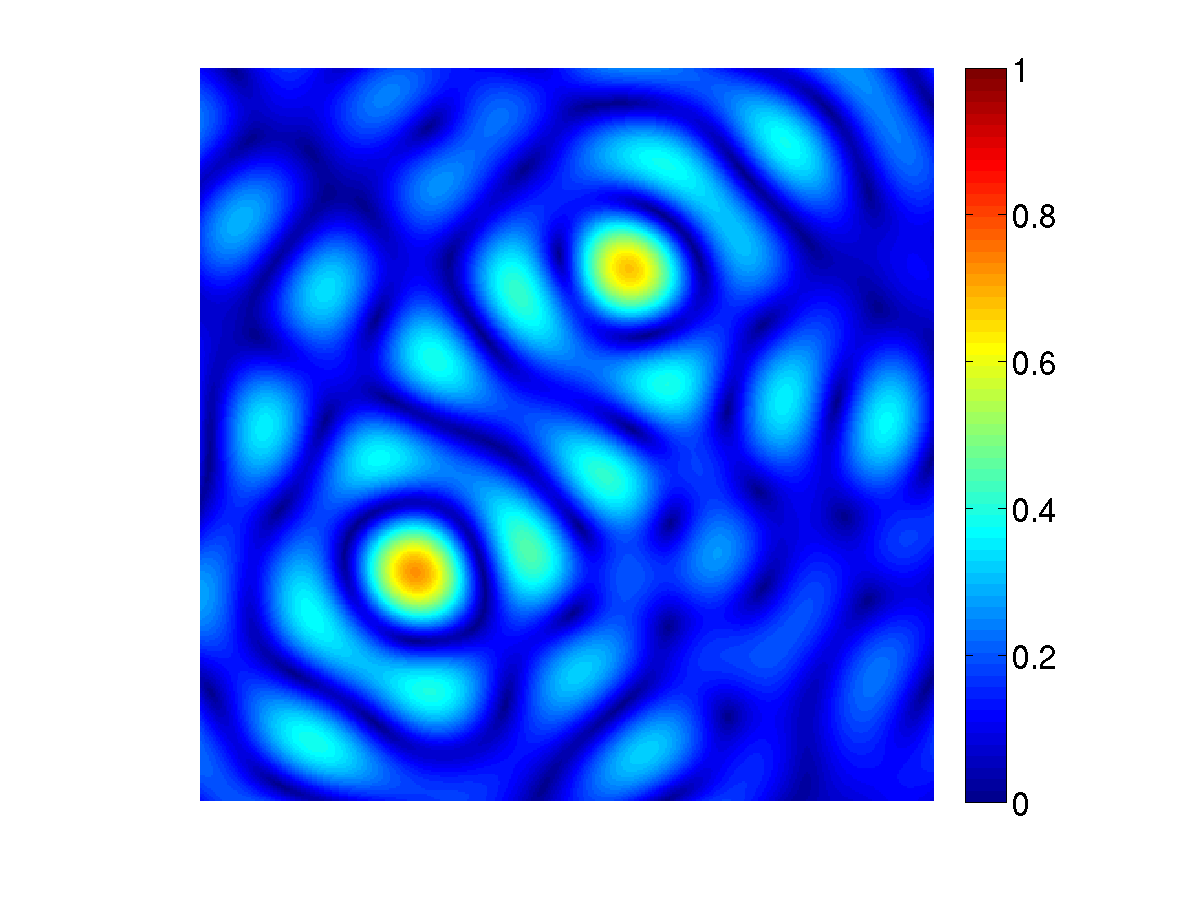}
    &\includegraphics[width=.25\textwidth]{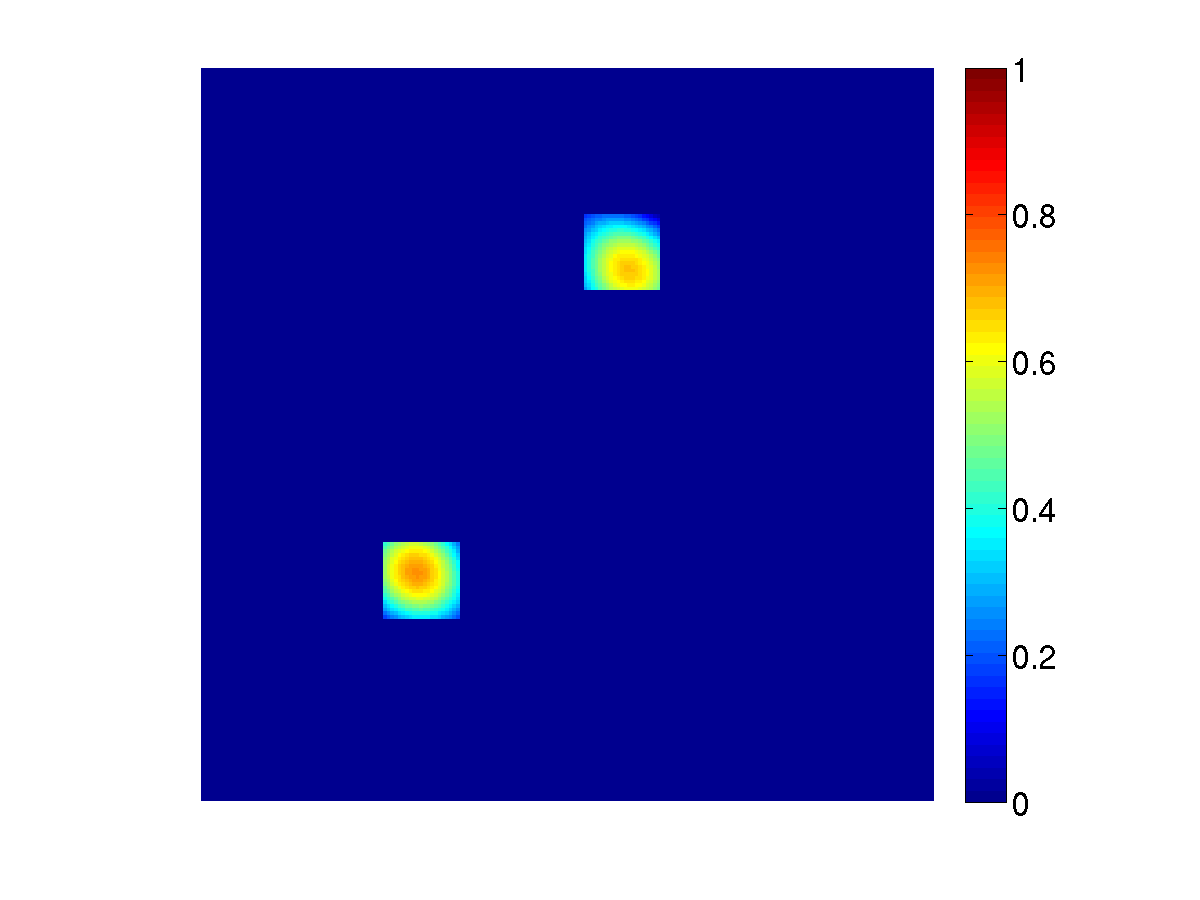}& \includegraphics[width=.25\textwidth]{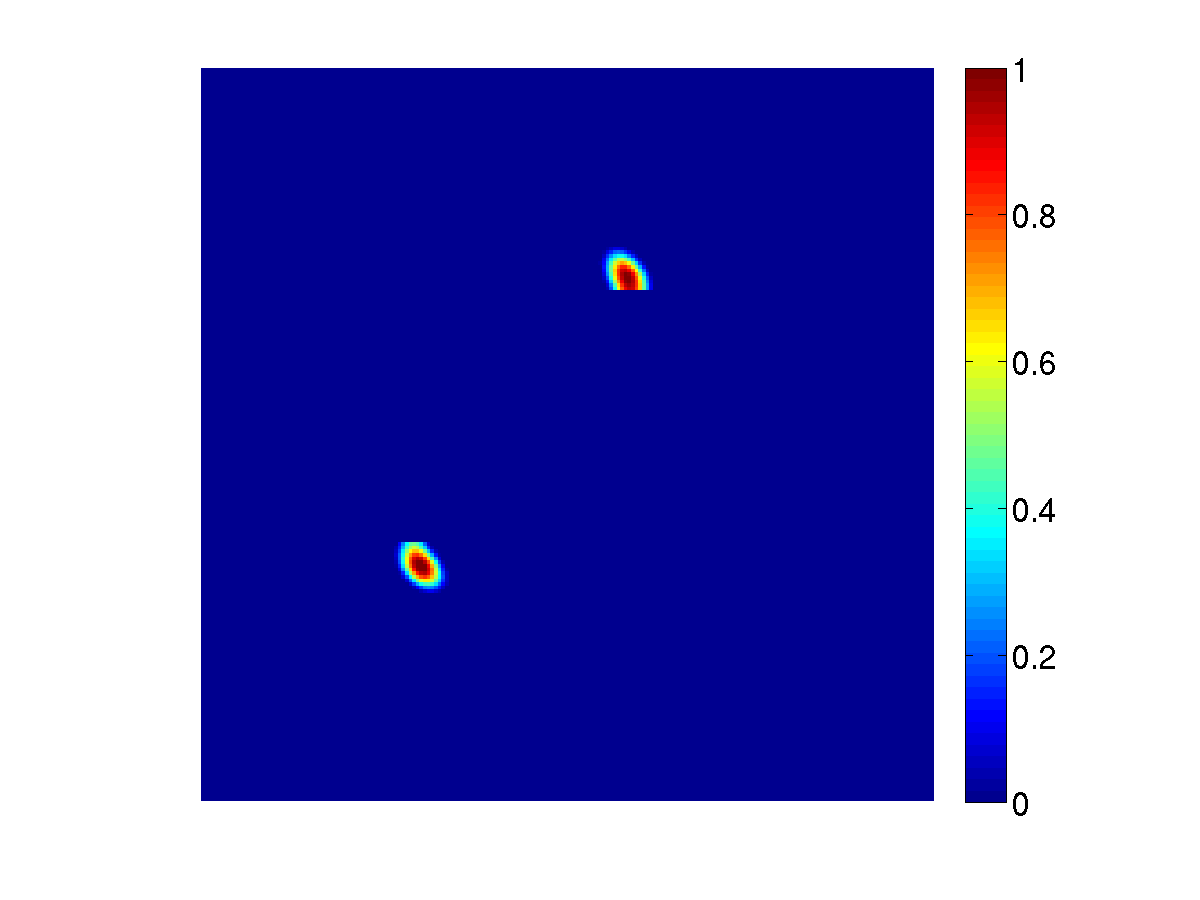}\\
    (a) true scatterer & (b) index $\Phi$ & (c) index $\Phi|_D$ & (d) sparse recon.
  \end{tabular}
  \caption{Numerical results for Example \ref{exam:2sc}(a): (a) true scatterer, (b) index
    $\Phi$, (c) index $\Phi|_D$ (restriction to the subdomain $D$) and (d) sparse reconstruction.
    The first and second rows refer to exact data and the data with $20\%$ noise, respectively.}
  \label{fig:2sca}
\end{figure}

Our first example illustrates the method for two separate scatterers.
\begin{exam}\label{exam:2sc}
We consider two separate square scatterers in the following two scenarios
\begin{itemize}
 \item[(a)] The scatterers are of width $0.2$ and centered at $(-0.8,-0.7)$ and
     $(0.3,0.9)$, respectively, and the coefficient $\eta$ in both region is $1$.
 \item[(b)] The scatterers are of width $0.3$ and centered at $(-0.25,0)$ and
     $(0.25,0)$, respectively, and the coefficient $\eta$ in the former and latter is
     $1.5$ and $1$, respectively.
\end{itemize}
\end{exam}

The two scatterers in Example \ref{exam:2sc}(a) are well apart from each other. The
recovery of the scatterer locations by the index $\Phi$ is quite satisfactory, especially
upon noting that we have just used one incident wave. Two distinct scatterers are
observed for both exact data and the data with $20\%$ noise, cf. Fig. \ref{fig:2sca}(b).
However, the magnitudes are inaccurate, and the estimate suffers from spurious
oscillations in the homogeneous background, due to the ill-posed nature of the IMSP and
the oscillating behavior of fundamental solutions. Nonetheless, two localized square
subregions $D$ (each of width $0.4$) encompass the modes of the index $\Phi$, see Fig.
\ref{fig:2sca}(c), and may be taken as an approximate scatterer support. {In Fig.
\ref{fig:2sca}(c), the entire sampling domain $\widetilde{\Omega}$ is shown, and the two
small squares represent the approximate support $D$. Outside of the domain $D$, the index
$\Phi|_D$ is set to zero, i.e., identical with homogeneous background, and will not be
updated during the enhancing step.} The enhancing step is initialized with $\Phi|_D$, and
the results are shown in Fig. \ref{fig:2sca}(d). The regularization parameters for
getting the reconstructions, which are determined in a trial-and-error manner, are
presented in Table \ref{tab:regpara}. The enhancement {of the approximation $\Phi_D$ over
the domain $D$} is significant: the recovered background is now mostly homogeneous, and
the magnitudes and sizes of the recovered scatterers agree well with the exact ones. This
shows clearly the significant potentials of the proposed mixed regularization for inverse
scattering problems. {The numbers in Table \ref{tab:regpara} also sheds valuable insights
into the mixed model \eqref{mixed}: the value of the regularization parameter $\alpha$ is
much larger than that of $\beta$. Hence, the $L^1$ penalty plays a predominant role in
ensuring the sparsity of the solution, whereas the $H^1$ penalty yields a locally smooth
structure.}

\begin{table}
  \centering
  \caption{Regularization parameters $(\alpha,\beta)$ for the examples.}
  \begin{tabular}{ccccc}
  \hline
    example & \ref{exam:2sc}(a) & \ref{exam:2sc}(b) & \ref{exam:ring} & \ref{exam:cube}\\
  \hline
  $\epsilon = 0\%$ &(2.0e-6,1.5e-9) & (8.0e-6,1.4e-8) & (7.0e-6,1.0e-9) &(2.5e-9,4.0e-14)\\
  $\epsilon = 20\%$&(3.0e-6,2.0e-9) & (8.5e-6,9.0e-9) & (7.0e-6,5.0e-9) &(2.5e-9,5.0e-14)\\
  \hline
  \end{tabular}  \label{tab:regpara}
\end{table}

\begin{figure}[h!]
  \centering
  \begin{tabular}{cccc}
    \includegraphics[width=.25\textwidth]{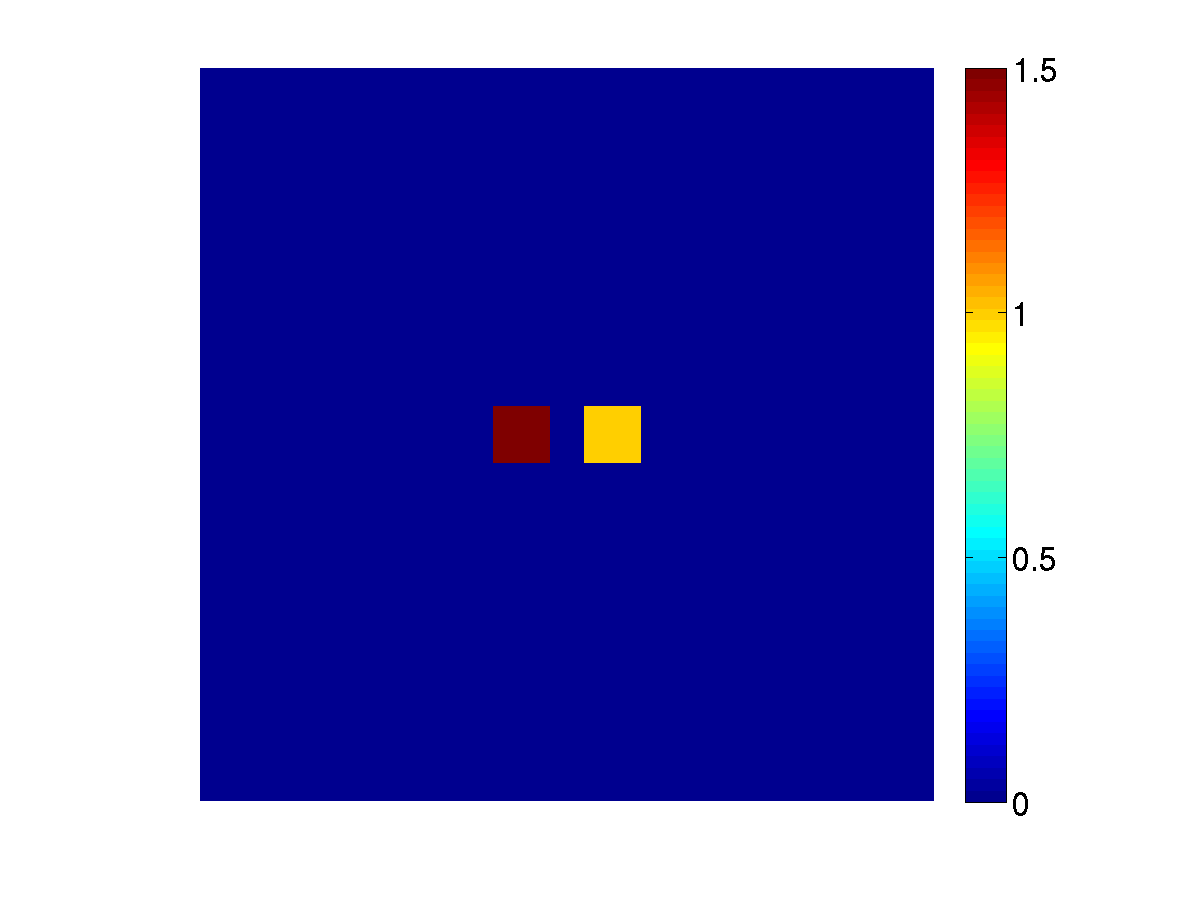} & \includegraphics[width=.25\textwidth]{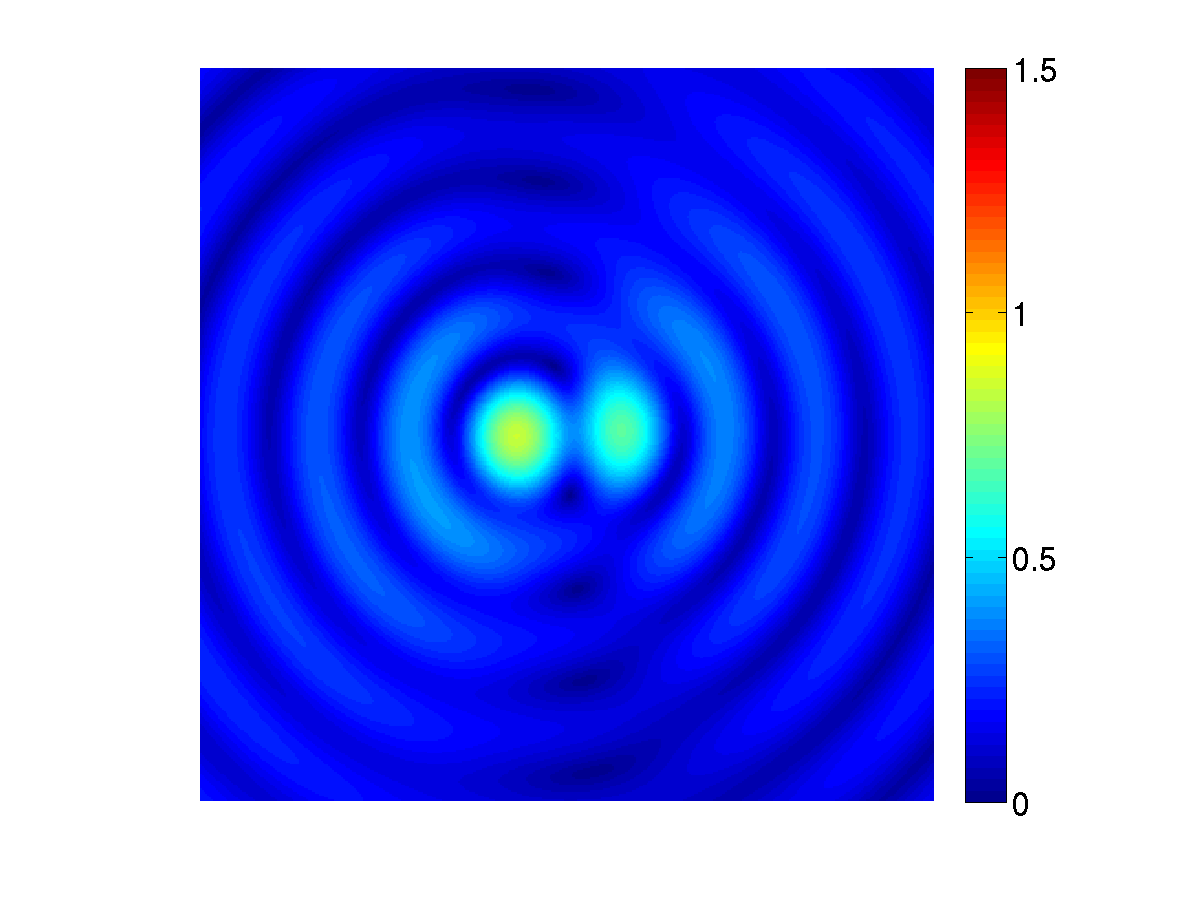}
   &\includegraphics[width=.25\textwidth]{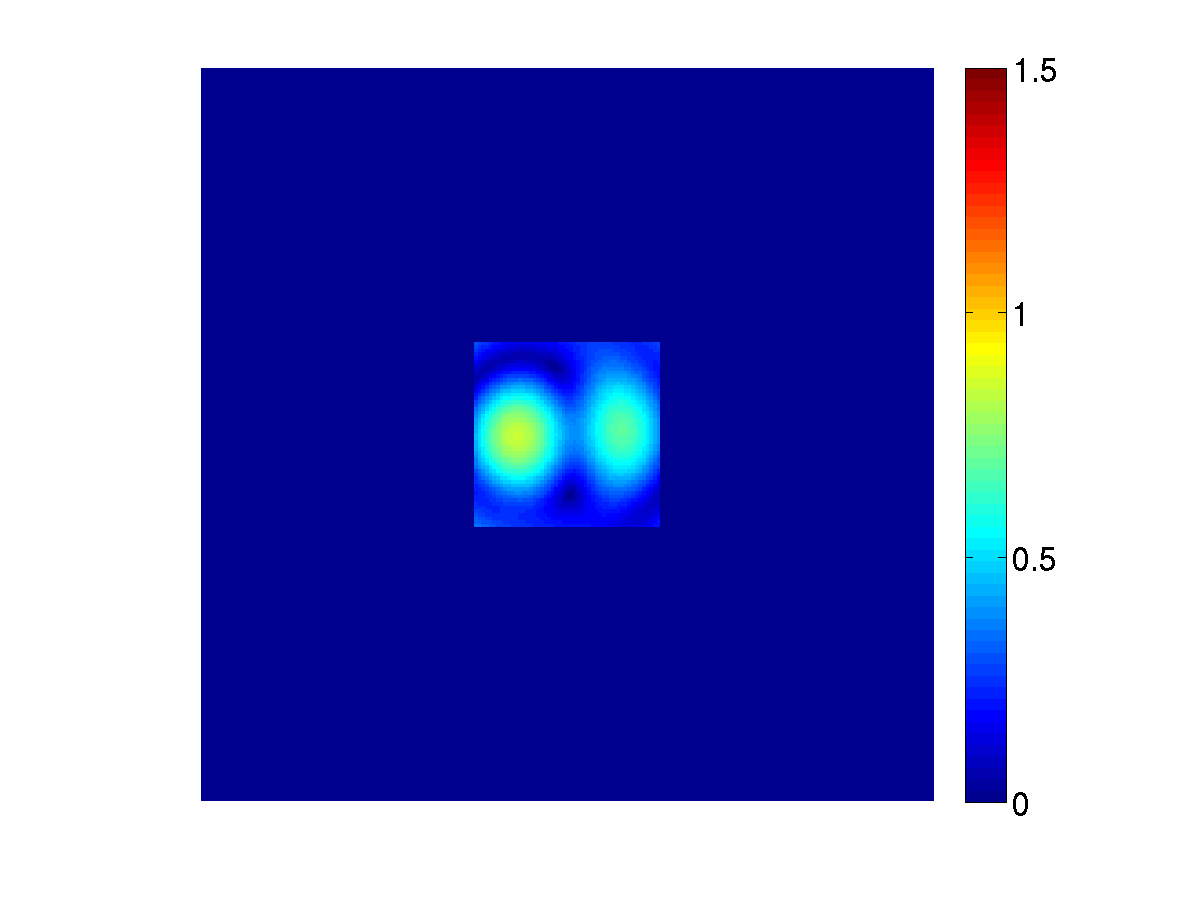} & \includegraphics[width=.25\textwidth]{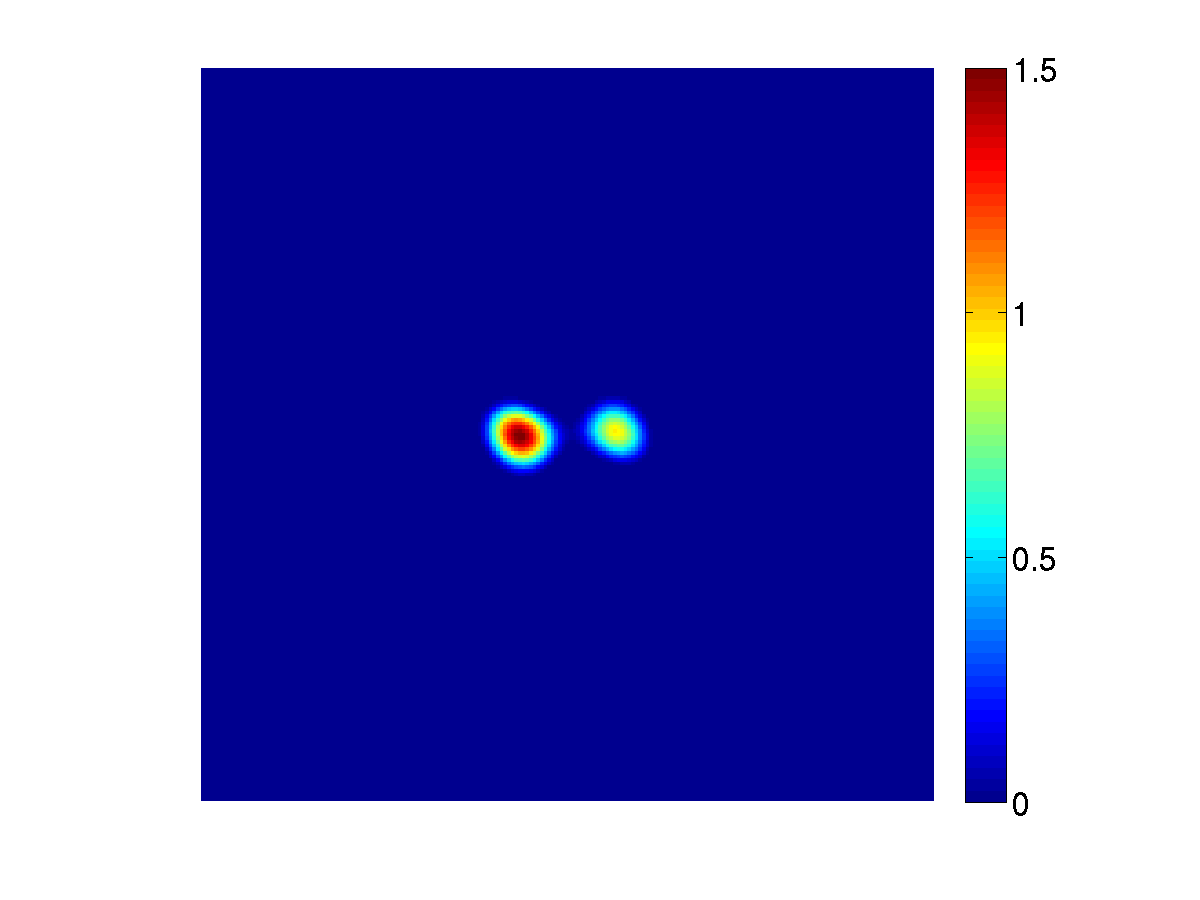}\\
    \includegraphics[width=.25\textwidth]{exam1_recon_true.png} & \includegraphics[width=.25\textwidth]{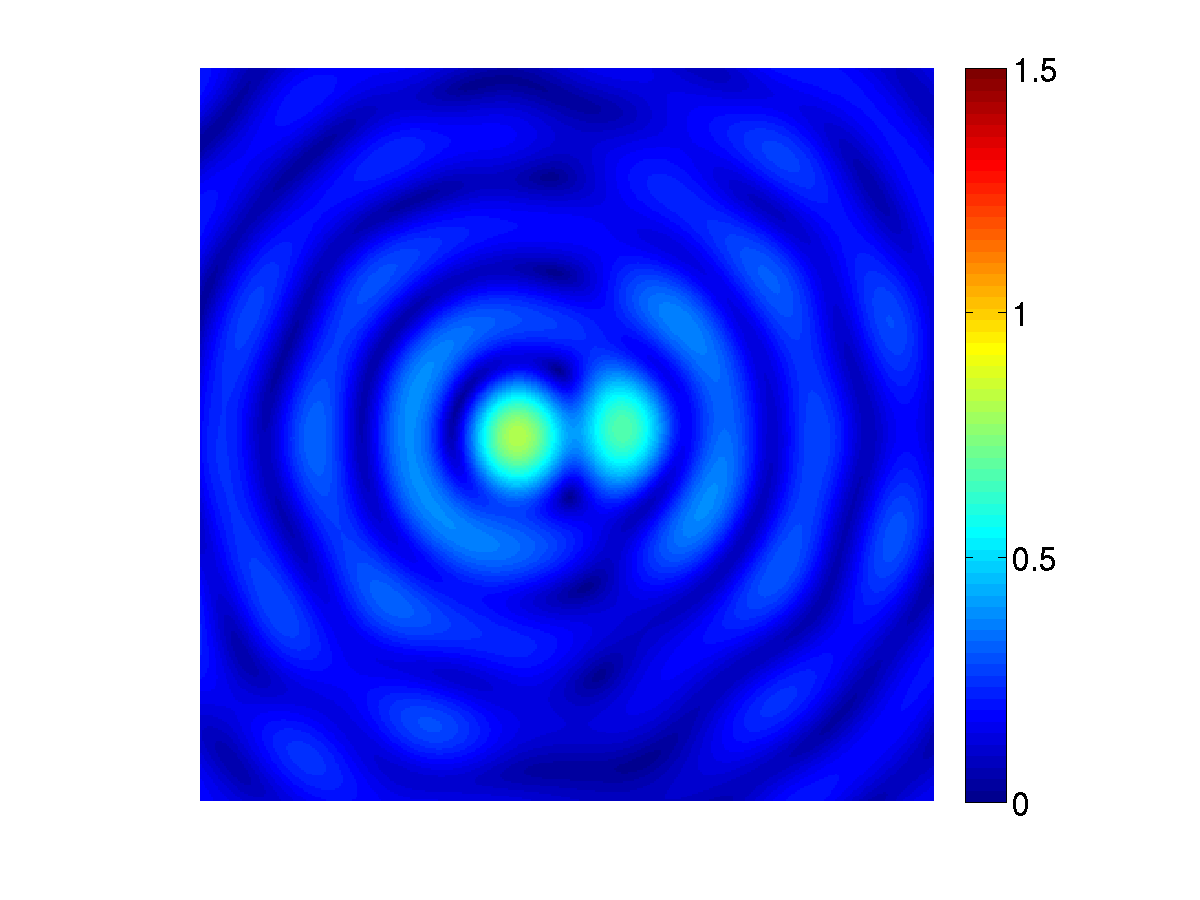}
   &\includegraphics[width=.25\textwidth]{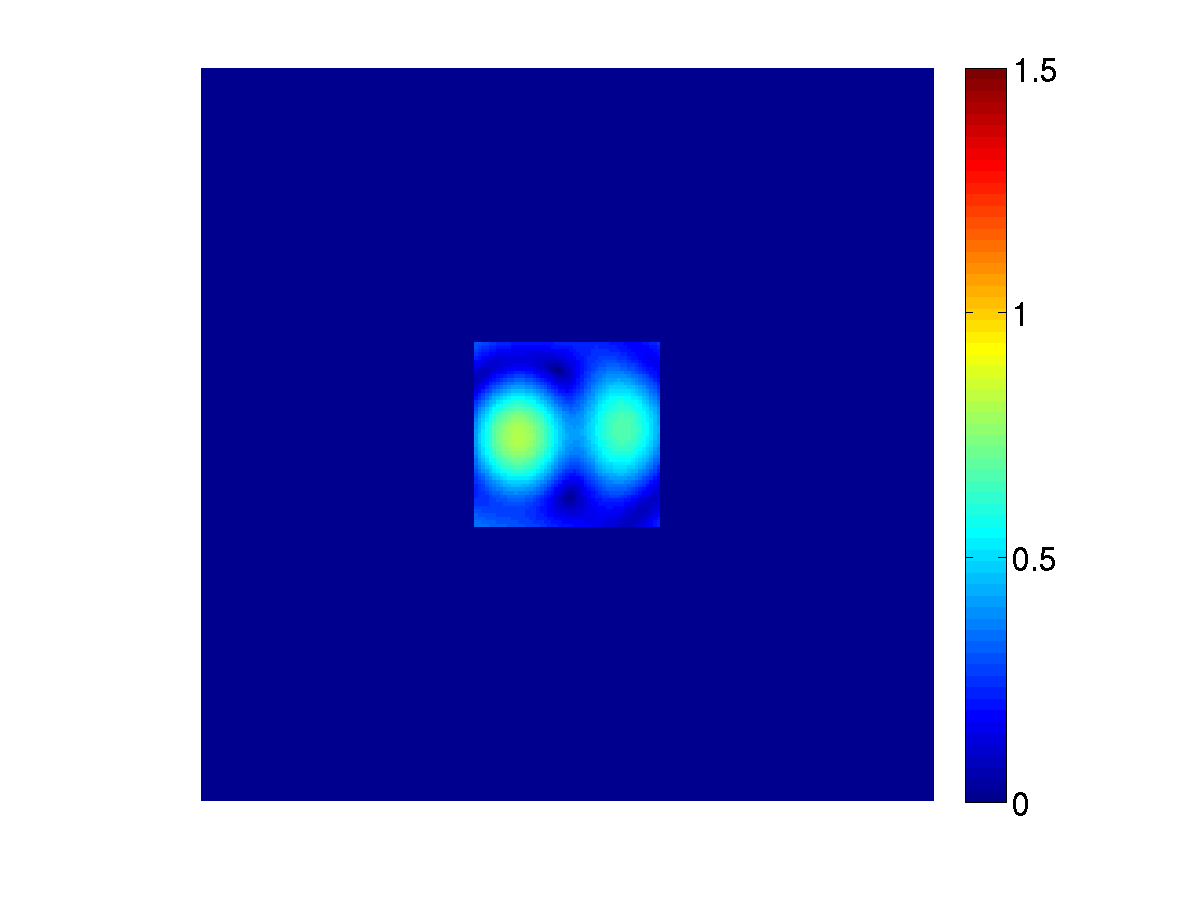}& \includegraphics[width=.25\textwidth]{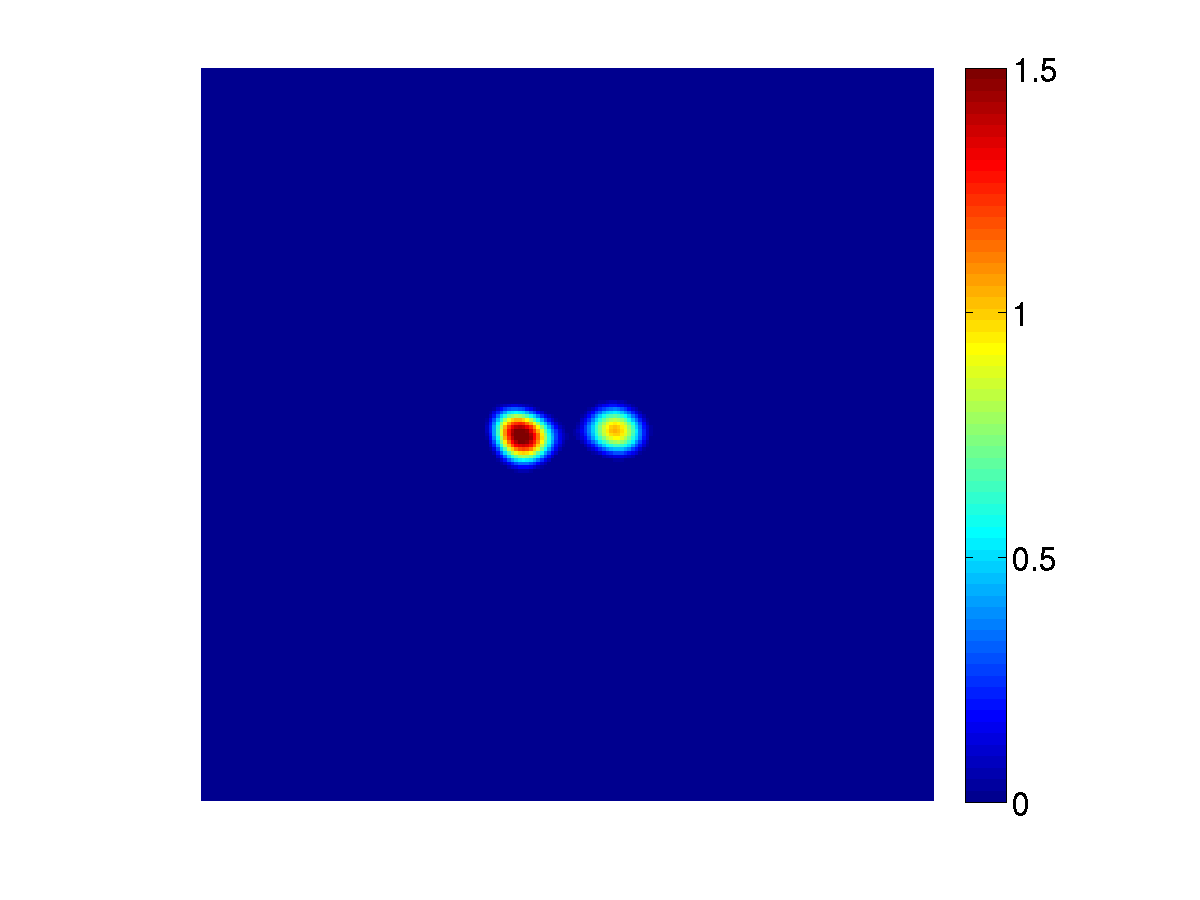}\\
    \includegraphics[width=.25\textwidth]{exam1_recon_true.png} & \includegraphics[width=.25\textwidth]{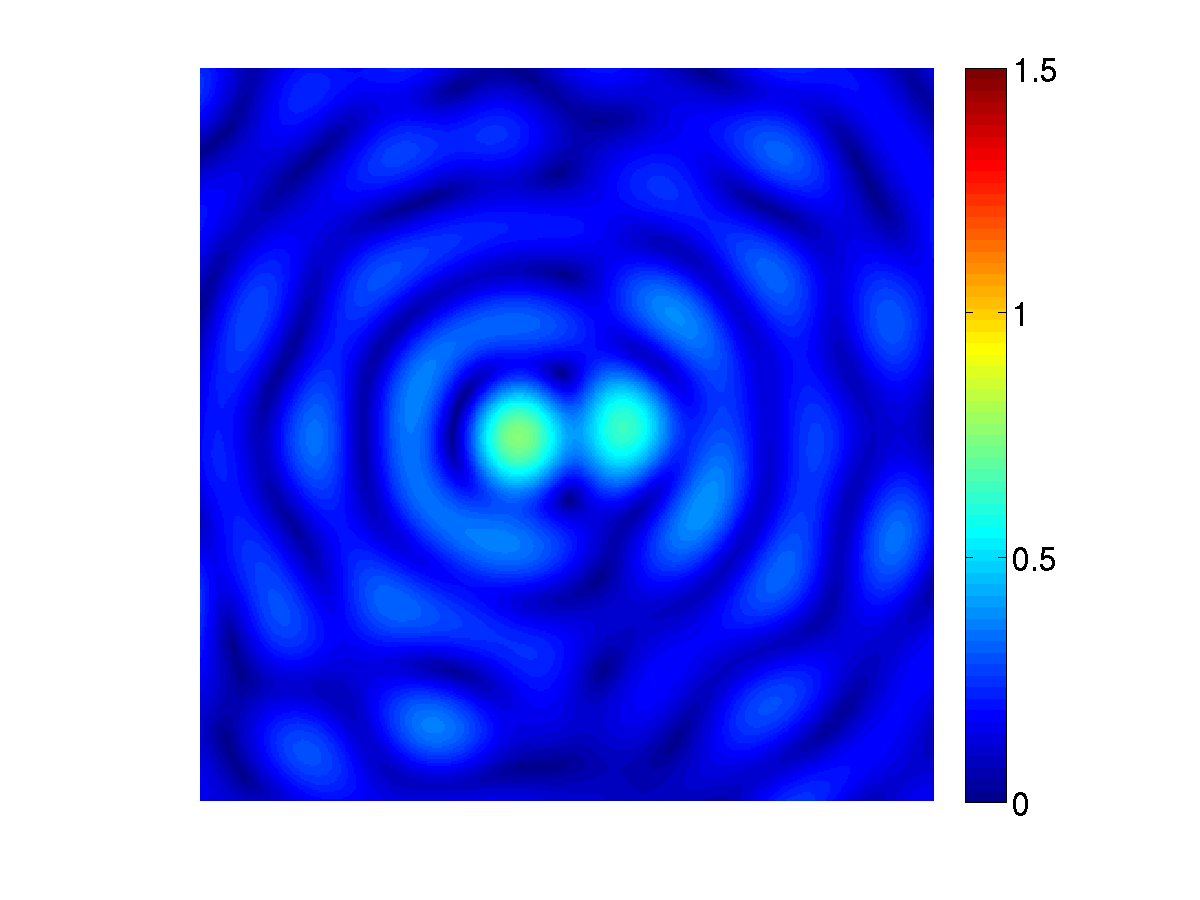}
   &\includegraphics[width=.25\textwidth]{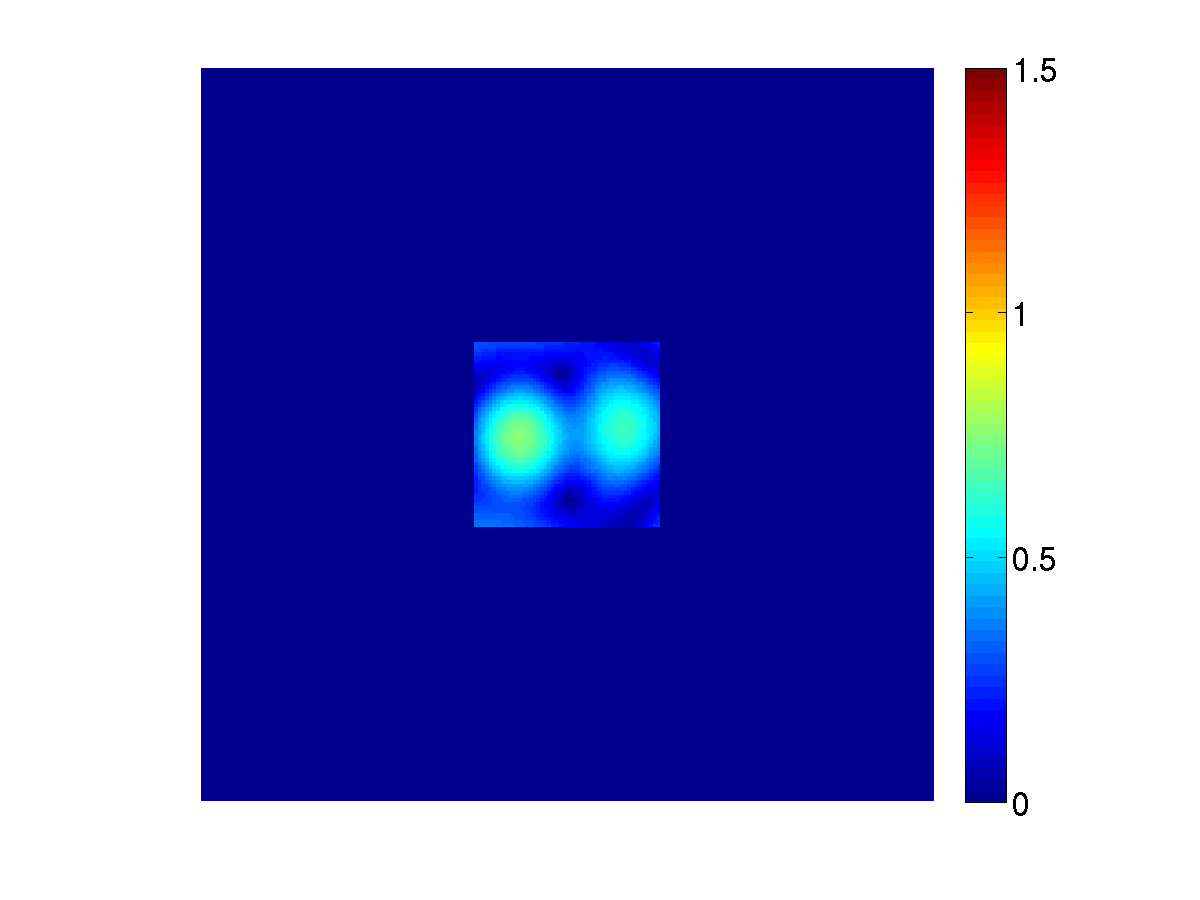}& \includegraphics[width=.25\textwidth]{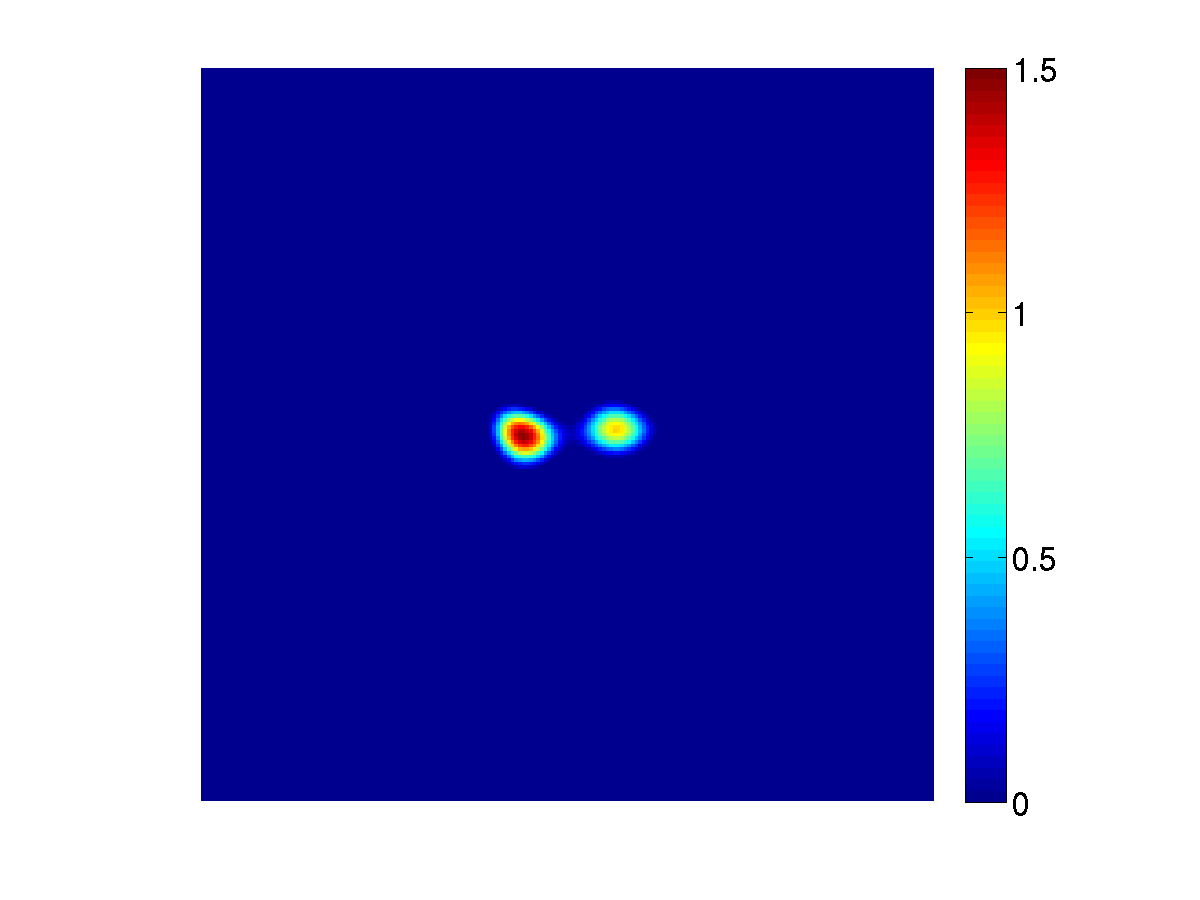}\\
   (a) true scatterer & (b) index $\Phi$ & (c) index $\Phi|_D$ & (d) sparse recon.
  \end{tabular}
  \caption{Numerical results for Example \ref{exam:2sc}(b): (a) true scatterer, (b) index
   $\Phi$, (c) index $\Phi|_D$ (restriction to the subdomain $D$) and (d) sparse reconstruction. From top to bottom,
   the rows refer to exact data, noisy data with $10\%$ and $20\%$, respectively.}
  \label{fig:2scb}
\end{figure}

The two scatterers in Example \ref{exam:2sc}(b) stay very close to each other, and thus
it is rather challenging for precise numerical reconstruction. The detection of the
scatterer locations by the index $\Phi$, see Fig. \ref{fig:2scb}(b), is still very
impressive. In particular, it clearly distinguishes two separate scatterers with their
locations correctly retrieved, and this remains stable for data with up to $20\%$ noise.
The mixed regularization is performed on a square $D$ (of width $1$) enclosing the two
modes in the index $\Phi$, see Fig. \ref{fig:2scb}(c). This choice of the inversion
domain $D$ allows possibly connecting of the modes. However, the enhancement correctly
recognizes two separate scatterers, with their magnitudes and sizes in excellent
agreement with the exact ones. Also the estimated background is very crispy.
Surprisingly, the estimate deteriorates only slightly in that the right scatterer
elongates a little bit towards the left scatterer as the noise level $\epsilon$ increases
from $0$ to $20\%$. Although not presented, we would like to note that for this
particular example, the reconstructions are still reasonable for data with $30\%$ noise.
Hence the proposed inverse scattering method is very robust with respect to the data
noise.

Next we consider a ring-shaped scatterer.
\begin{exam}\label{exam:ring}
The scatterer is one ring-shaped square located at the origin, with the outer and inner
side lengths being $0.6$ and $0.4$, respectively. The coefficient $\eta$ of the scatterer
is $1$. Two incident directions $d_1=\frac{1}{\sqrt{2}}(1,1)^\mathrm{T}$ and
$d_2=\frac{1}{\sqrt{2}}(1,-1)^\mathrm{T}$ are considered.
\end{exam}

\begin{figure}
  \centering
  \begin{tabular}{cccc}
    \includegraphics[width=.25\textwidth]{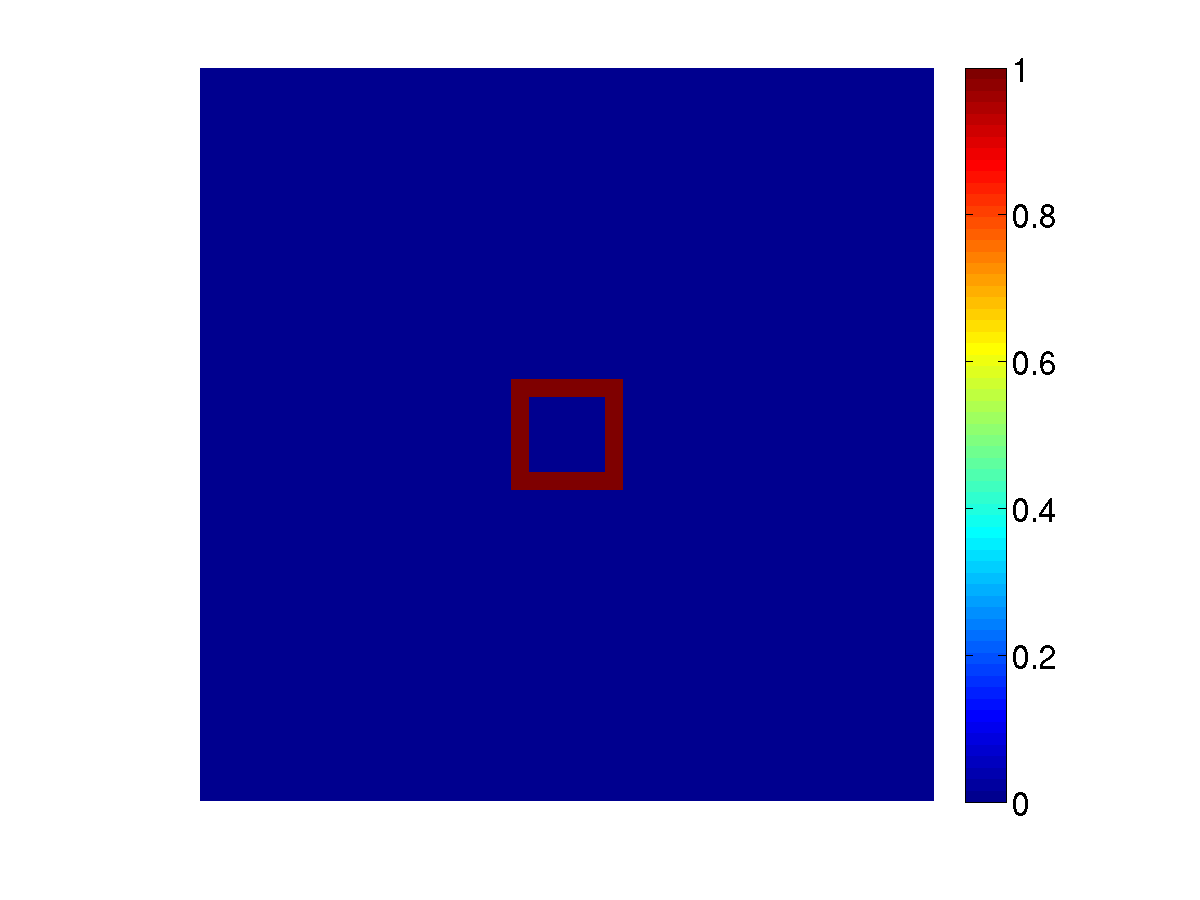} & \includegraphics[width=.25\textwidth]{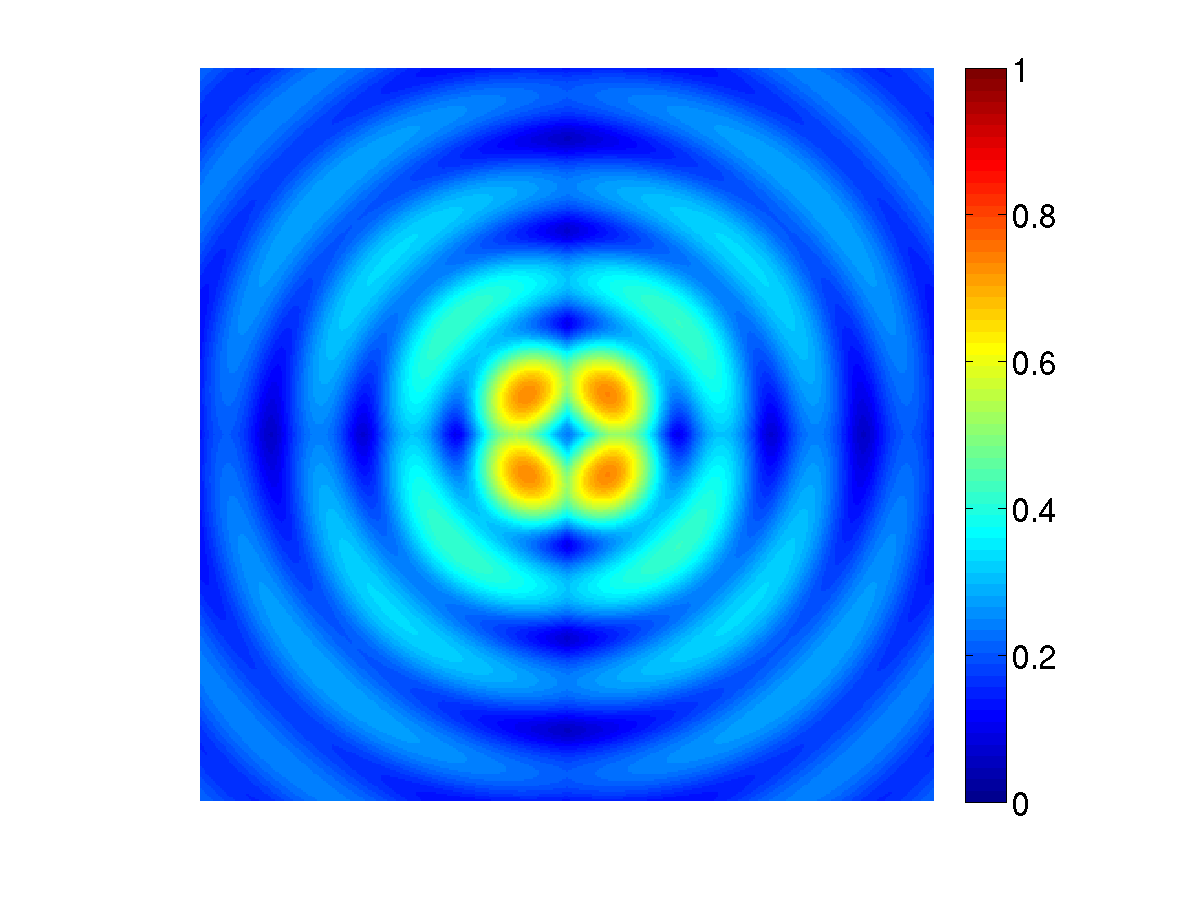}
   &\includegraphics[width=.25\textwidth]{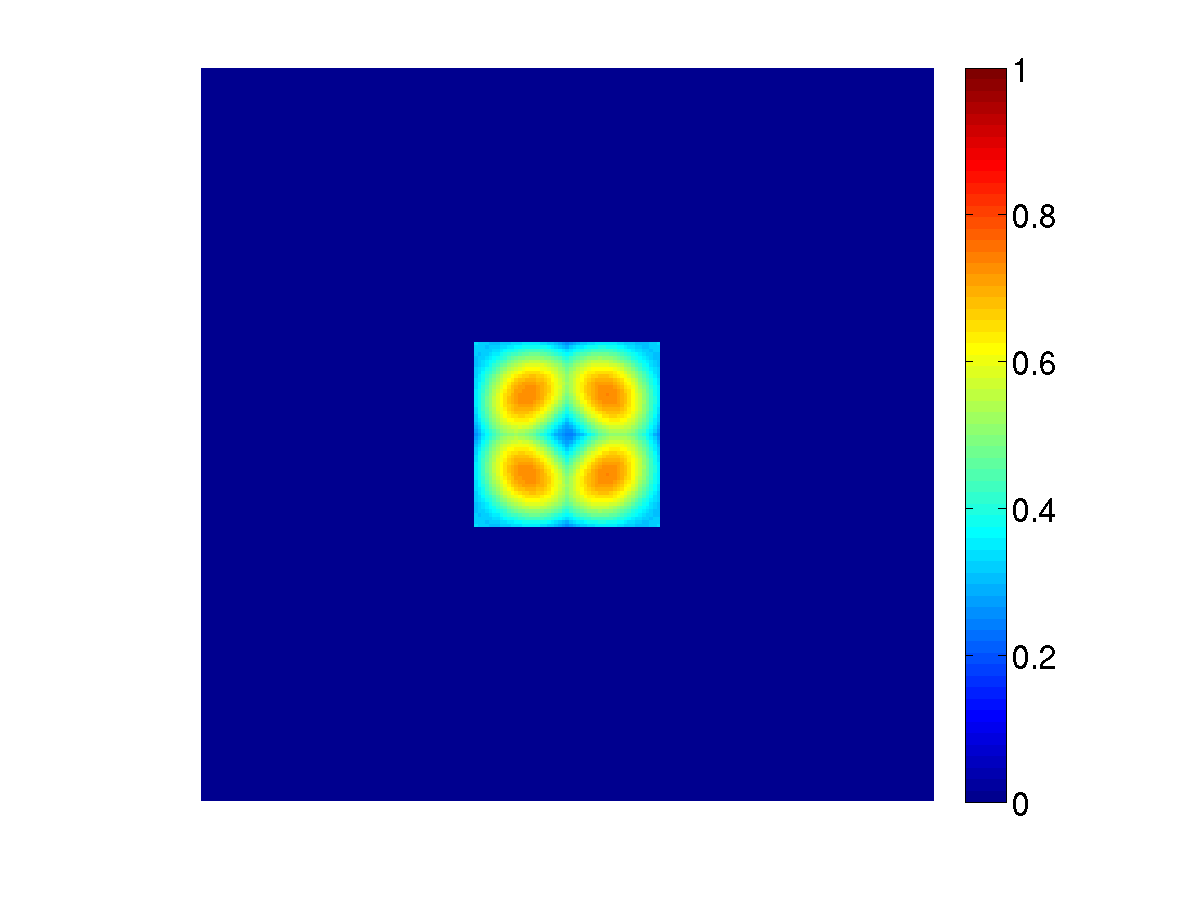} & \includegraphics[width=.25\textwidth]{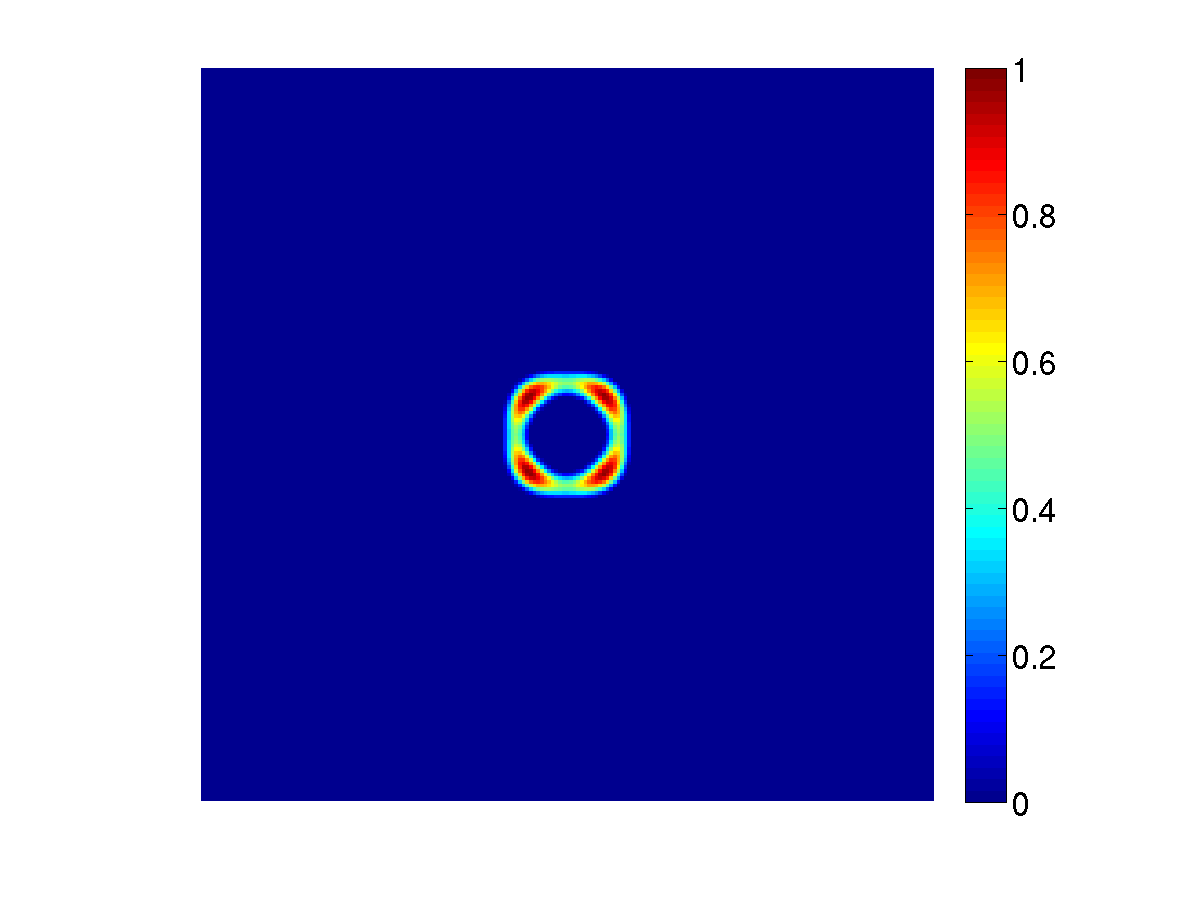}\\
    \includegraphics[width=.25\textwidth]{exam2_recon_true.png} & \includegraphics[width=.25\textwidth]{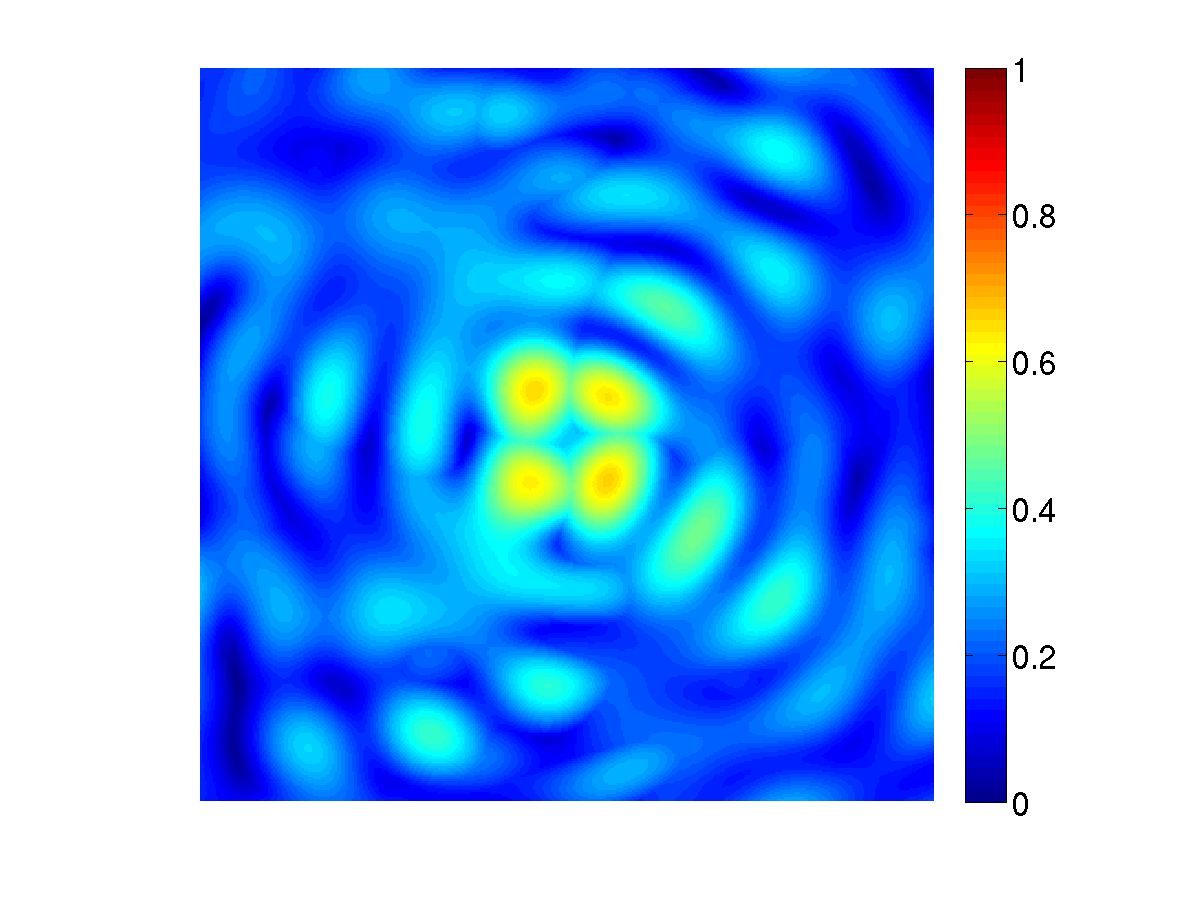}
   &\includegraphics[width=.25\textwidth]{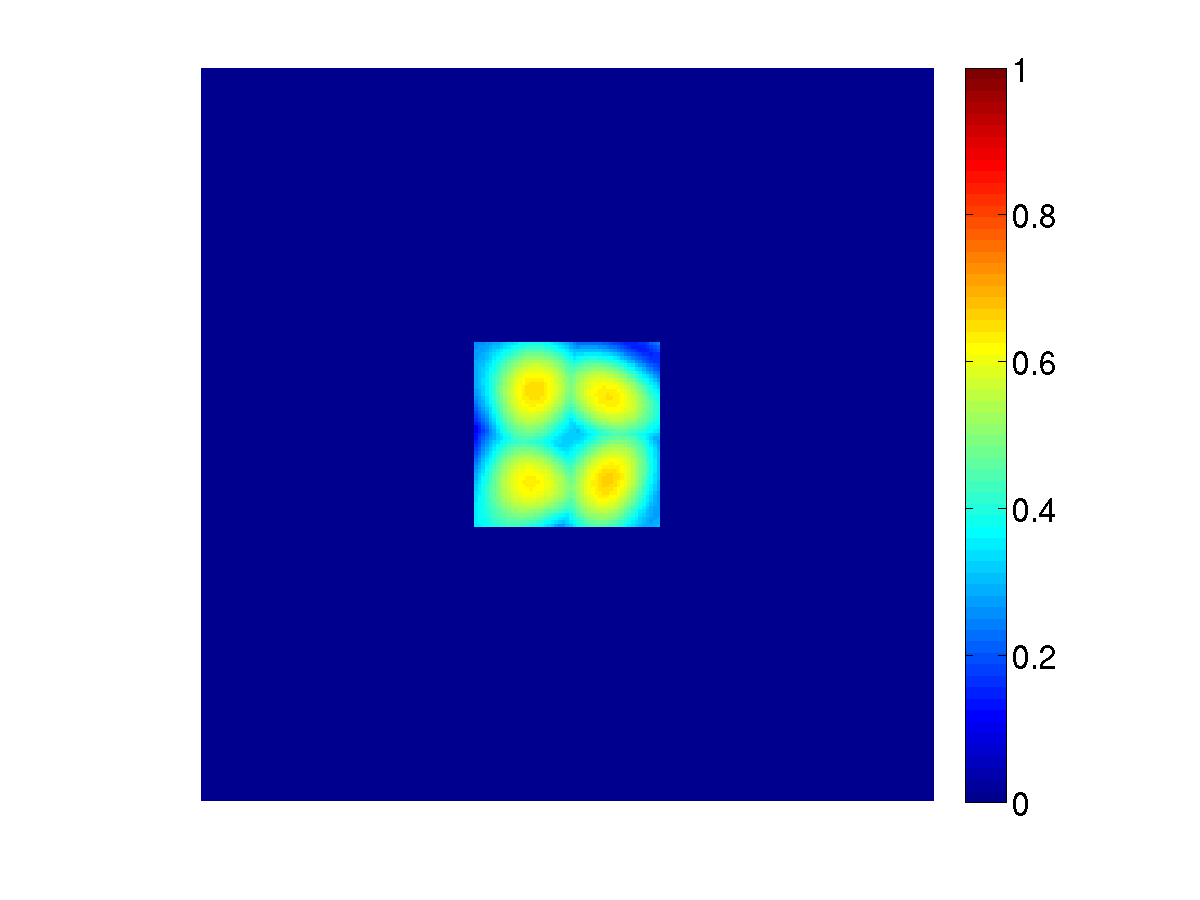}& \includegraphics[width=.25\textwidth]{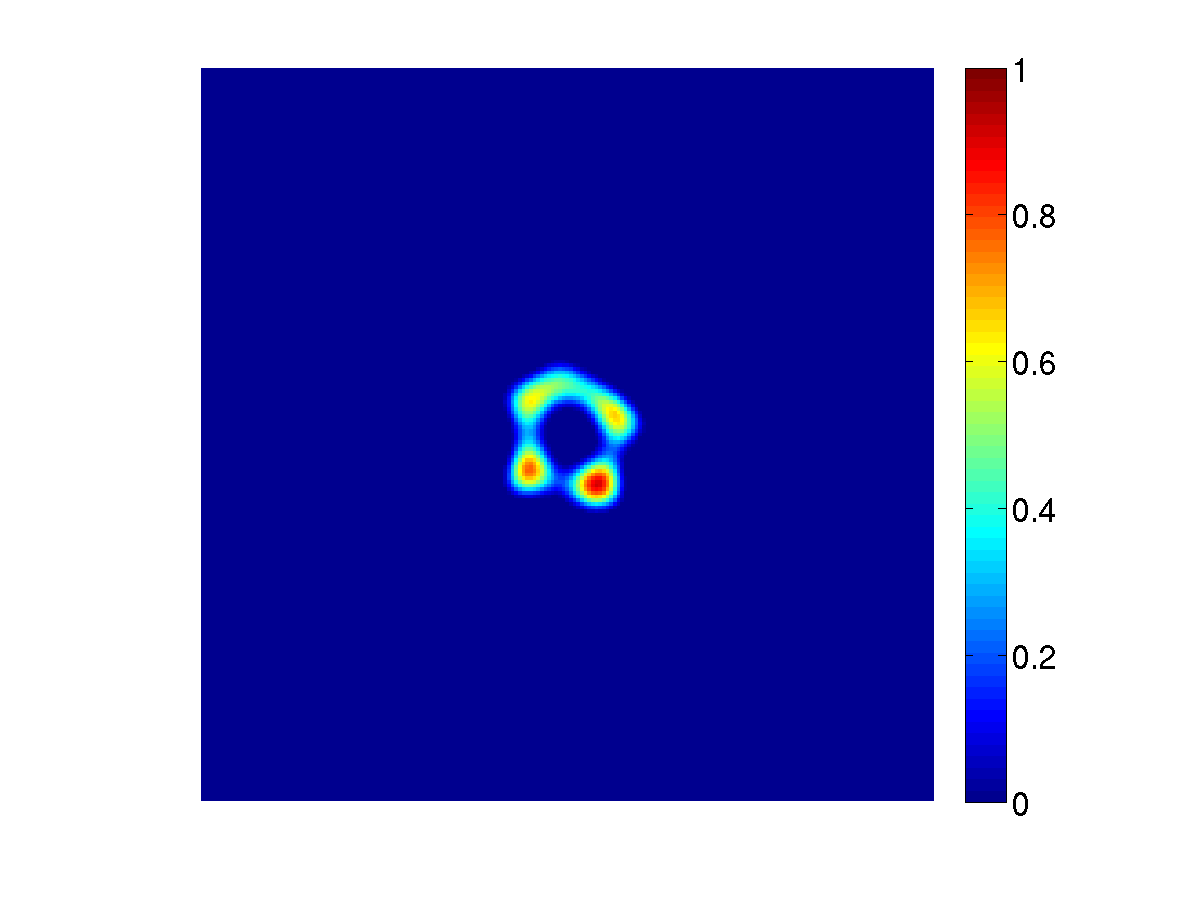}\\
   (a) true scatterer & (b) index $\Phi$ & (c) index $\Phi|_D$ & (d) sparse recon.
  \end{tabular}
  \caption{Numerical results for Example \ref{exam:ring}: (a) true scatterer, (b) index
   $\Phi$, (c) index $\Phi|_D$ (restriction to the subdomain $D$) and (d) sparse reconstruction.
   The first and second rows refer to exact data and the data with $20\%$ noise, respectively.}
  \label{fig:ring}
\end{figure}

Ring-shaped scatterer represents one of most challenging objects to recover, and it is
highly nontrivial even with multiple scattered field data sets, especially noting the
ring has a small thickness. It has been observed that one single incident field is
insufficient to completely resolve the ring structure, and only some parts of the ring
can be resolved, depending on the incident direction $d$ \cite{ItoJinZou:2011a}. Hence we
employ two incident waves in the directions $d_1=\frac{1}{\sqrt{2}}(1,1)^\mathrm{T}$ and
$d_2=\frac{1}{\sqrt{2}}(1,-1)^\mathrm{T}$ in order to yield sufficient amount of
information about the scatterer, and accordingly, the index function $\Phi$ is defined as
follows
\begin{equation*}
  \Phi(x_p)=\max_i\{\Phi_i(x_p)\}\quad \forall x_p\in\widetilde{\Omega},
\end{equation*}
where the function $\Phi_i$ refers to the index for the $i$th data set. The numerical
results with the exact data and $20\%$ noise in the data are shown in Fig.
\ref{fig:ring}. With just two incident waves, the index $\Phi$ can provide a quite
reasonable estimate of the ring shape. Despite some small oscillations, the overall
profile stands out clearly, and remains very stable for up to $20\%$ noise in the data.
The enhancing step via mixed regularization provides a very crispy estimate of the ring
structure: the recovered scatterer has a clear ring structure, which agrees excellently
with the exact one in terms of both magnitude and size. The presence of $20\%$ data noise
causes visible deterioration to the reconstruction, see Fig. \ref{fig:ring}(d).
Nonetheless, the enhanced reconstruction still exhibits a clear ring shape, and it
represents a very good approximation to the true scatterer upon noting the large amount
of data noise.

\subsection{Three-dimensional example}

Our last example shows the feasibility of the method for three-dimensional problems.
\begin{exam}\label{exam:cube}
We consider two cubic scatterers of width $0.1$ centered at $(0.35,0.15,0.15)$ and
$(-0.35,0.15,0.15)$, respectively. One single incident field with direction
$d=\frac{1}{\sqrt{3}}(1,1,1)^\mathrm{T}$ is used, and the coefficient $\eta$ of the
scatterers is taken to be $1$.
\end{exam}

The scattered field $u^s$ is measured at $600$ points uniformly distributed on the
surface $\Gamma$ of a cubic of width $5$, (i.e., $10$ points in each direction). To
simulate the scattered field data, we take the sampling domain $\widetilde{\Omega}$ to be
the cubic $[-1,1]^3$, which is divided into a uniform mesh consisting of small cubes of
width $h=0.01$. The inversion domain $D$ for the integral equation \eqref{eqn:linint} is
divided into a coarser mesh consisting of small cubes of width $0.03$.

\begin{figure}
  \centering
  \begin{tabular}{cccc}
     \includegraphics[width=.25\textwidth]{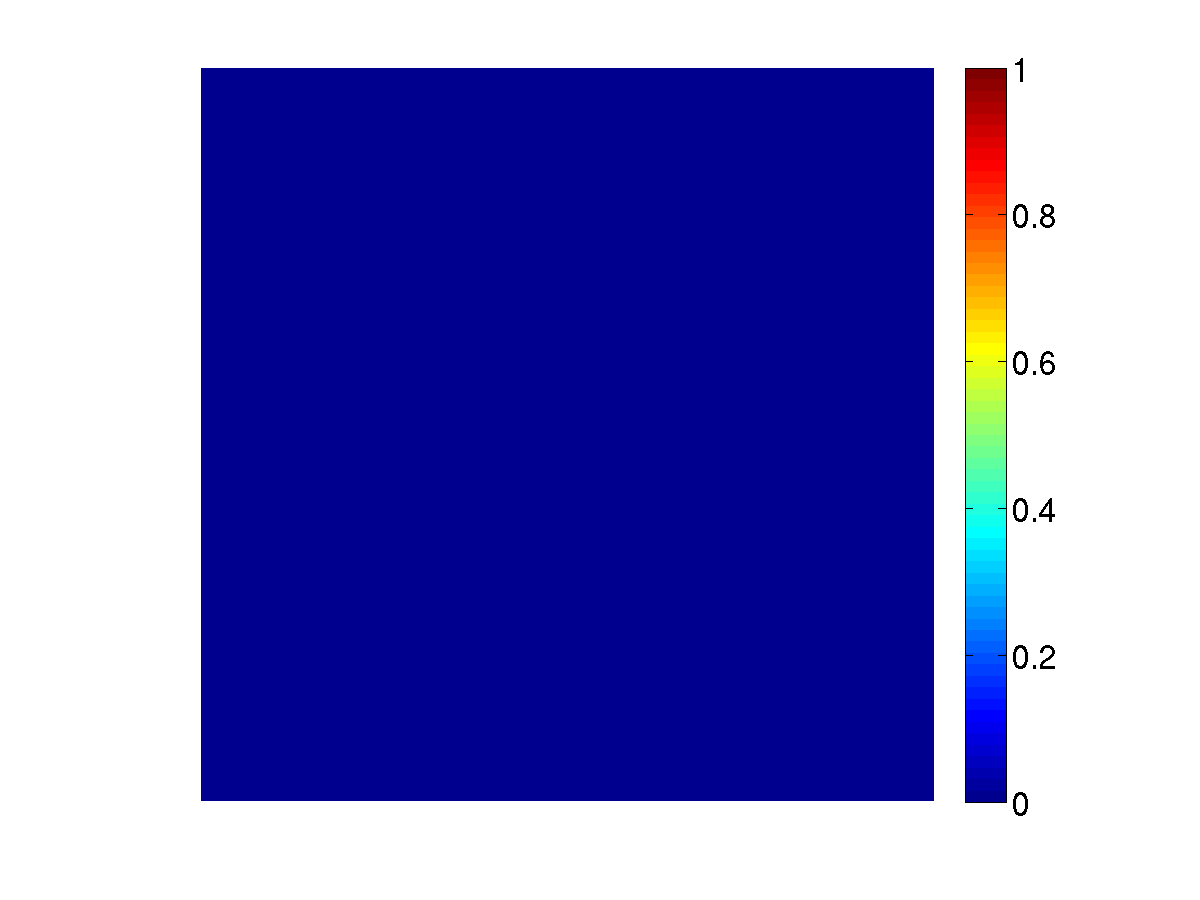} & \includegraphics[width=.25\textwidth]{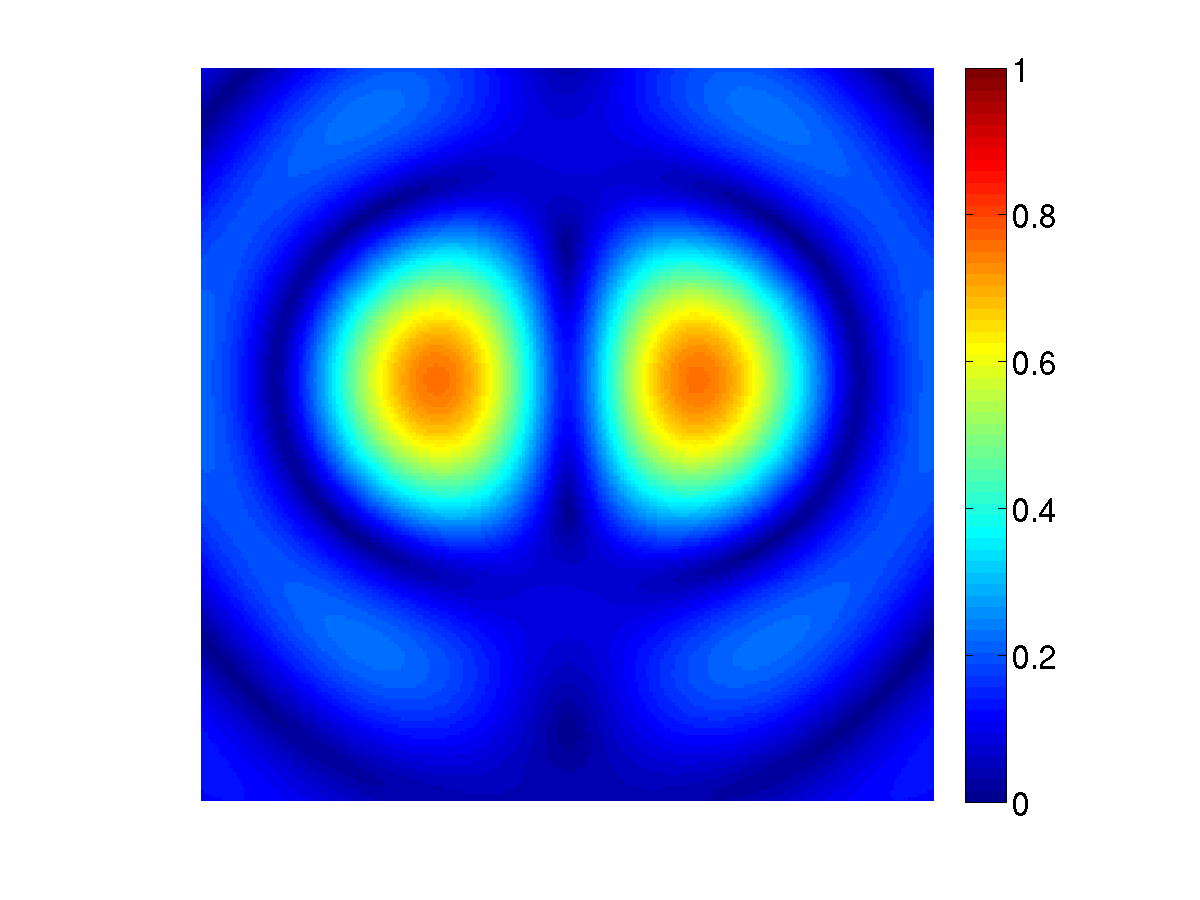}
    &\includegraphics[width=.25\textwidth]{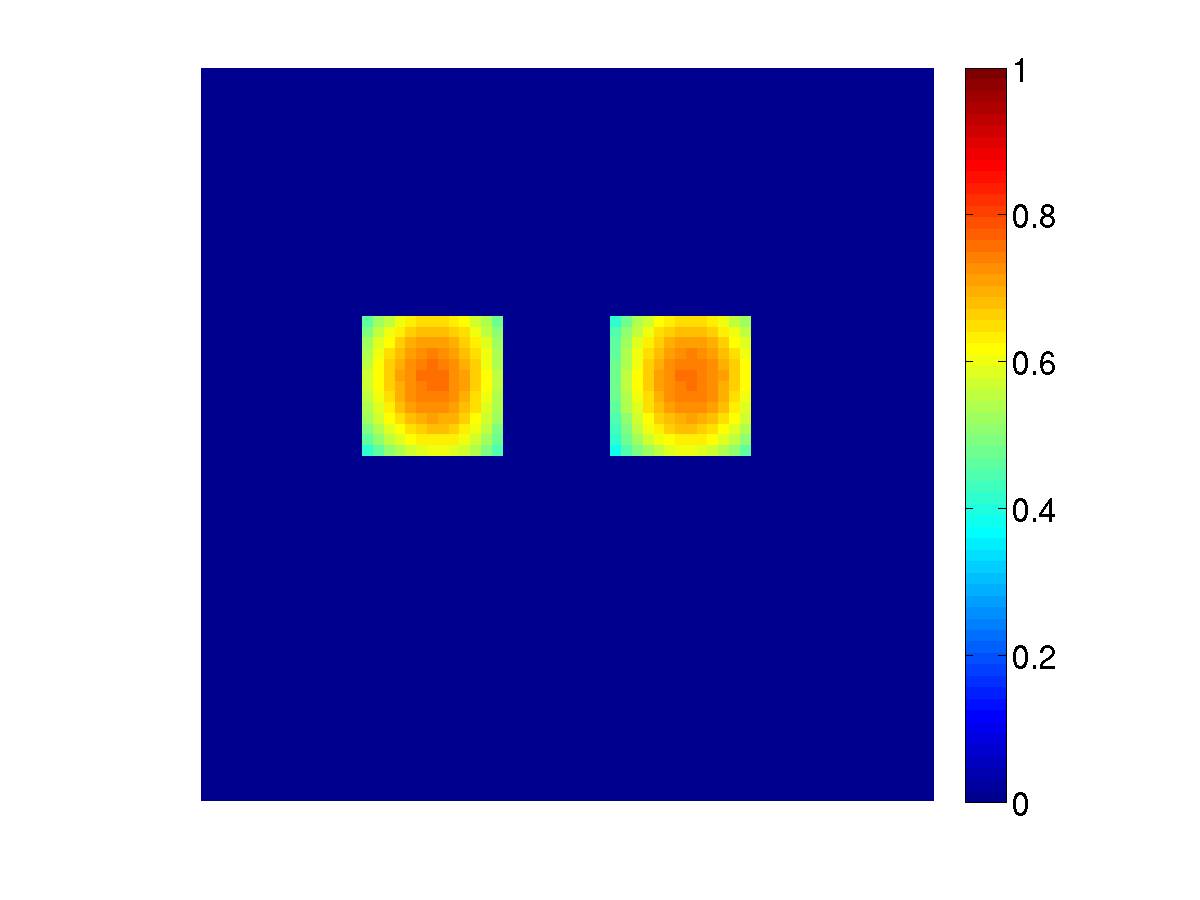} & \includegraphics[width=.25\textwidth]{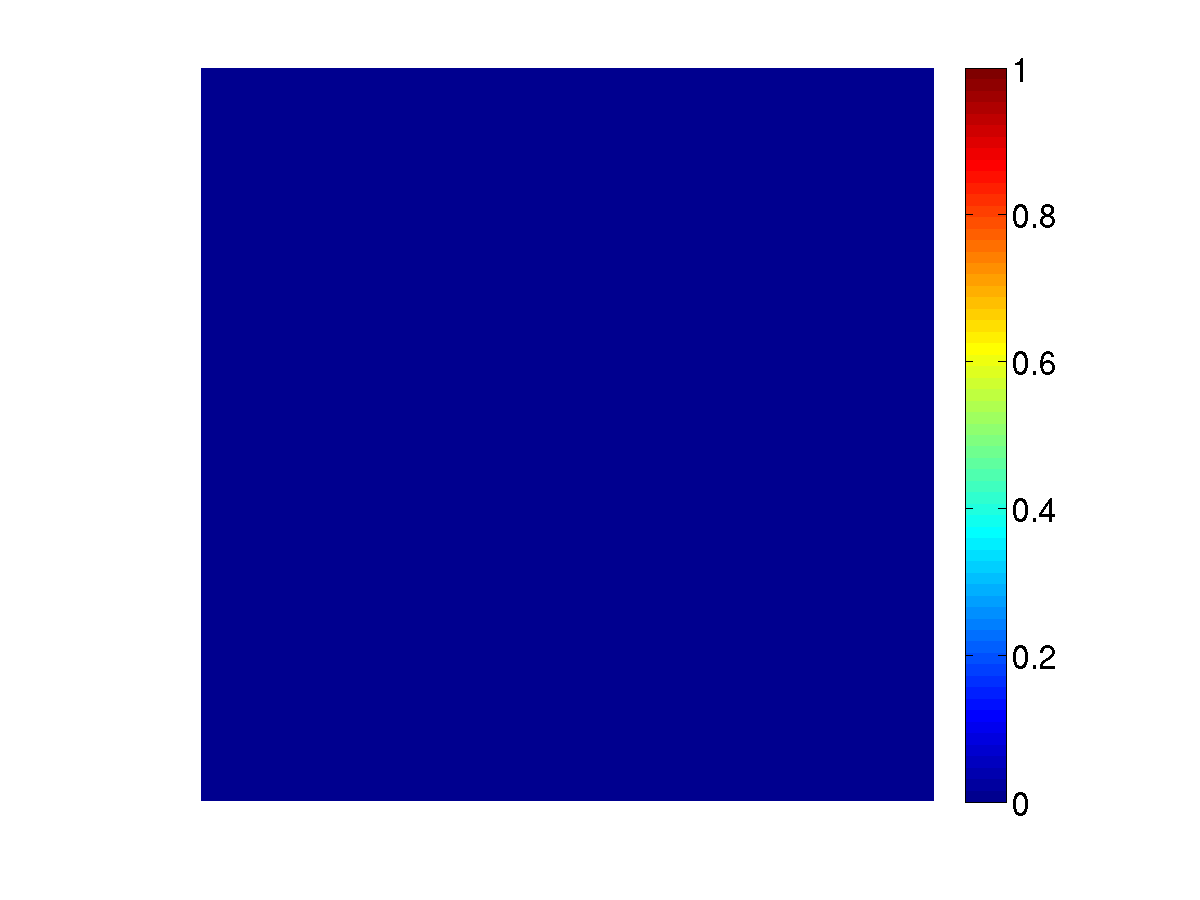}\\
     \includegraphics[width=.25\textwidth]{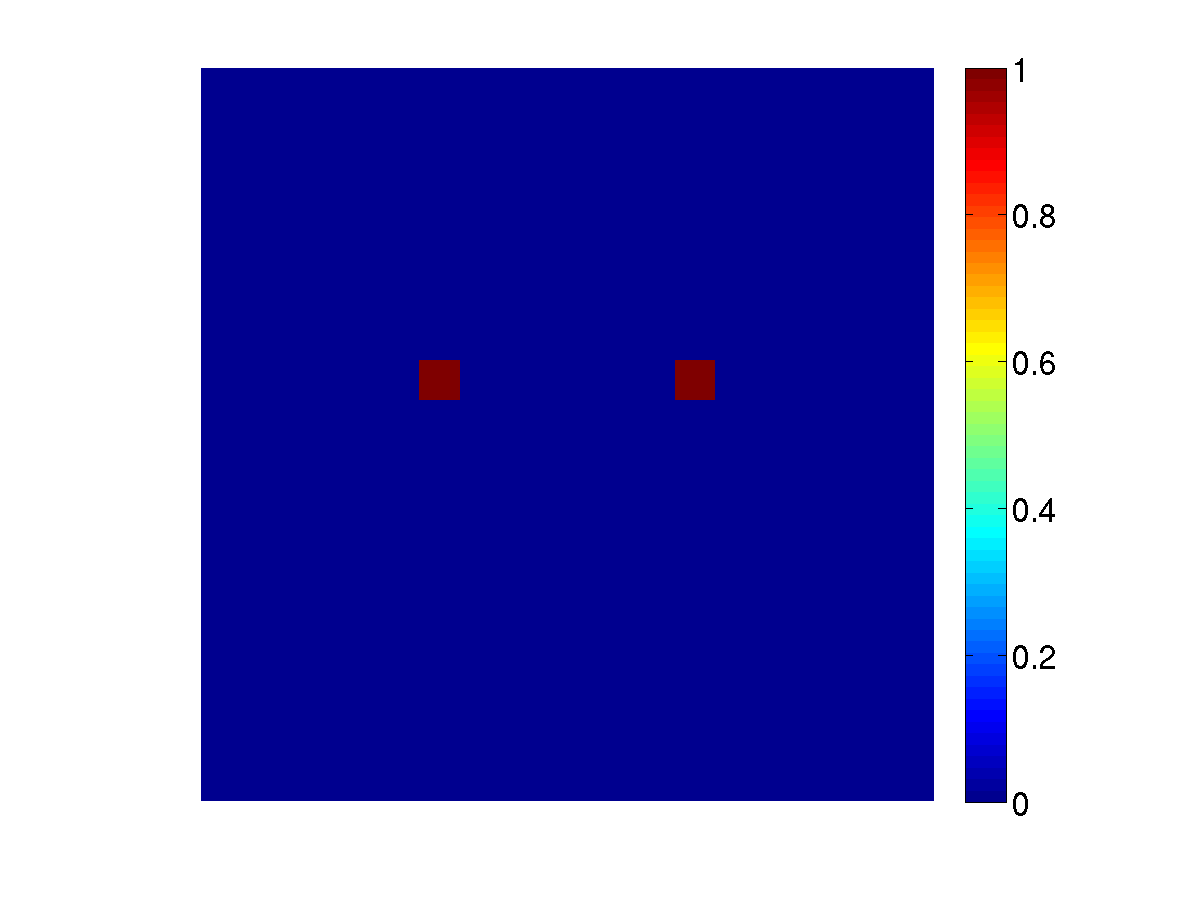} & \includegraphics[width=.25\textwidth]{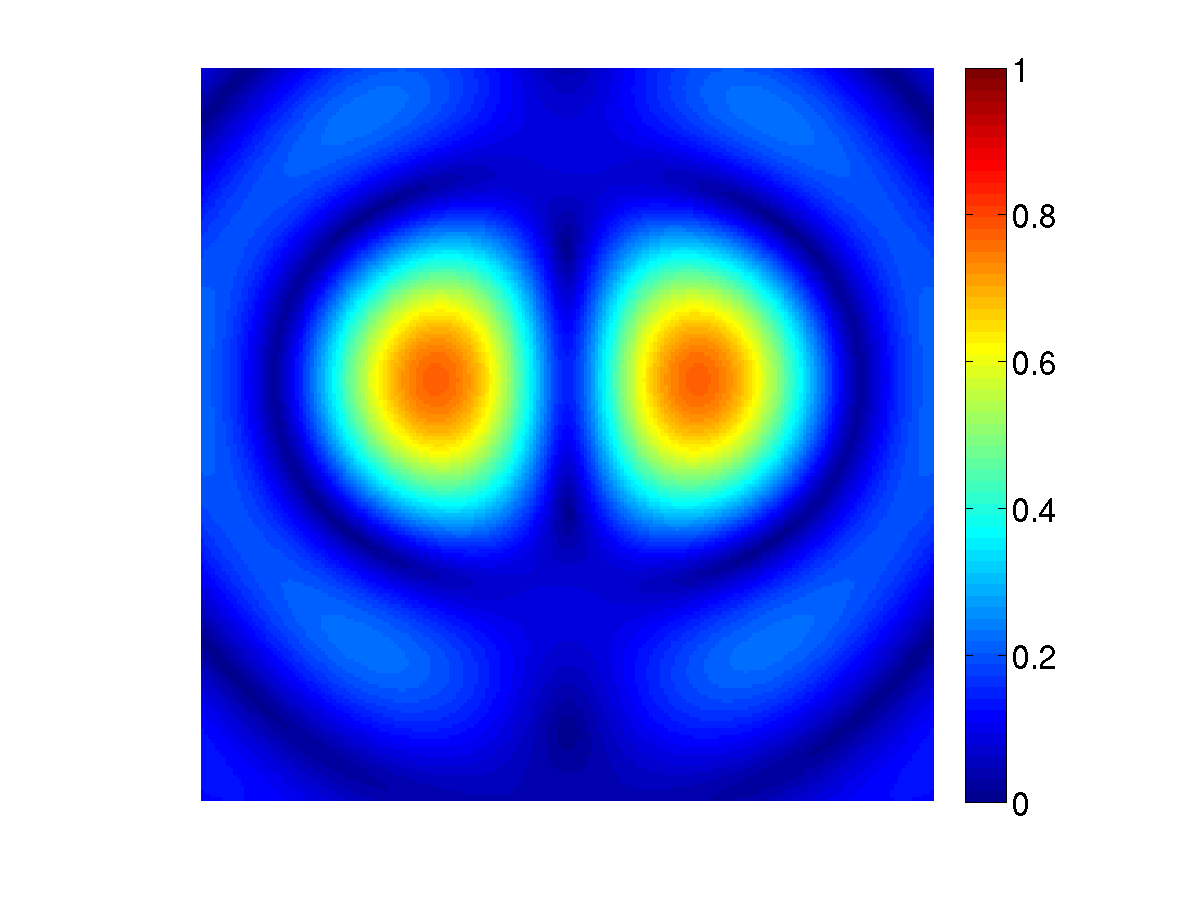}
    &\includegraphics[width=.25\textwidth]{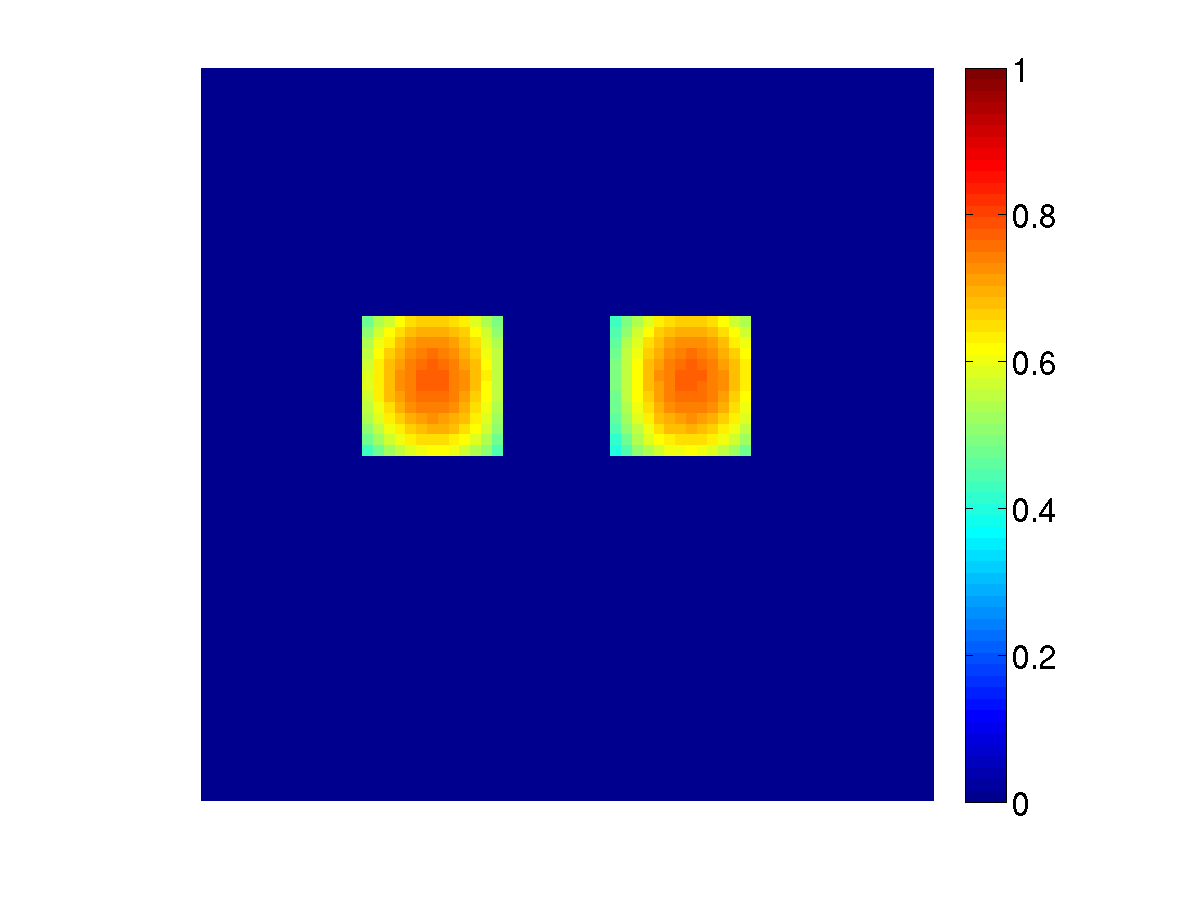} & \includegraphics[width=.25\textwidth]{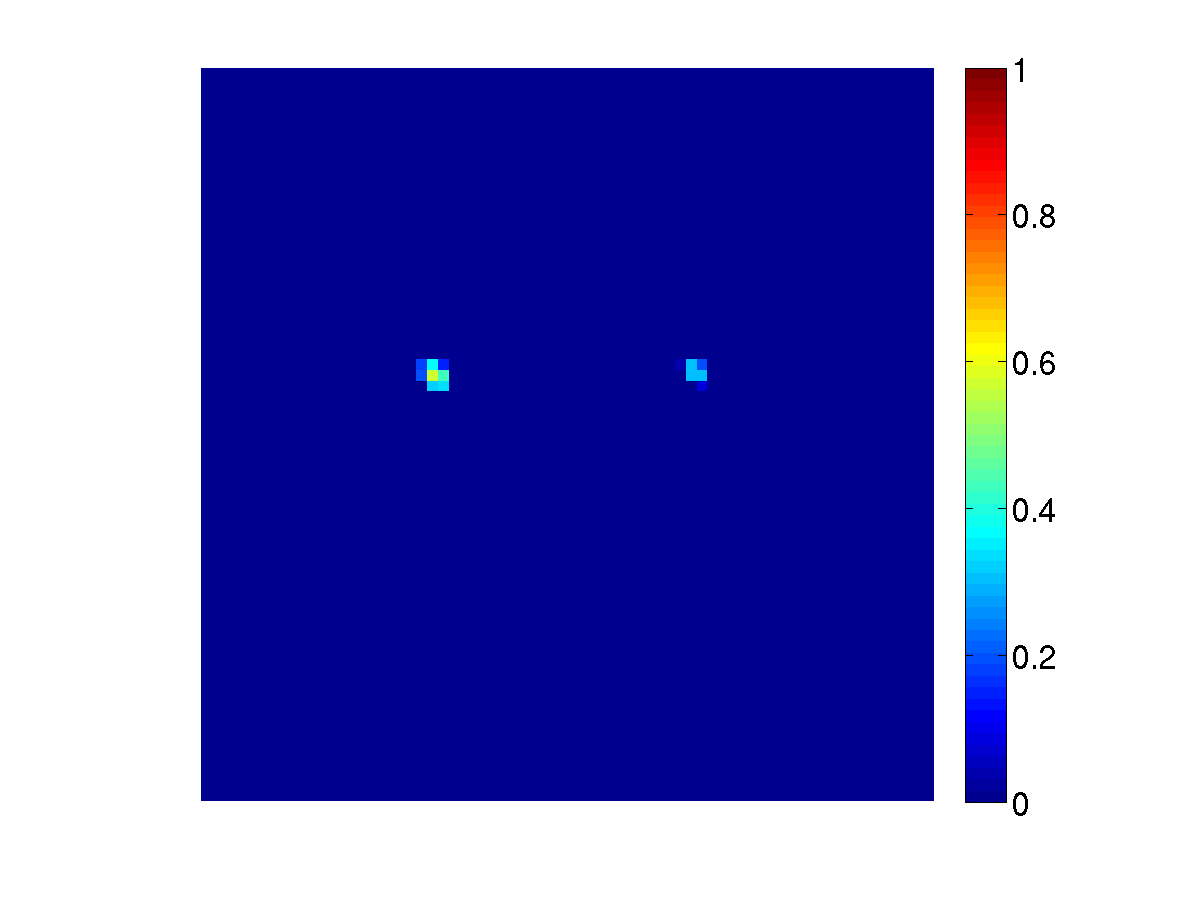}\\
     \includegraphics[width=.25\textwidth]{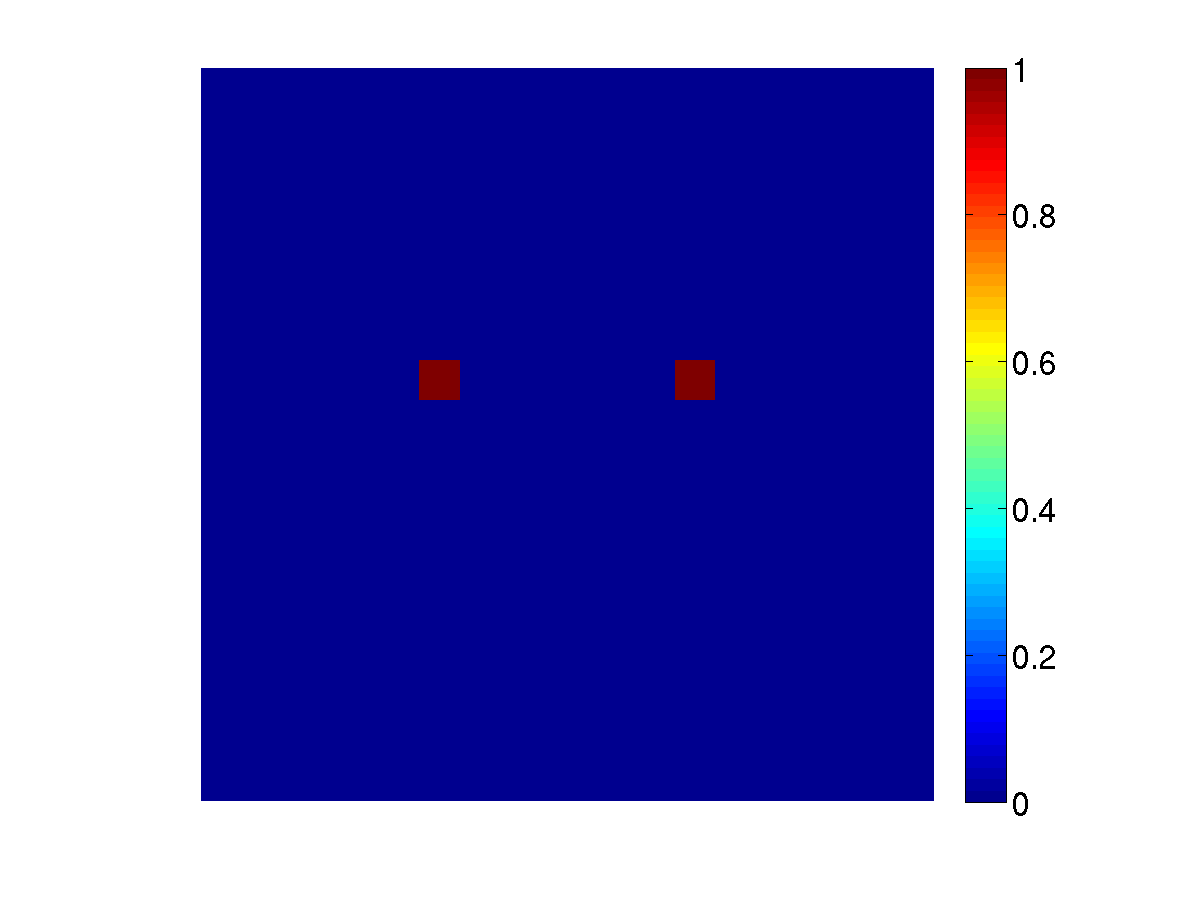} & \includegraphics[width=.25\textwidth]{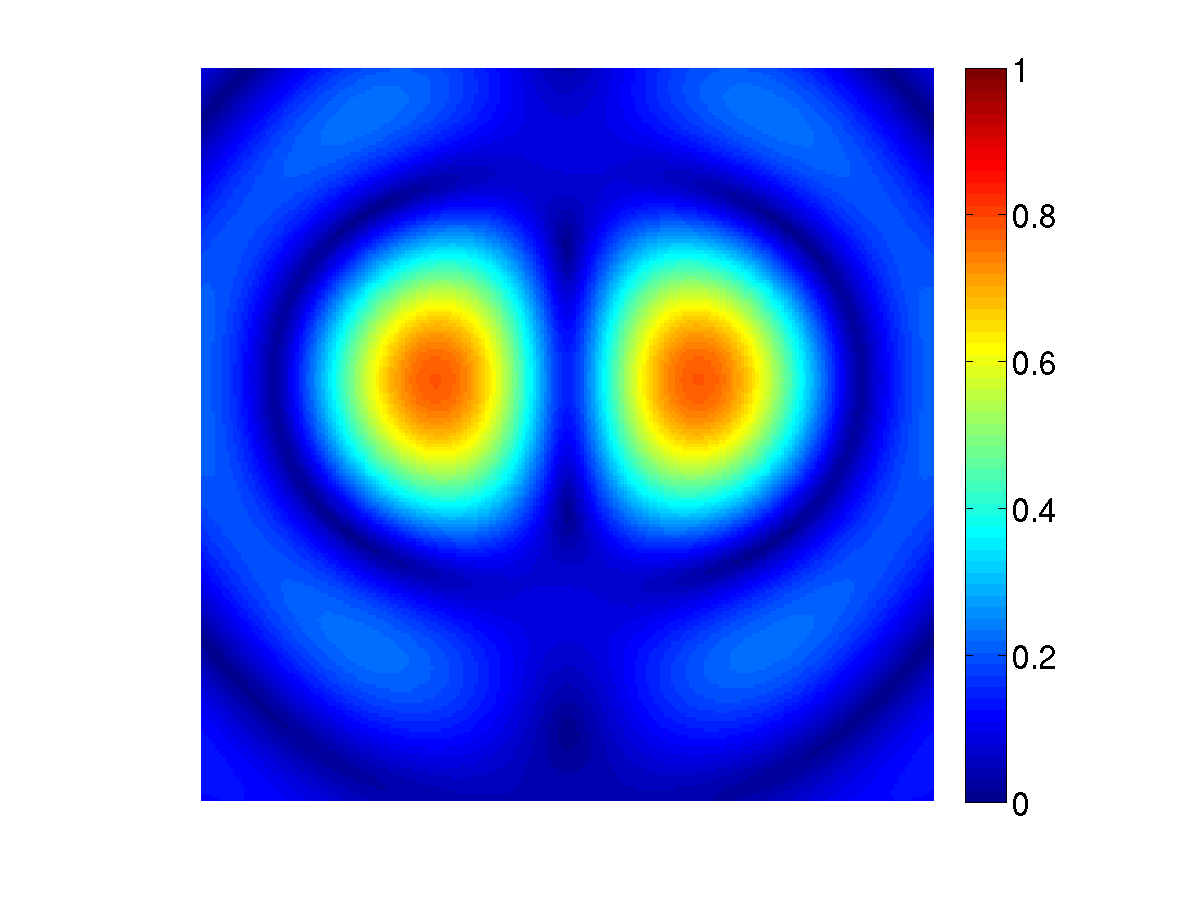}
    &\includegraphics[width=.25\textwidth]{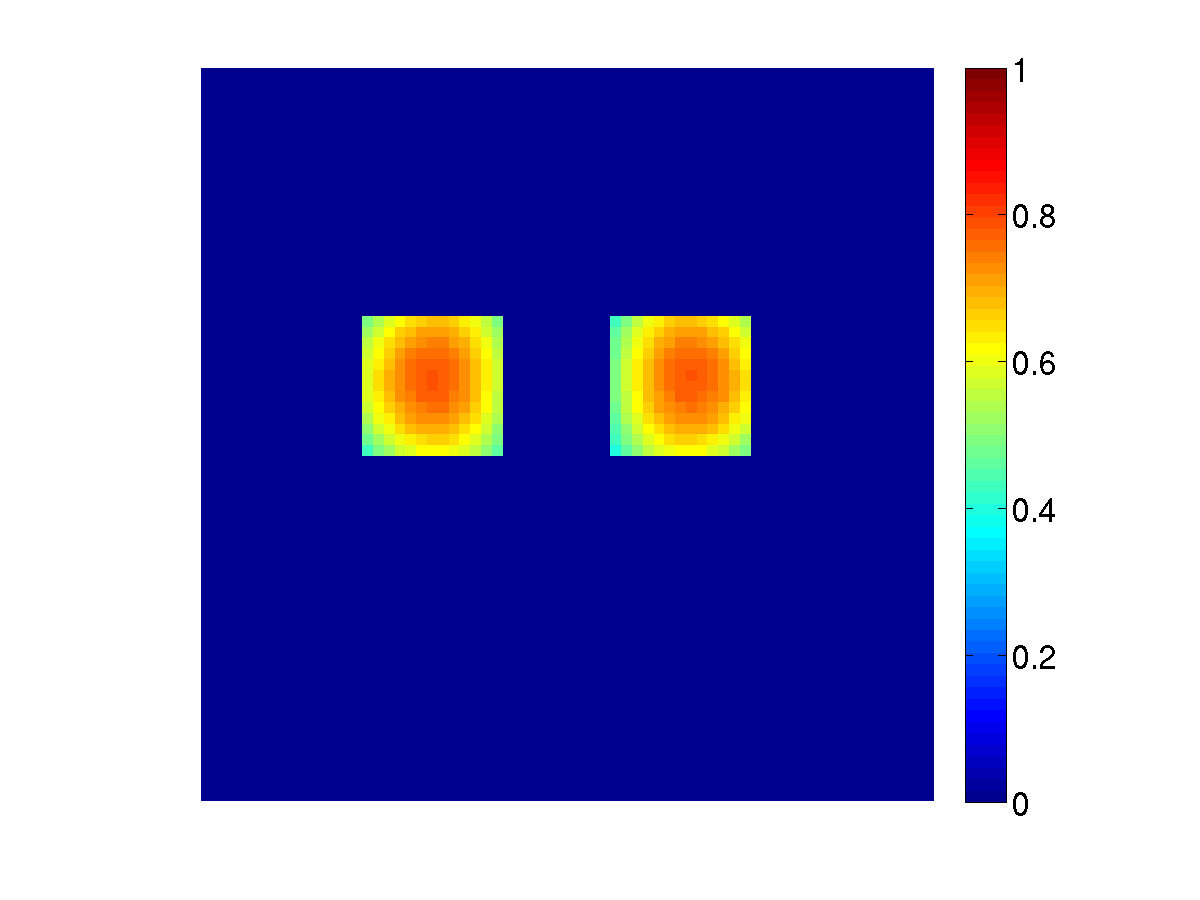} & \includegraphics[width=.25\textwidth]{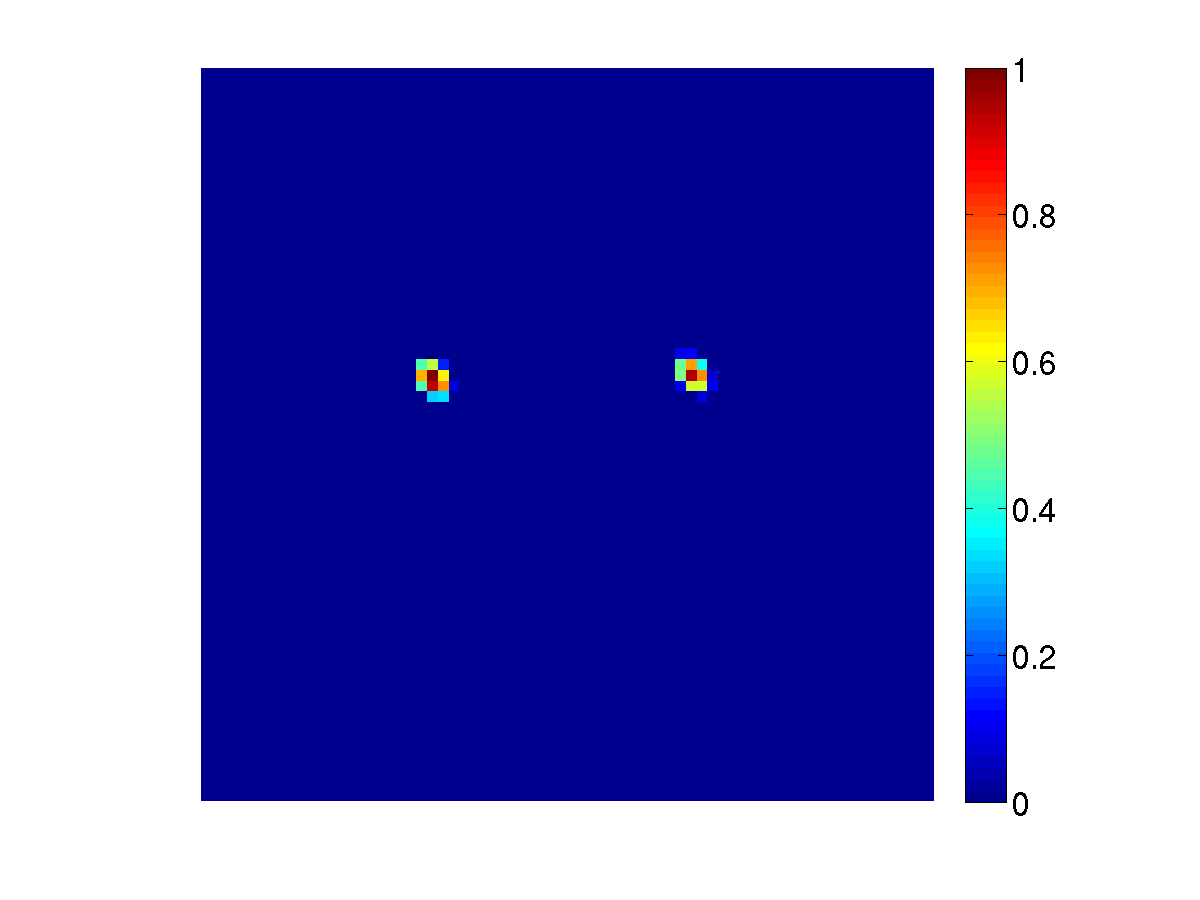}\\
     \includegraphics[width=.25\textwidth]{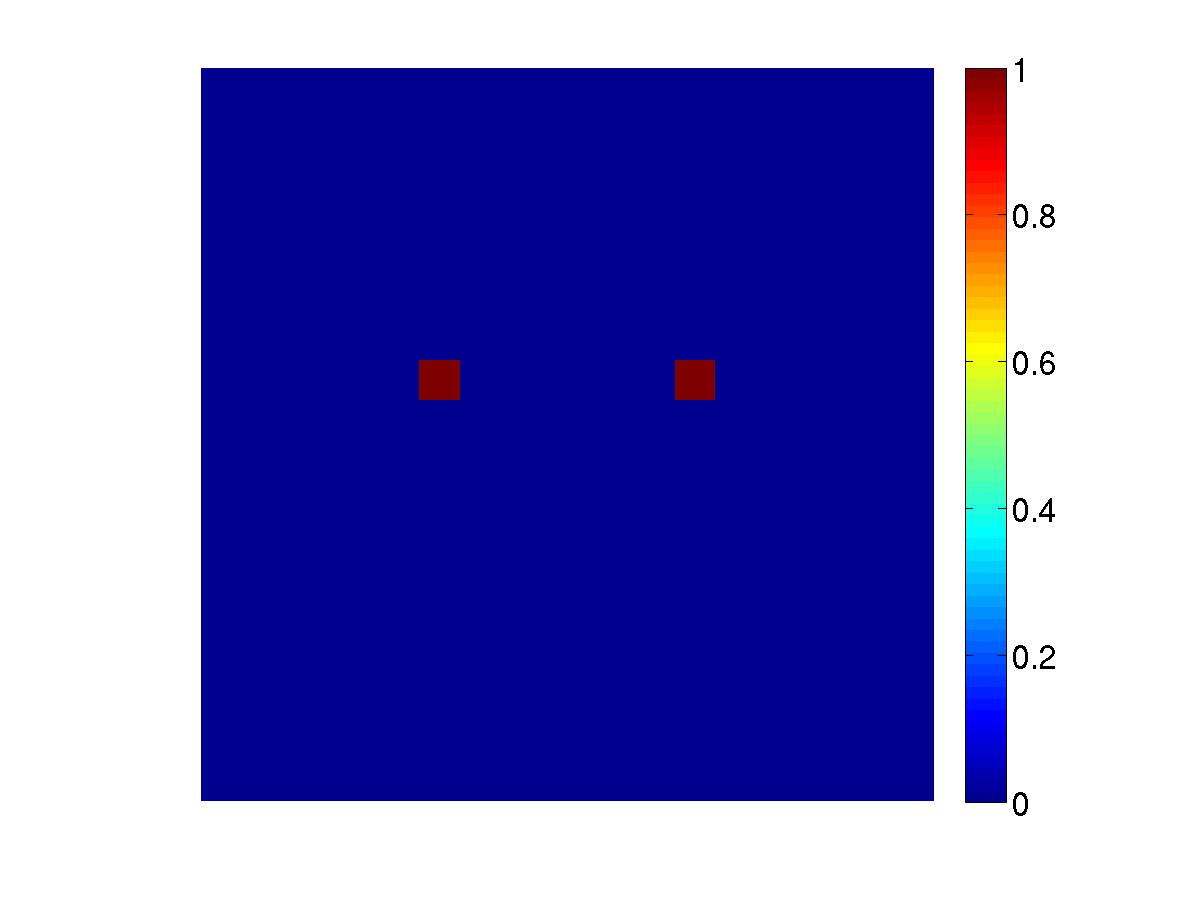} & \includegraphics[width=.25\textwidth]{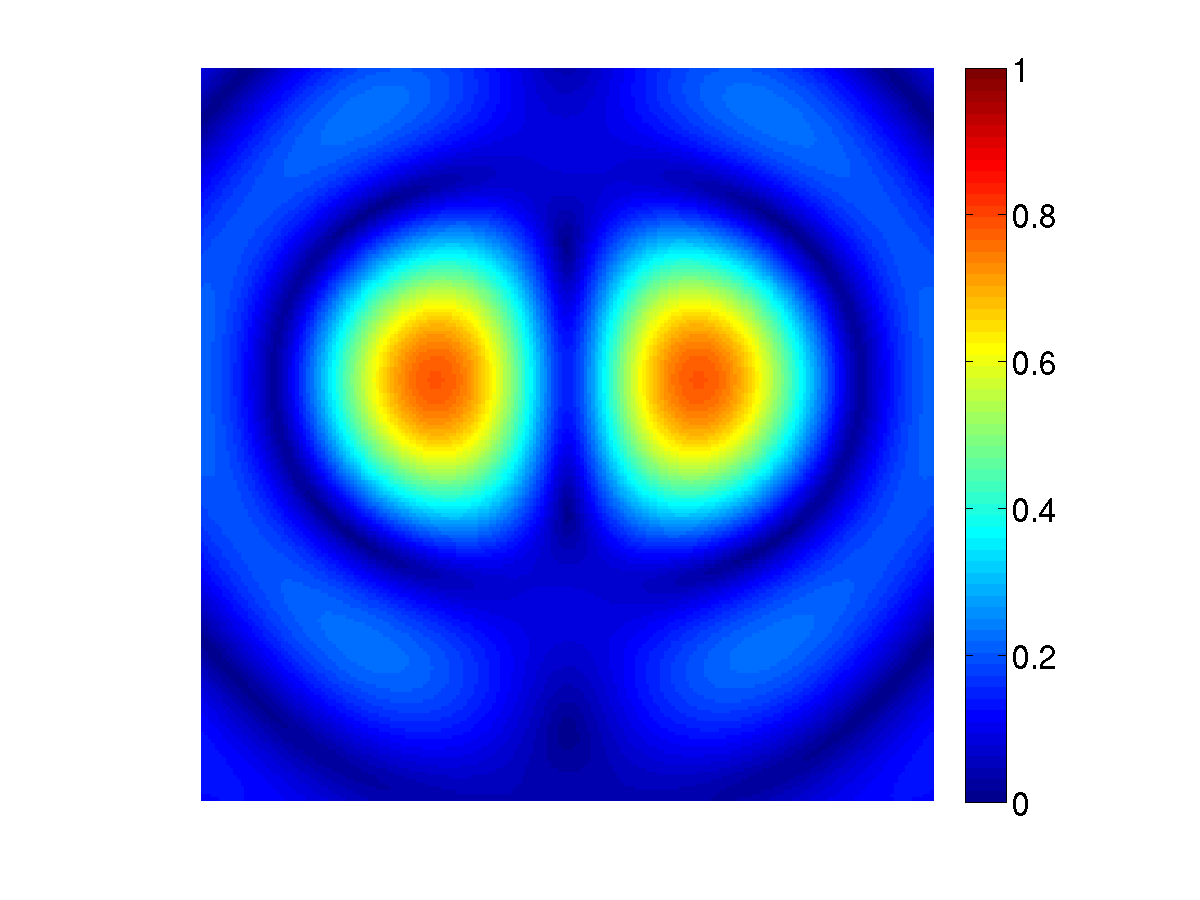}
    &\includegraphics[width=.25\textwidth]{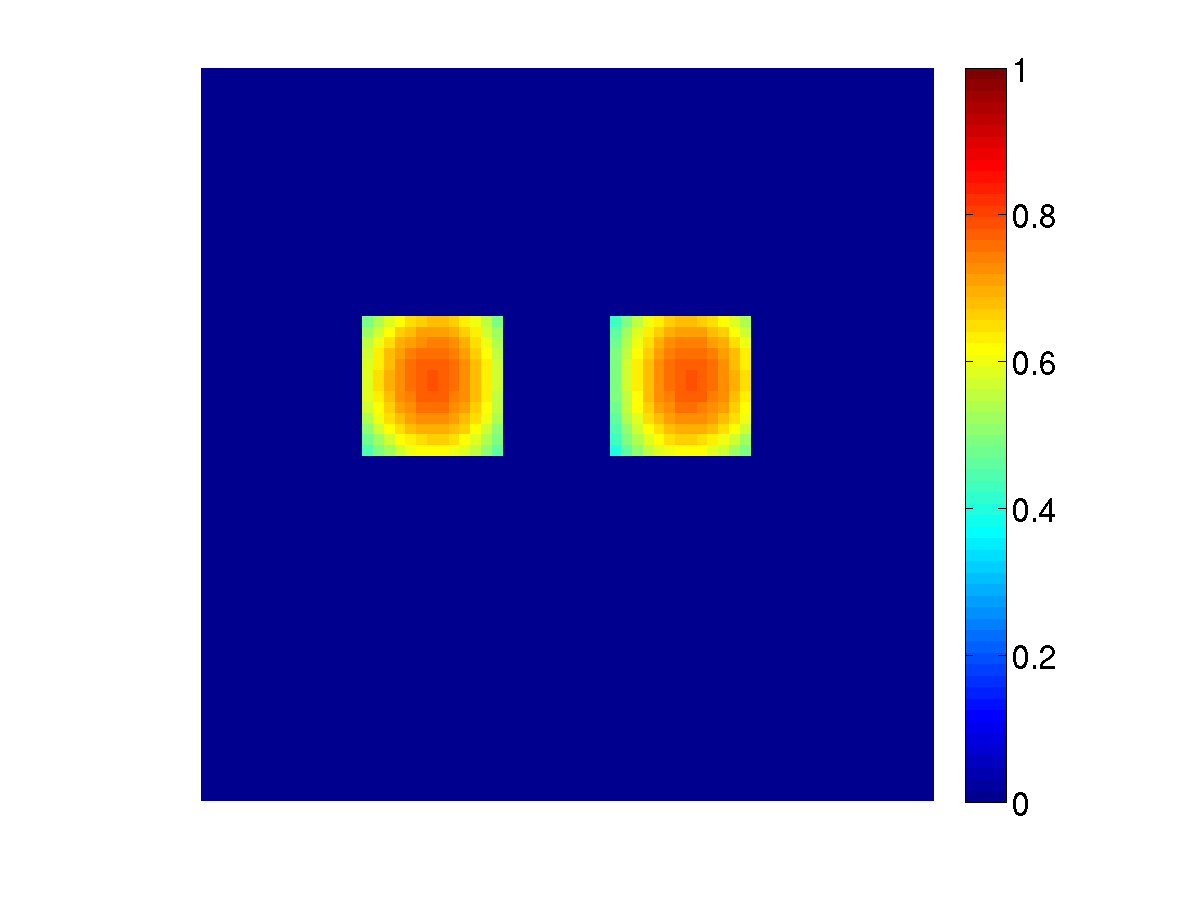} & \includegraphics[width=.25\textwidth]{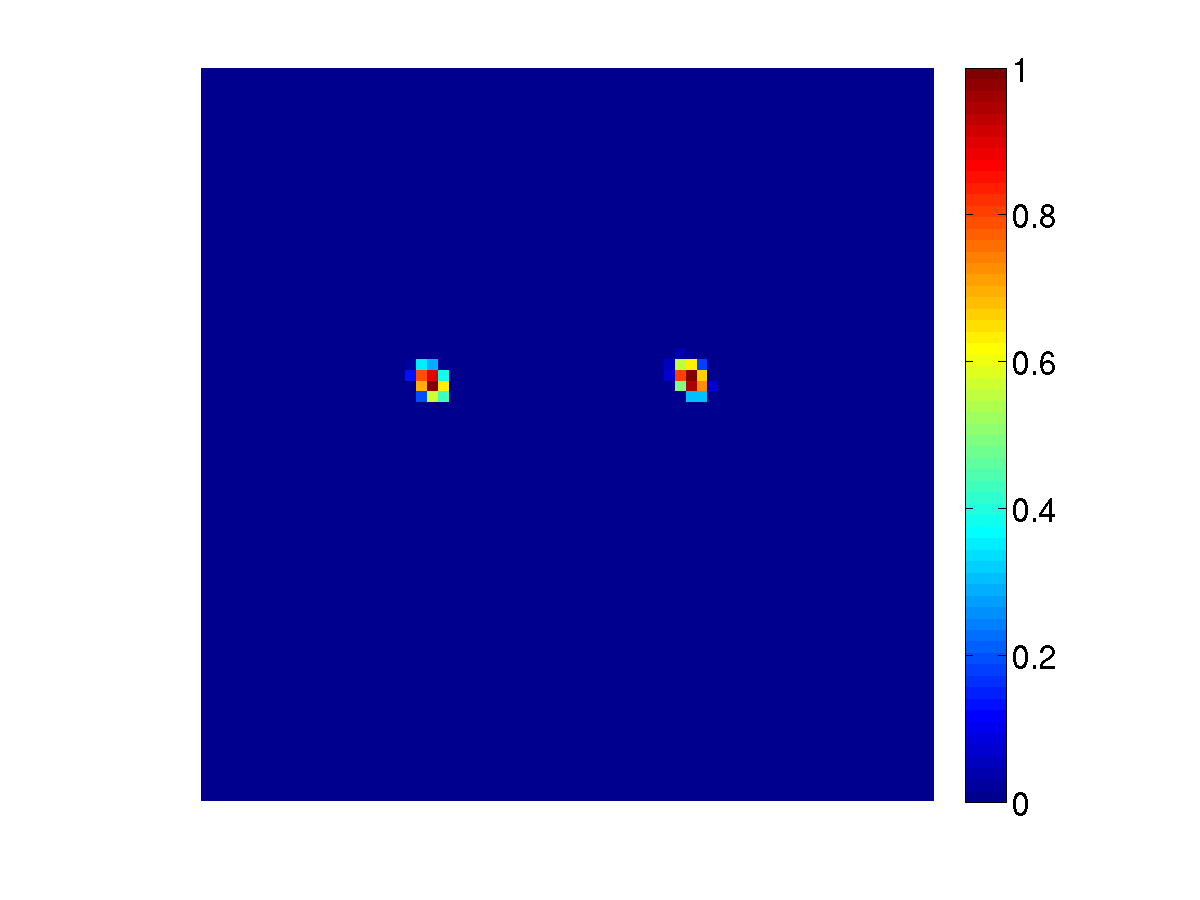}\\
     \includegraphics[width=.25\textwidth]{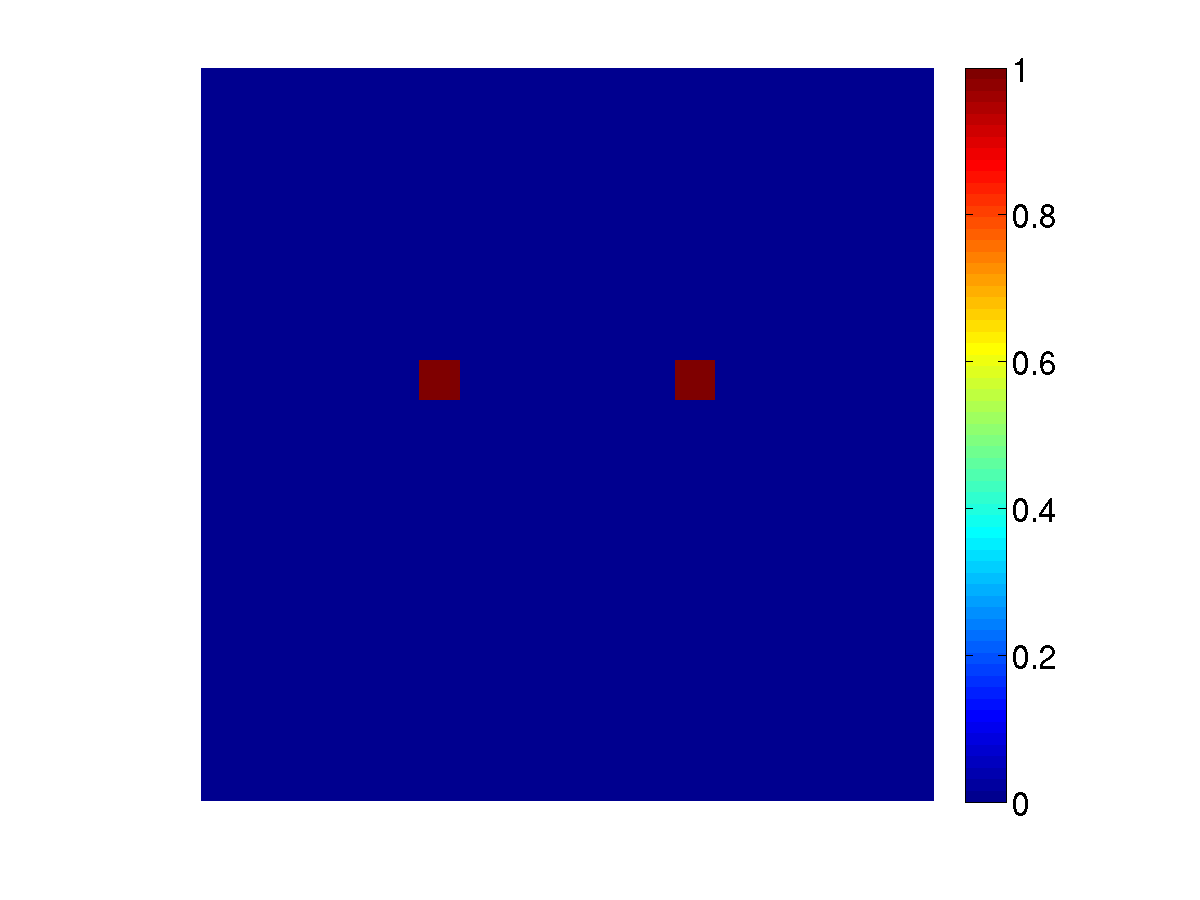} & \includegraphics[width=.25\textwidth]{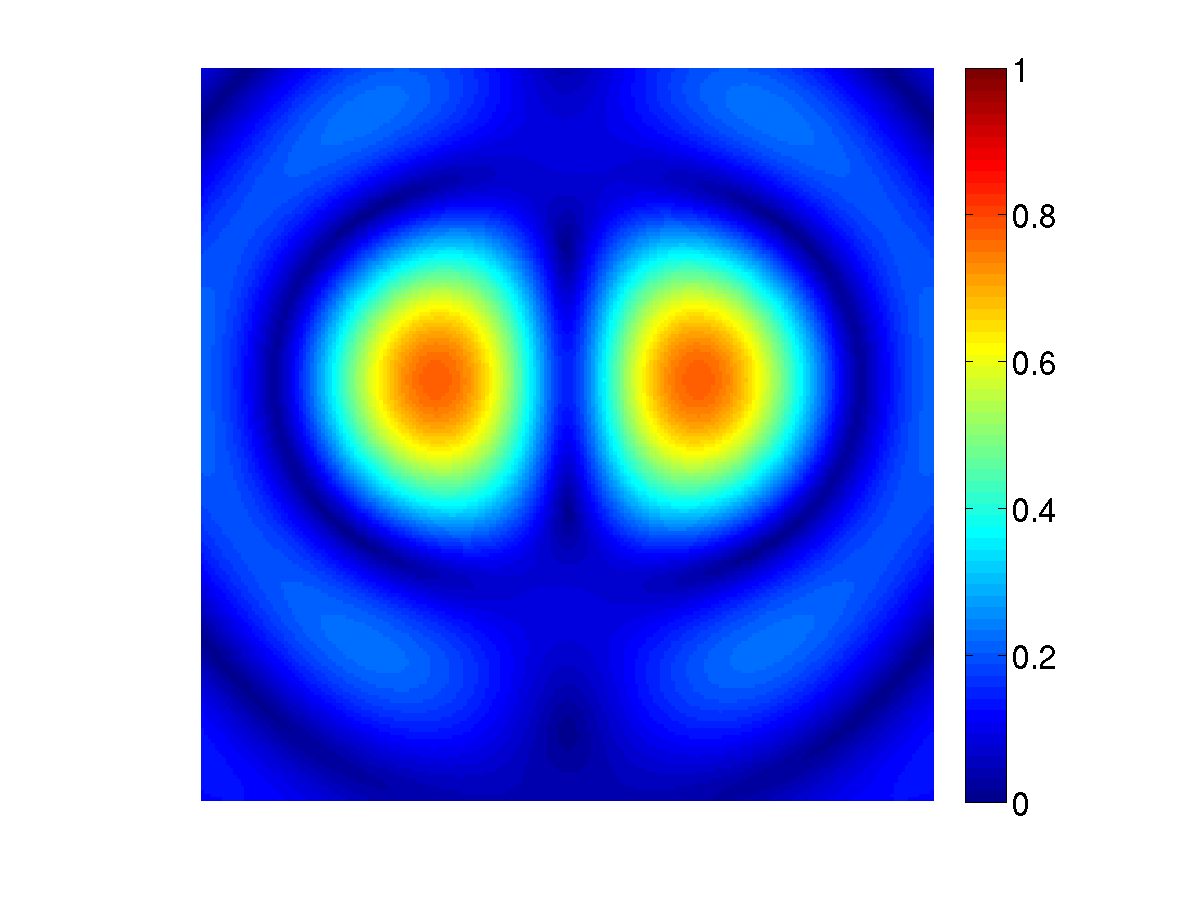}
    &\includegraphics[width=.25\textwidth]{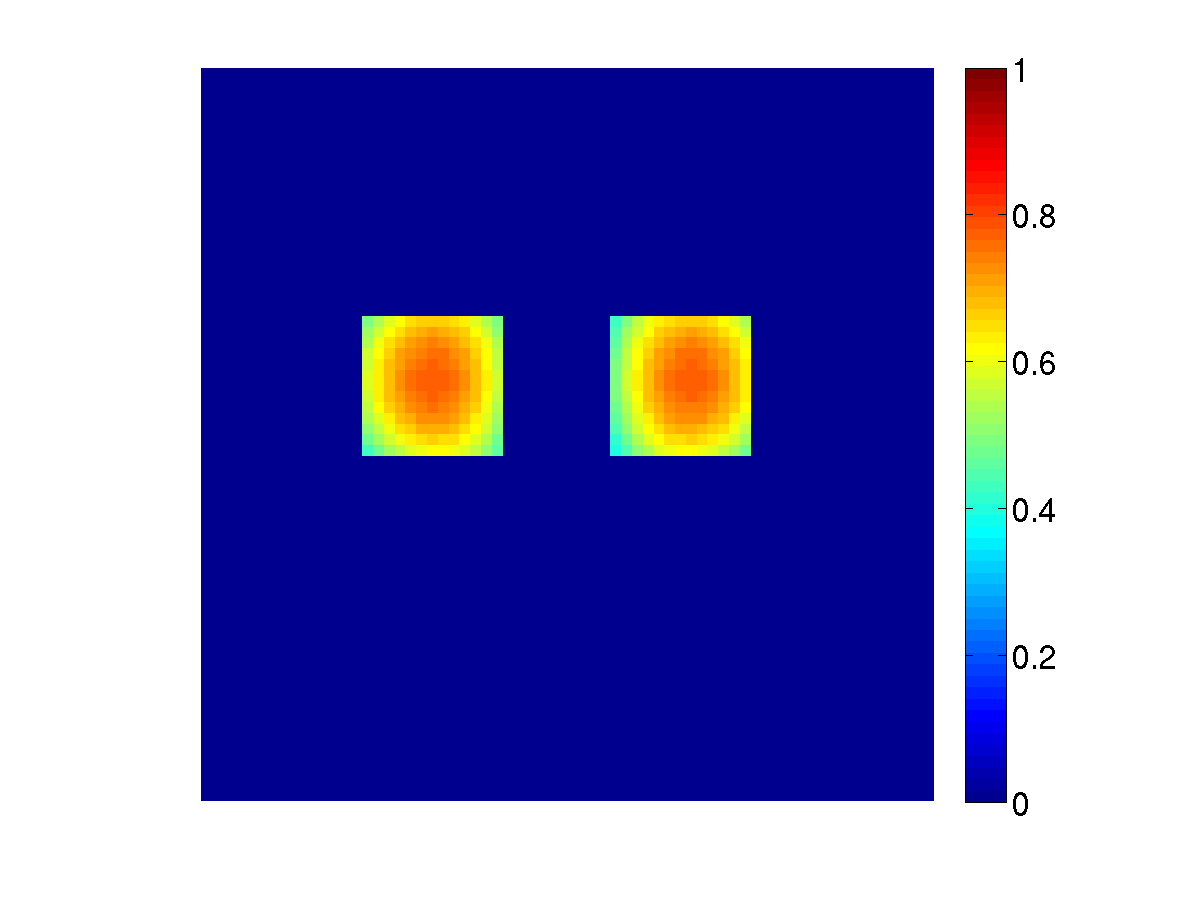} & \includegraphics[width=.25\textwidth]{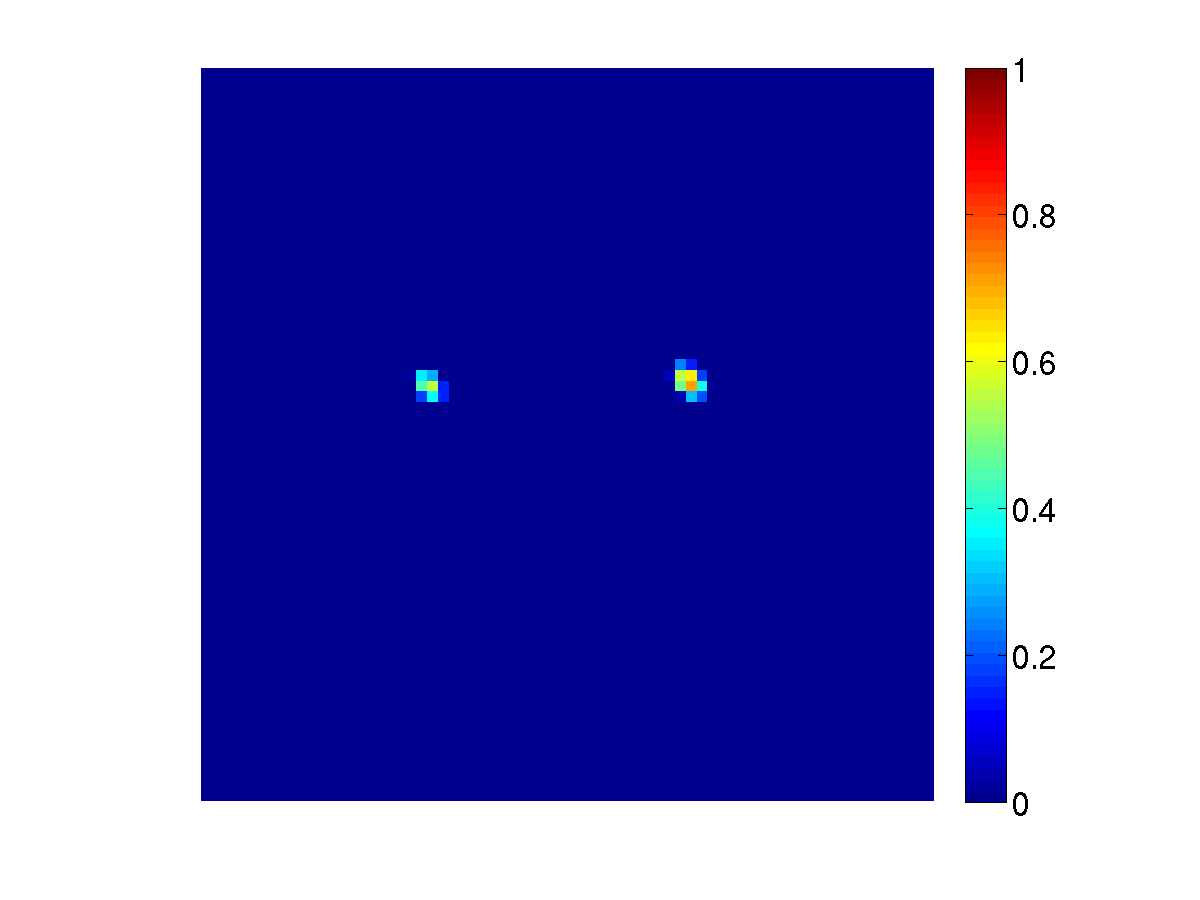}\\
     \includegraphics[width=.25\textwidth]{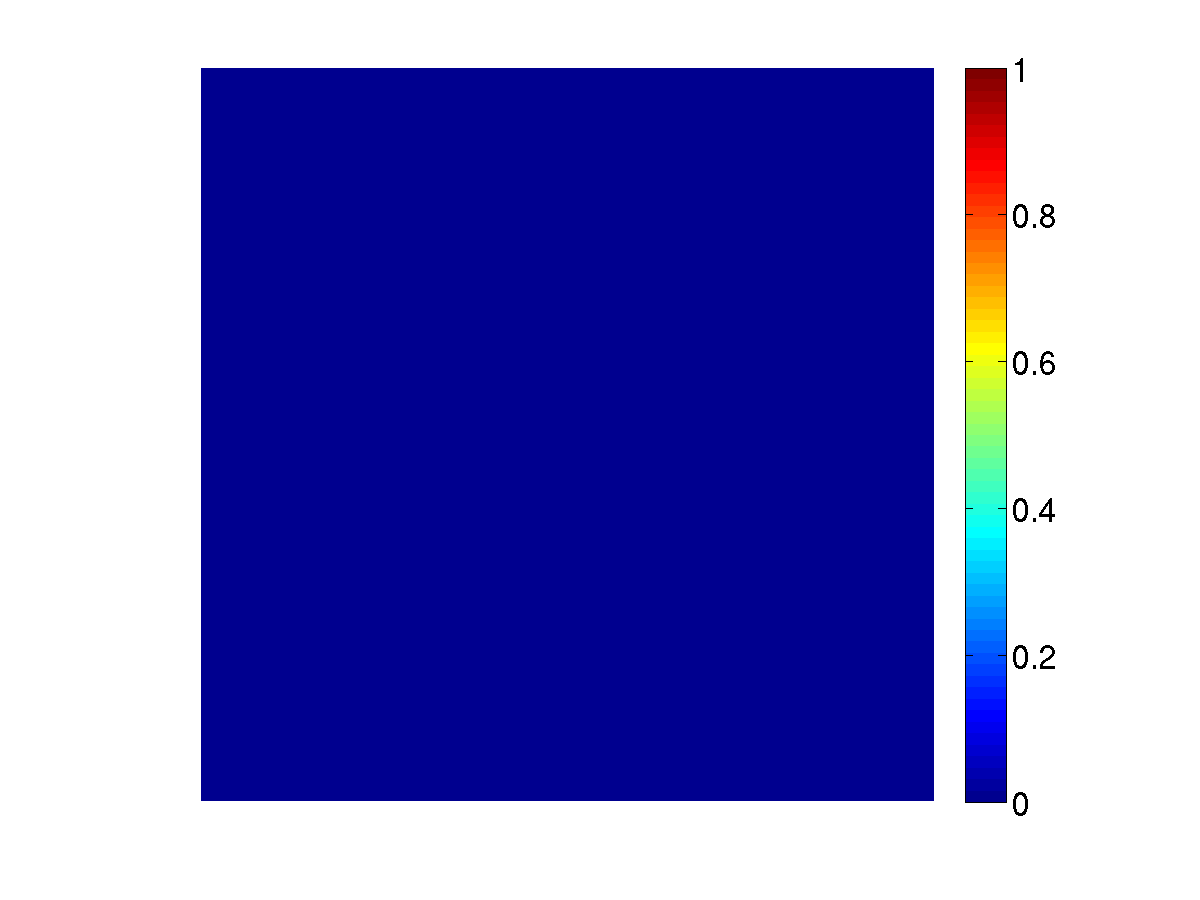} & \includegraphics[width=.25\textwidth]{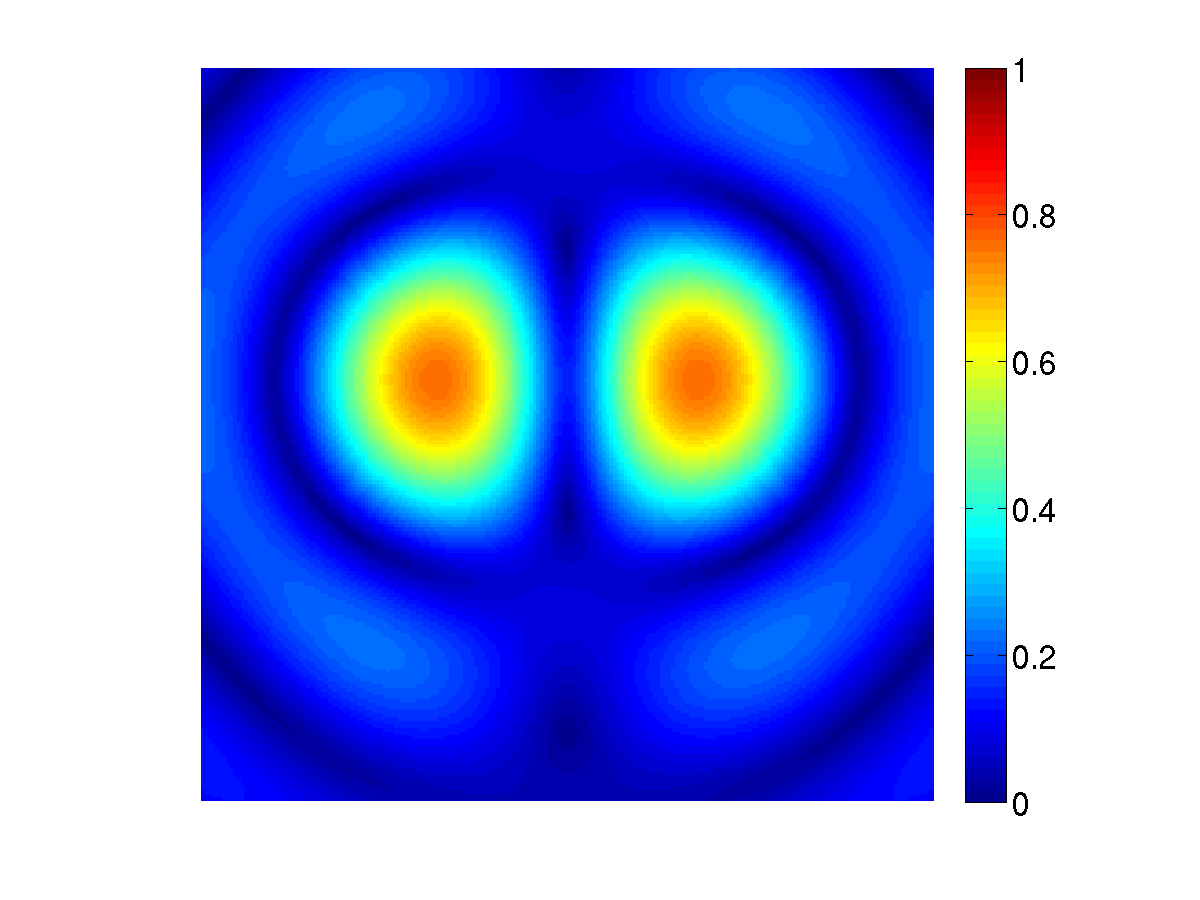}
    &\includegraphics[width=.25\textwidth]{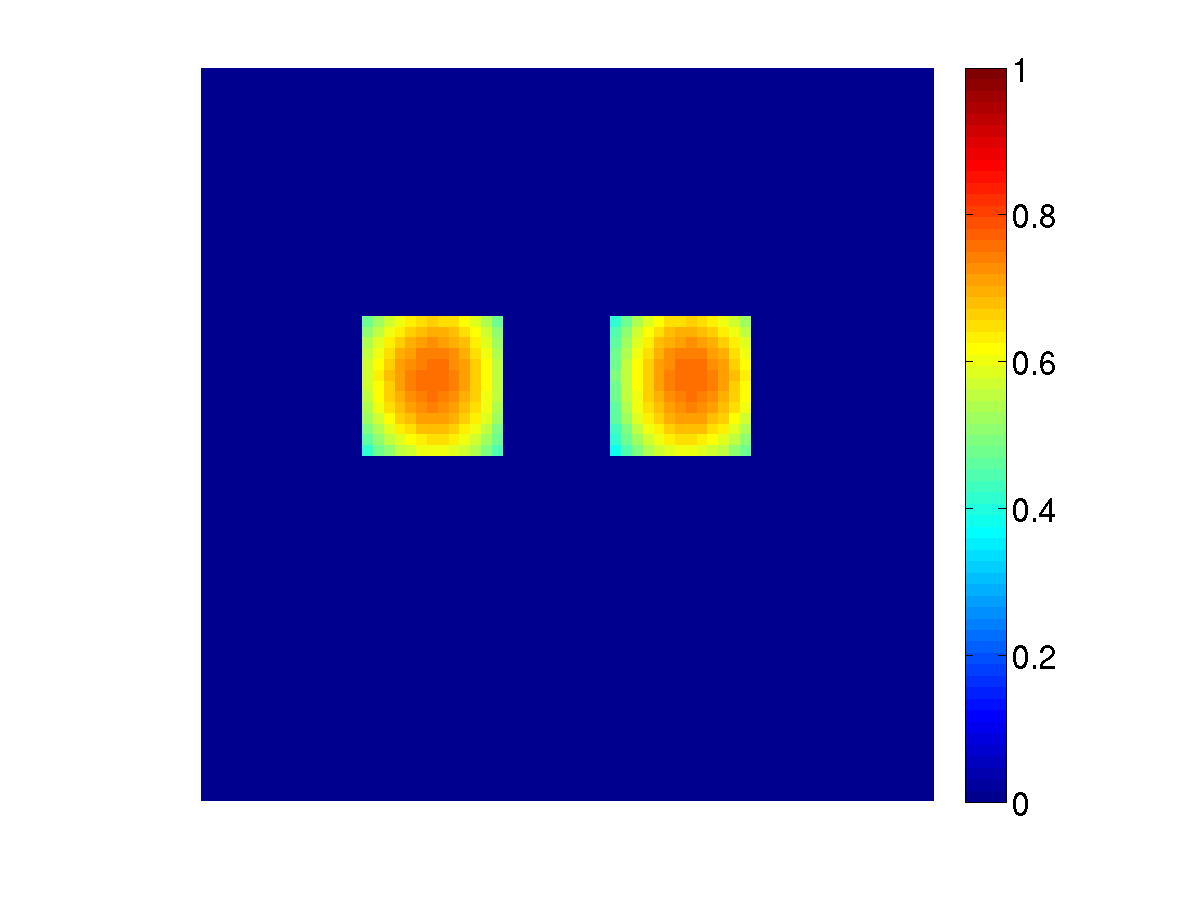}& \includegraphics[width=.25\textwidth]{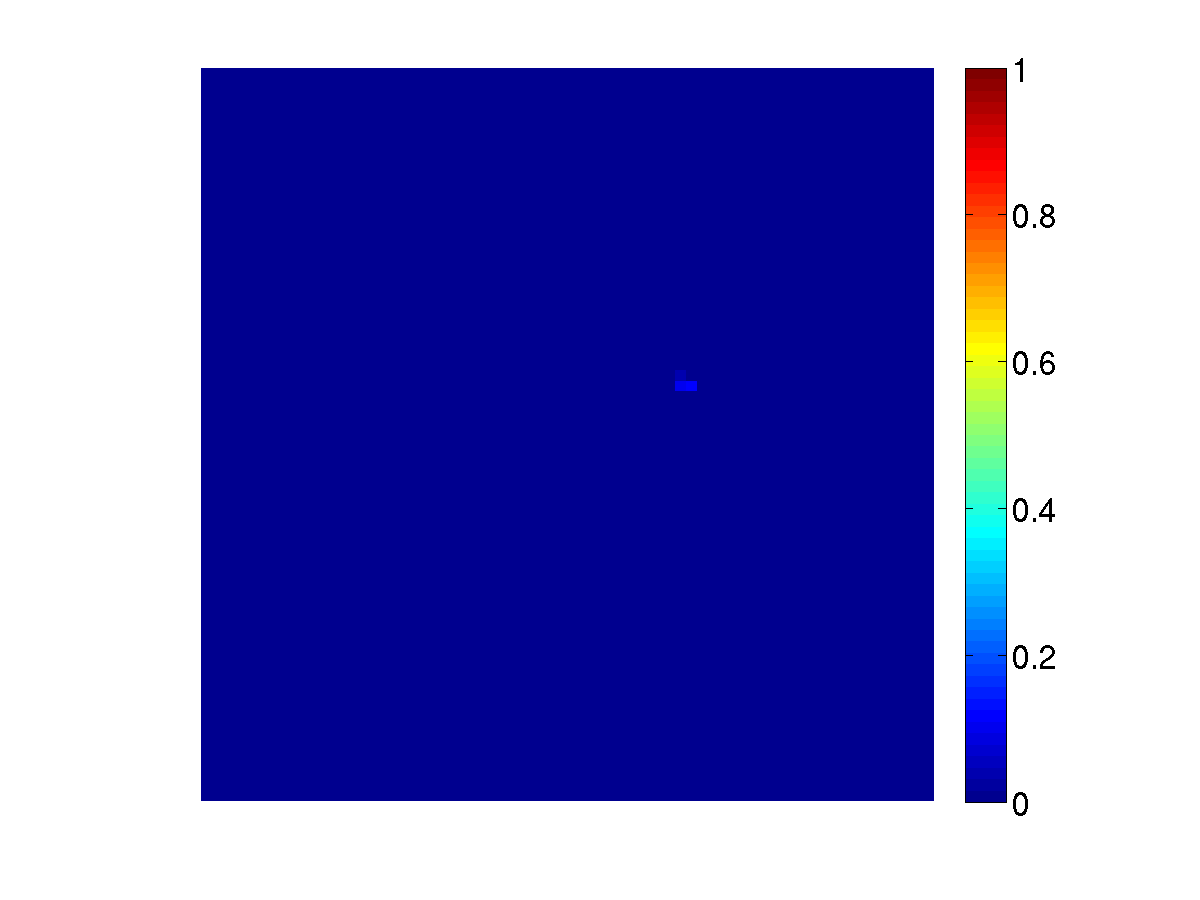}\\
    (a) true scatterer  & (b) index $\Phi$ & (c) index $\Phi|_D$ & (d) sparse recon.
  \end{tabular}
  \caption{Numerical results for Example \ref{exam:cube} with exact data: (a) true scatter, (b) index $\Phi$, (c)
  index $\Phi|_D$ (restriction to the subdomain $D$) and (d) sparse reconstruction. From the top to bottom: the cross sectional
  images at $x_2=1.07, 1.10,   1.13,   1.16$, $1.19$, and $1.22$, respectively.}
  \label{fig:cubex}
\end{figure}

\begin{figure}
  \centering
  \begin{tabular}{cccc}
     \includegraphics[width=.25\textwidth]{exam3d_rec_true_5.png} & \includegraphics[width=.25\textwidth]{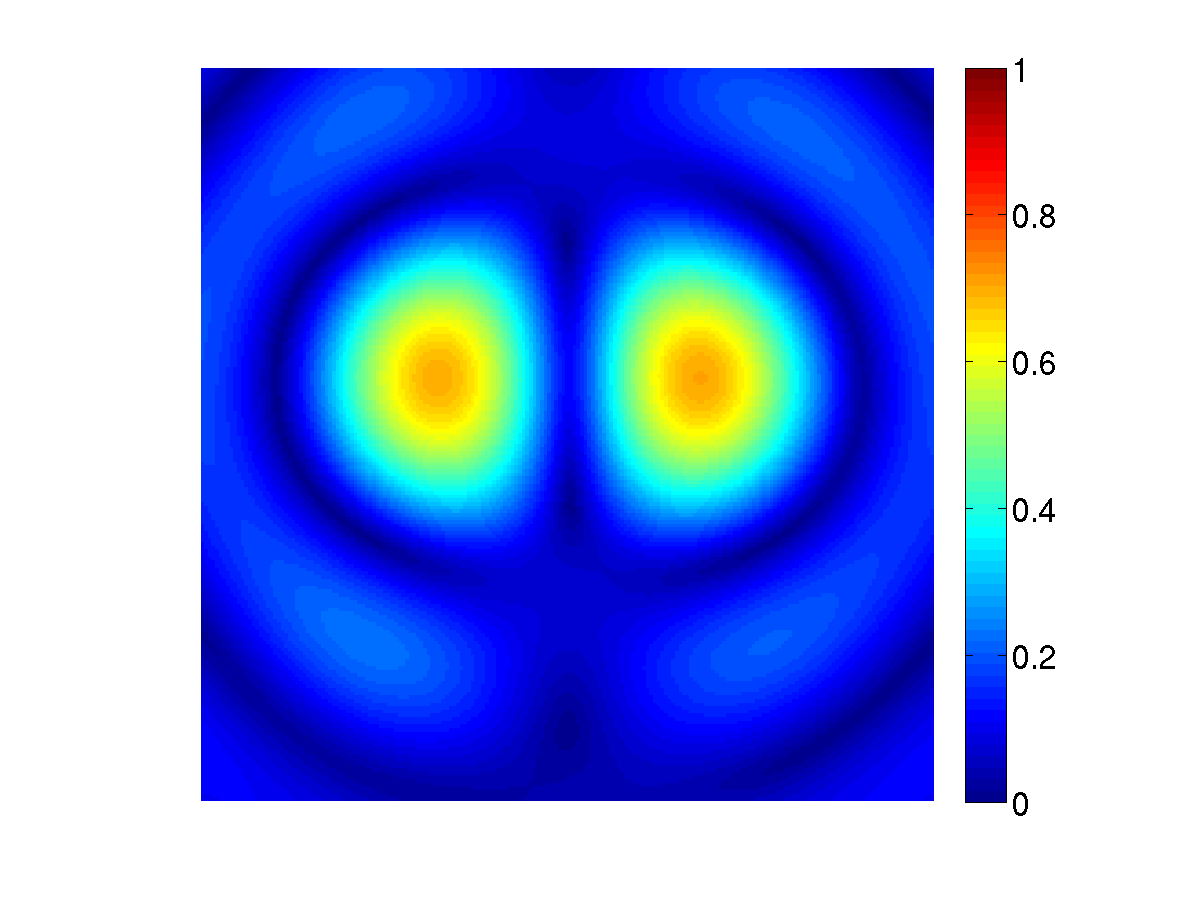}
    &\includegraphics[width=.25\textwidth]{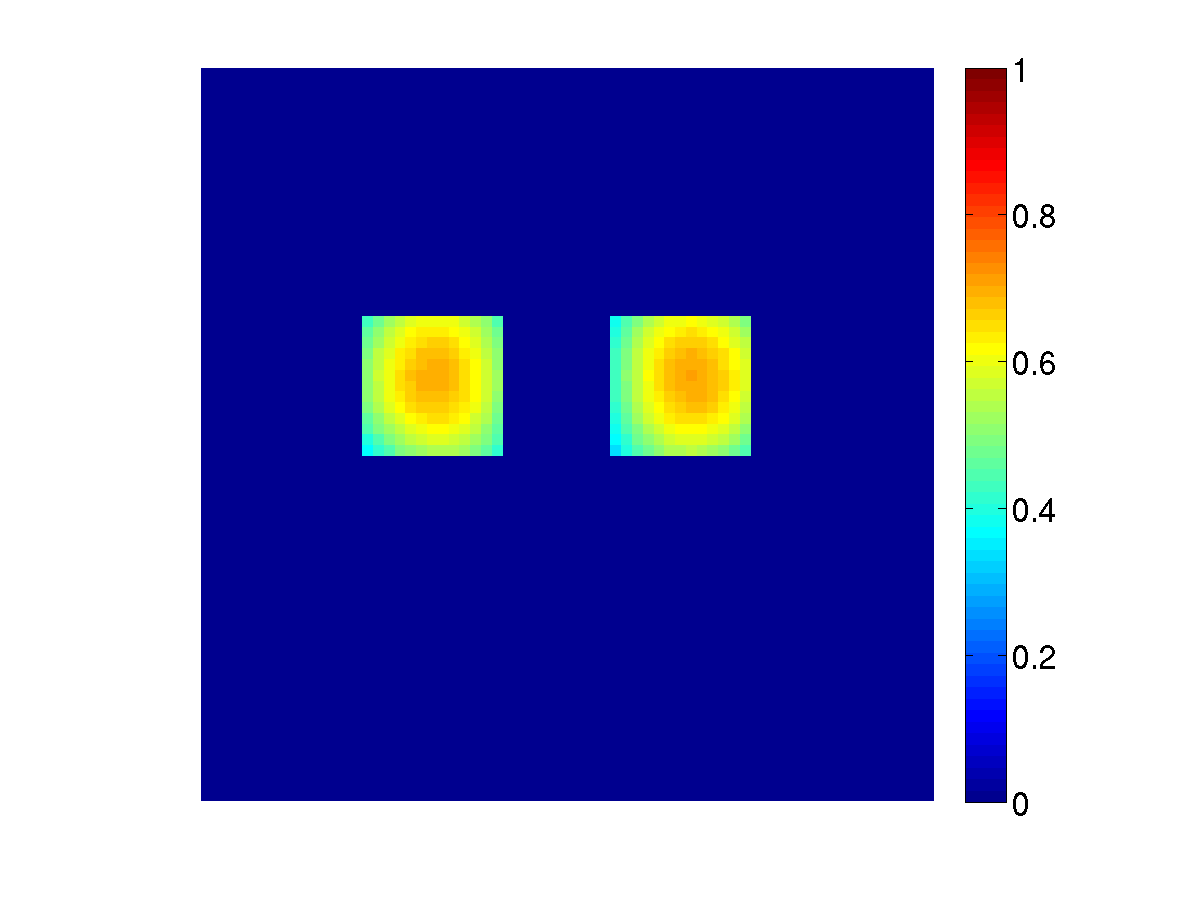} & \includegraphics[width=.25\textwidth]{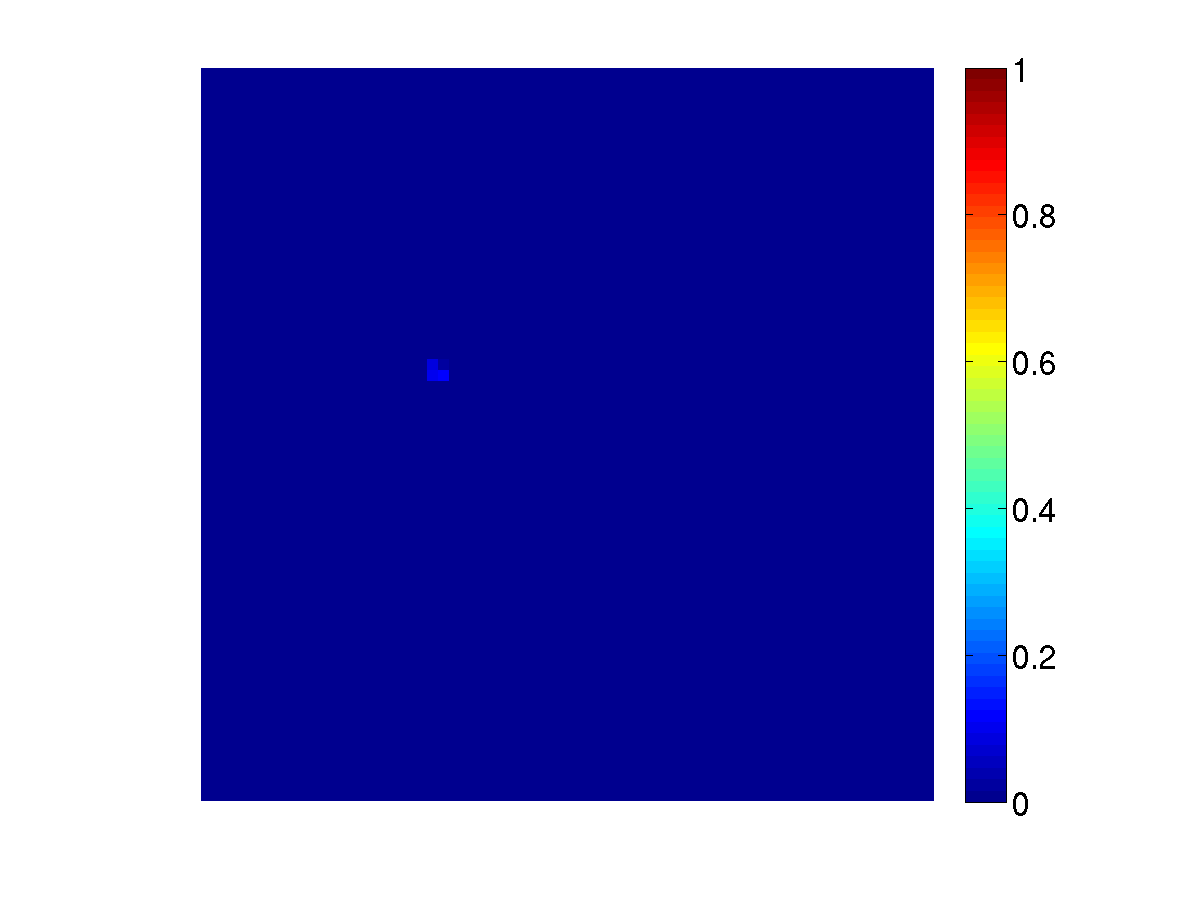}\\
     \includegraphics[width=.25\textwidth]{exam3d_rec_true_6.png} & \includegraphics[width=.25\textwidth]{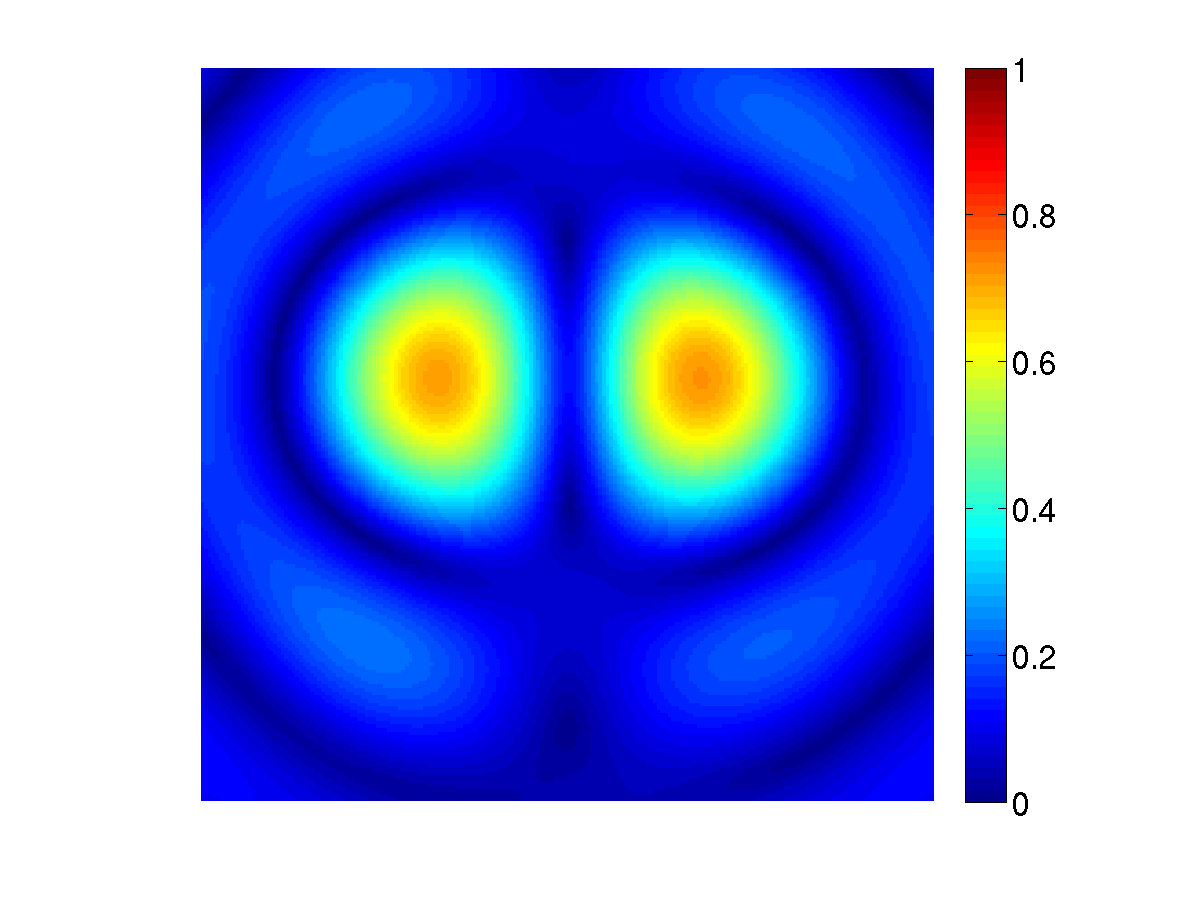}
    &\includegraphics[width=.25\textwidth]{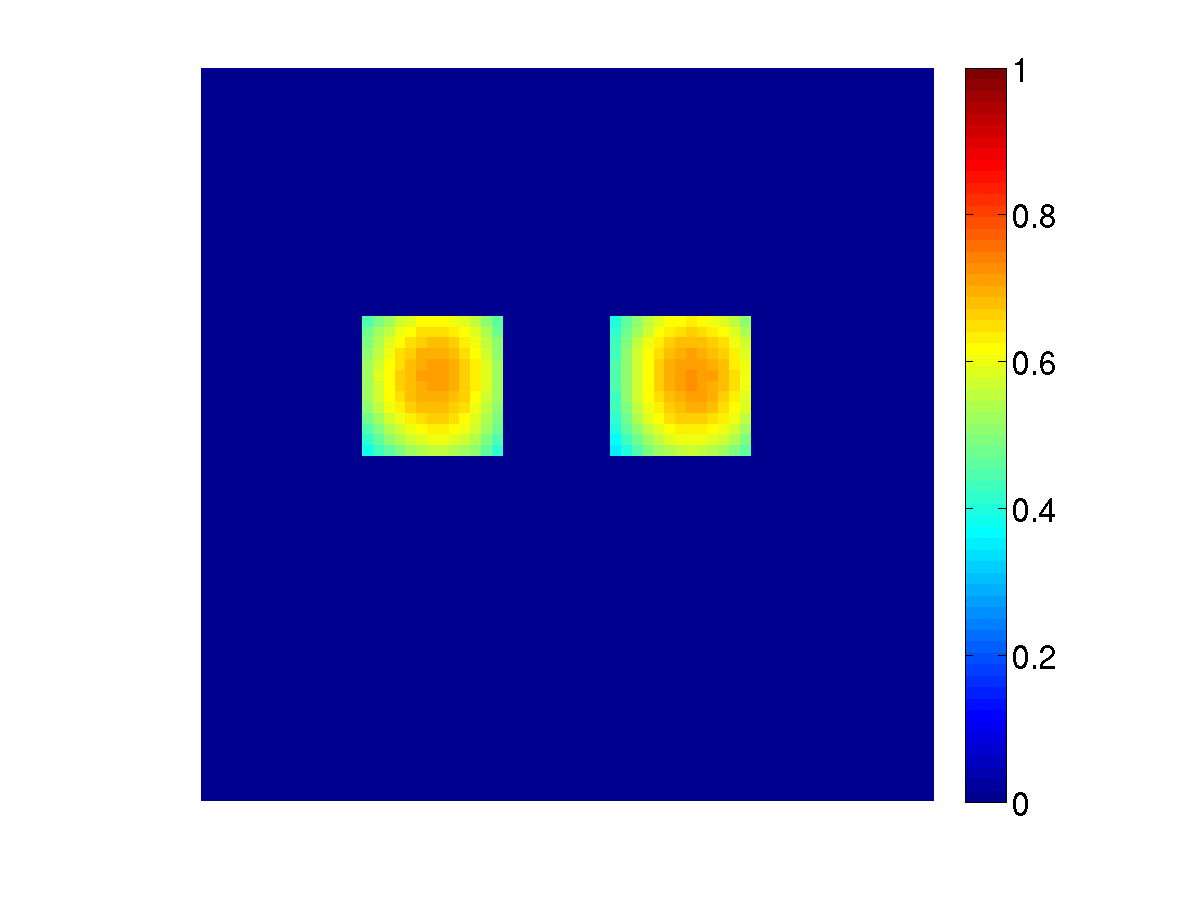} & \includegraphics[width=.25\textwidth]{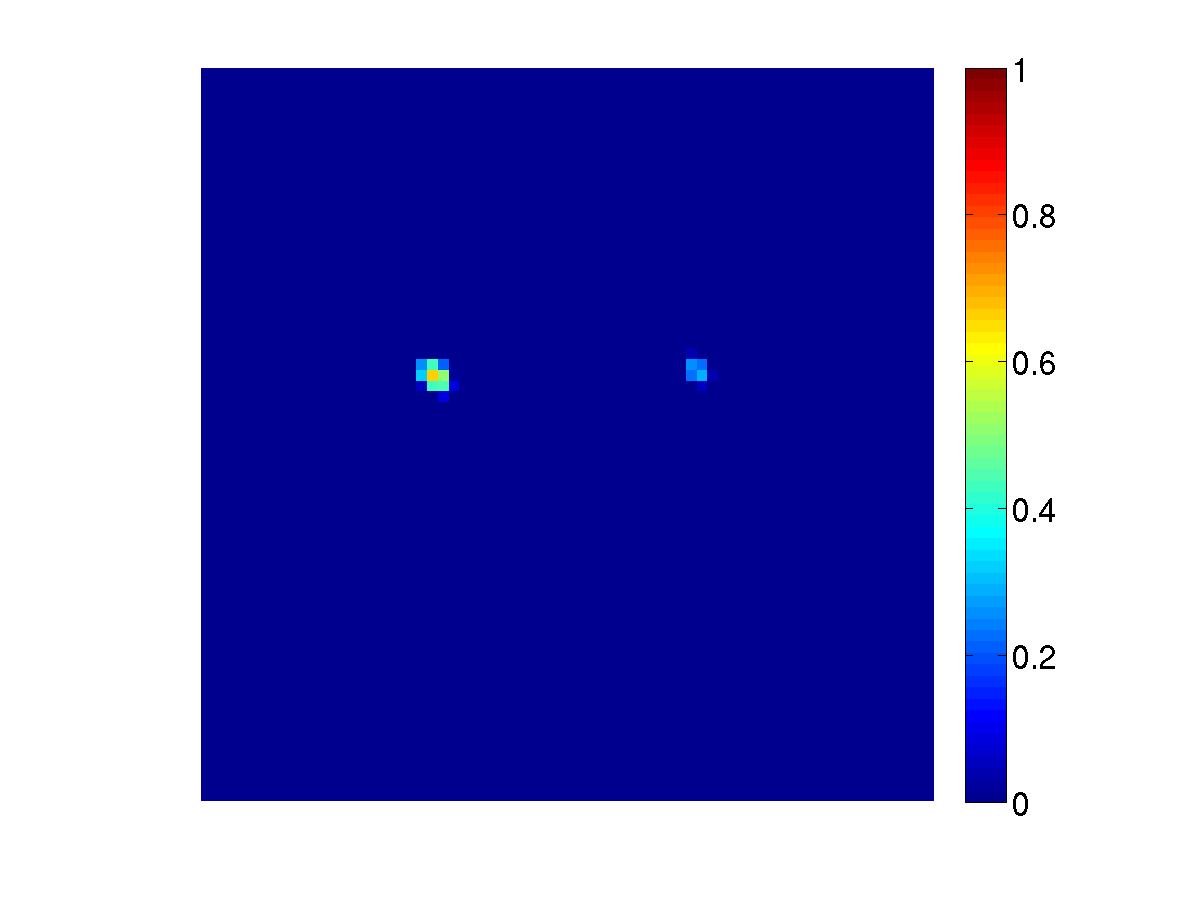}\\
     \includegraphics[width=.25\textwidth]{exam3d_rec_true_7.png} & \includegraphics[width=.25\textwidth]{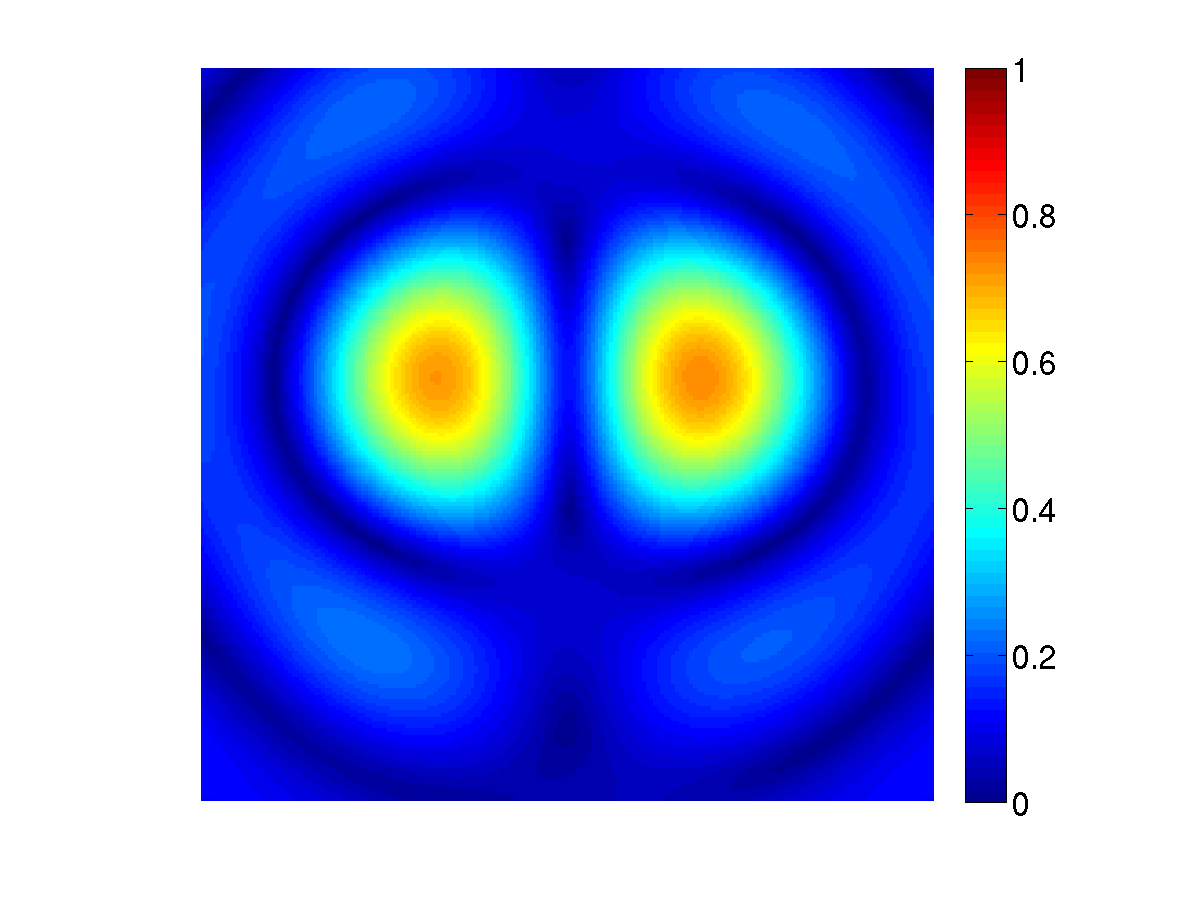}
    &\includegraphics[width=.25\textwidth]{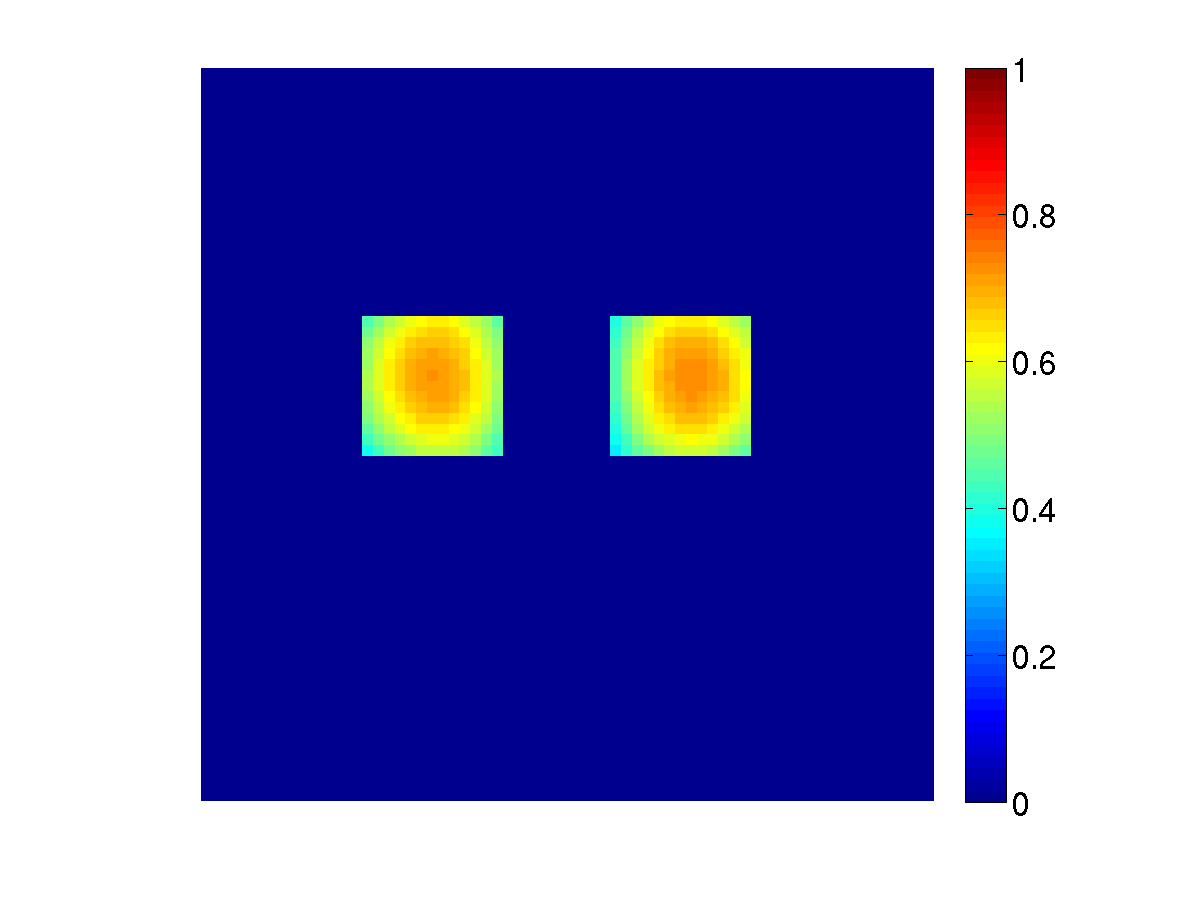} & \includegraphics[width=.25\textwidth]{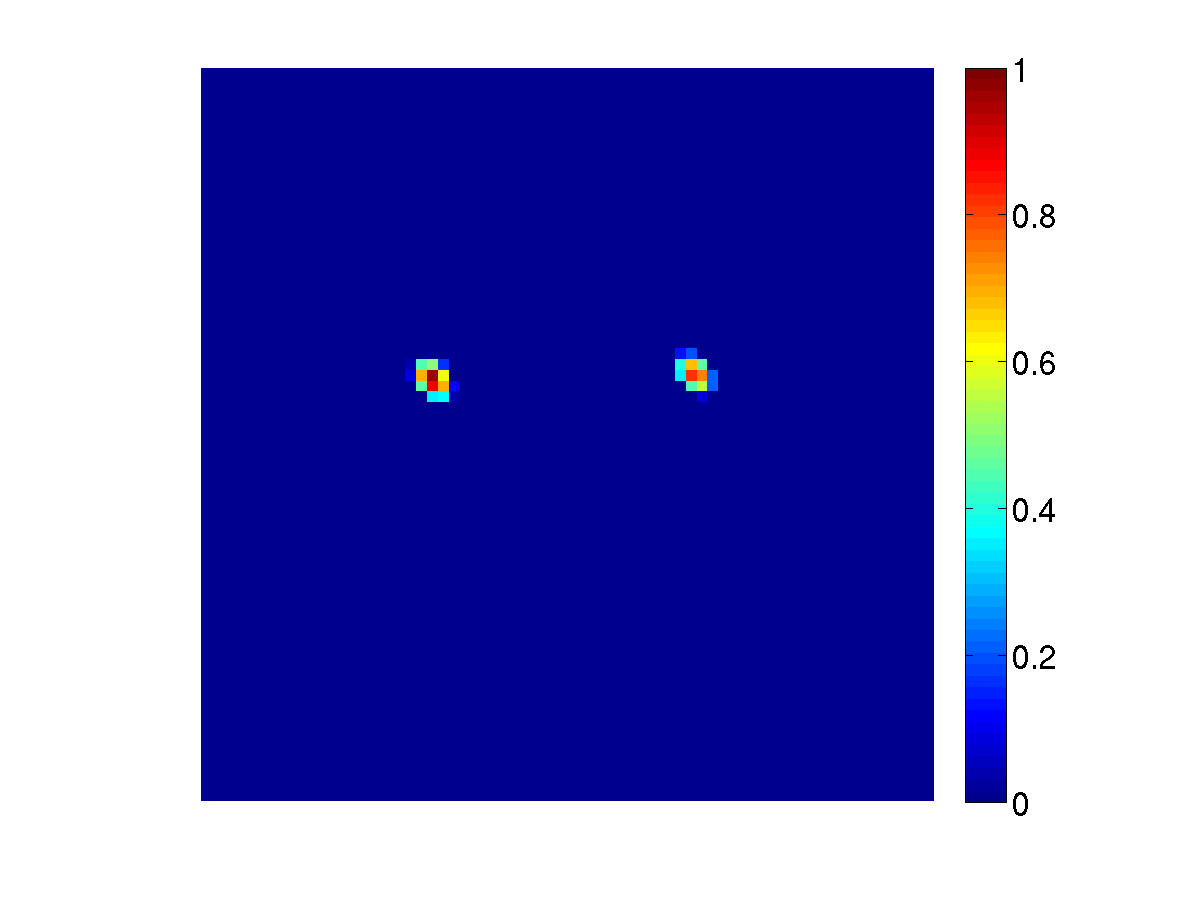}\\
     \includegraphics[width=.25\textwidth]{exam3d_rec_true_8.png} & \includegraphics[width=.25\textwidth]{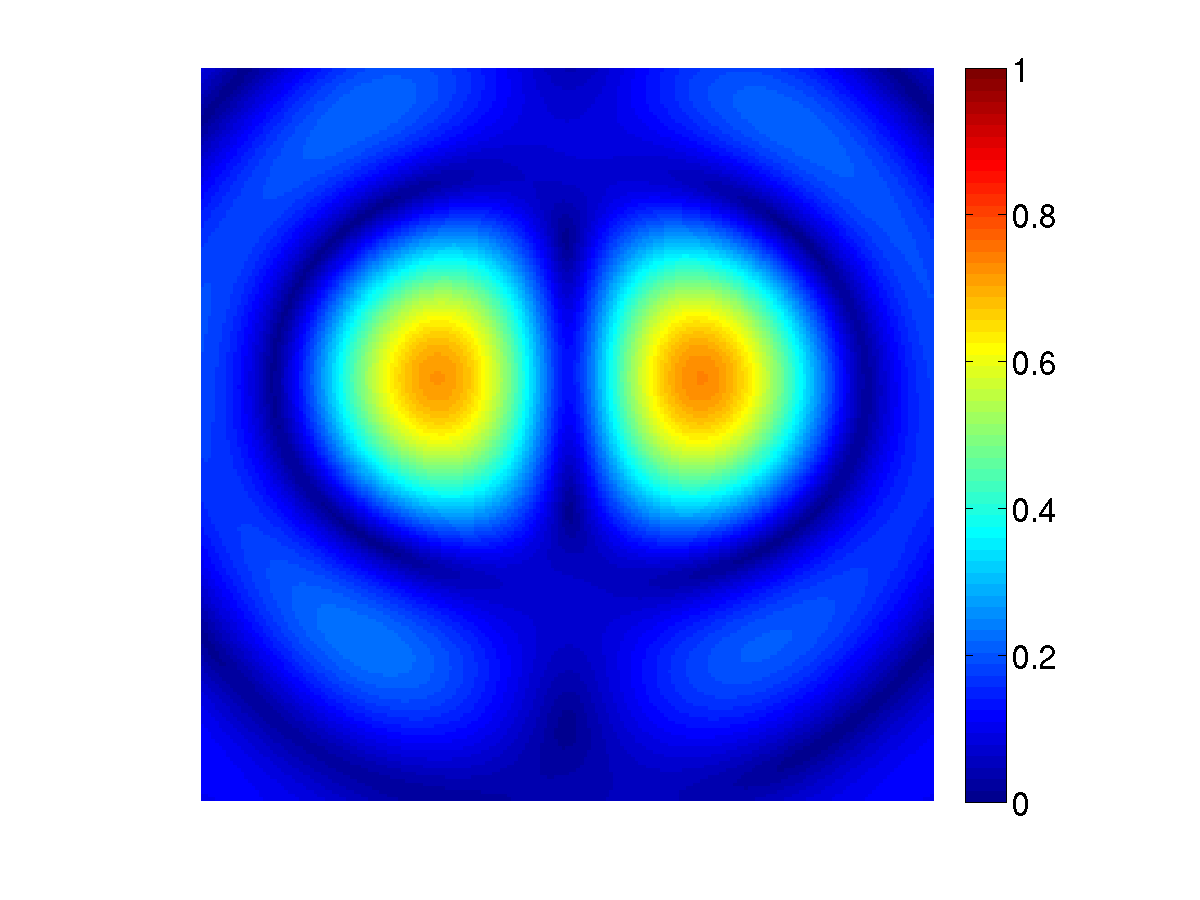}
    &\includegraphics[width=.25\textwidth]{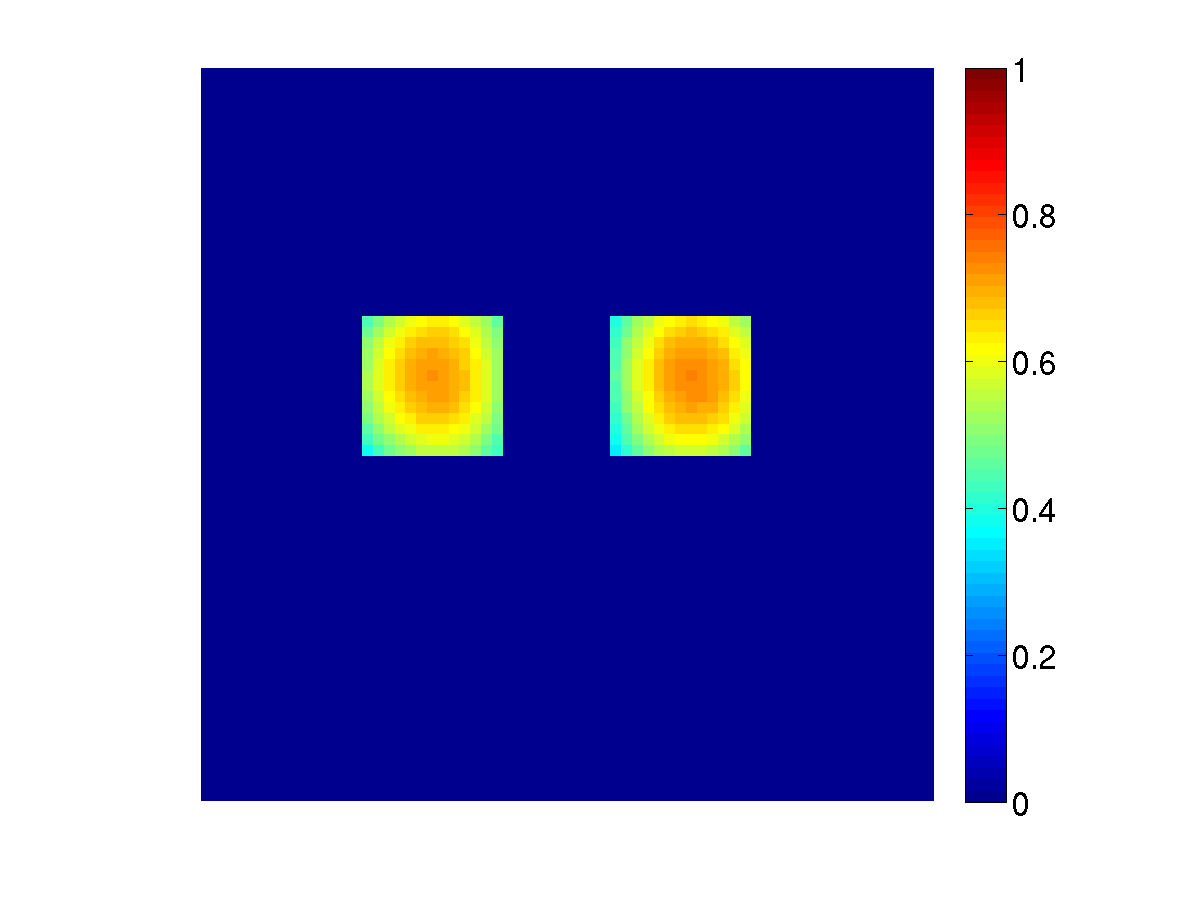} & \includegraphics[width=.25\textwidth]{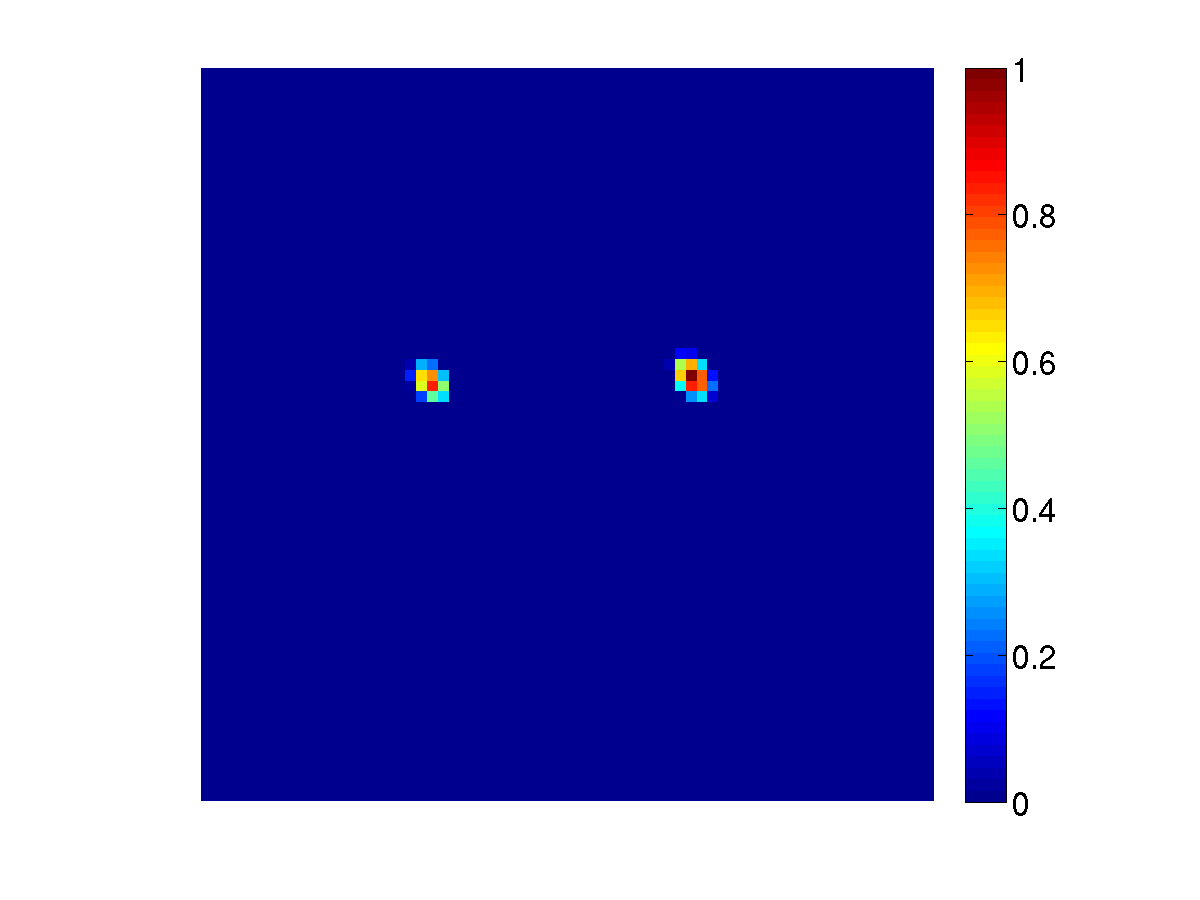}\\
     \includegraphics[width=.25\textwidth]{exam3d_rec_true_9.png} & \includegraphics[width=.25\textwidth]{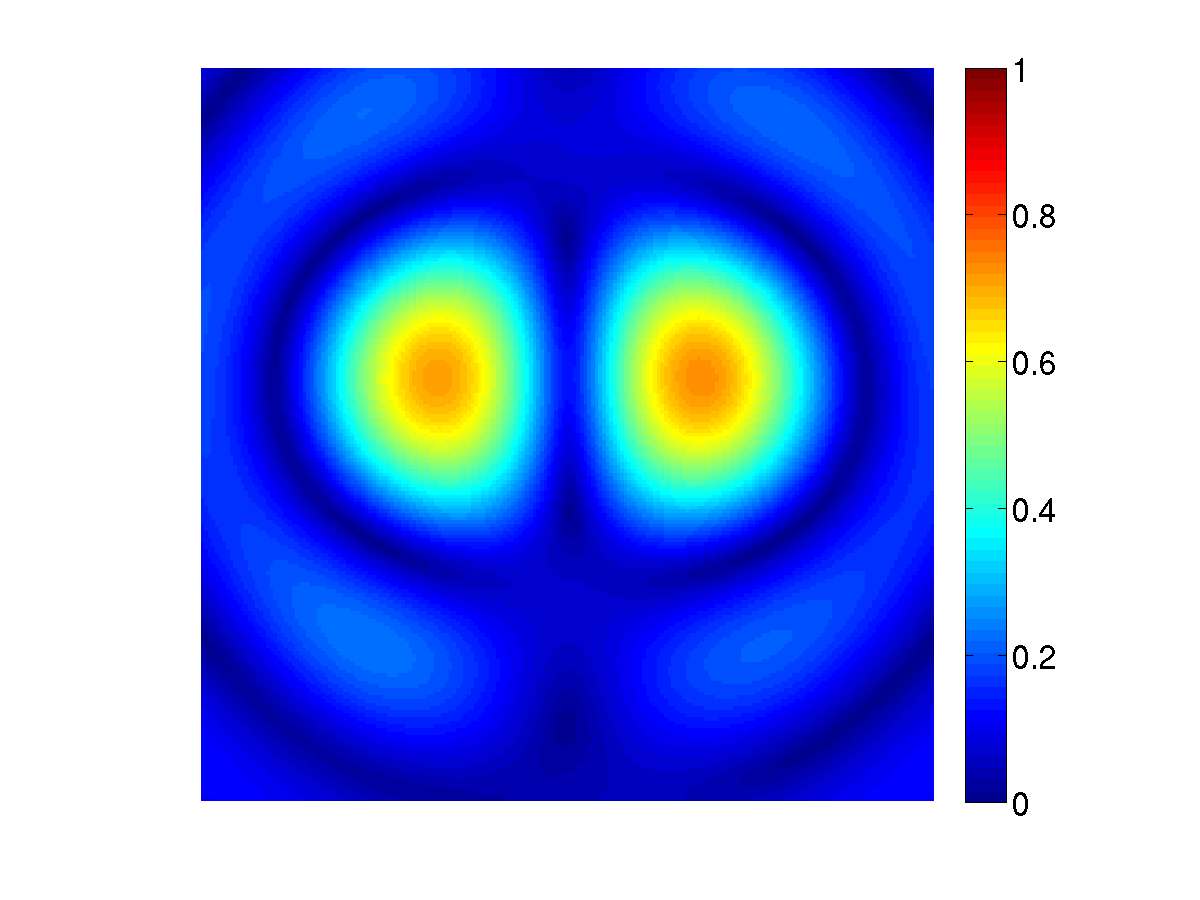}
    &\includegraphics[width=.25\textwidth]{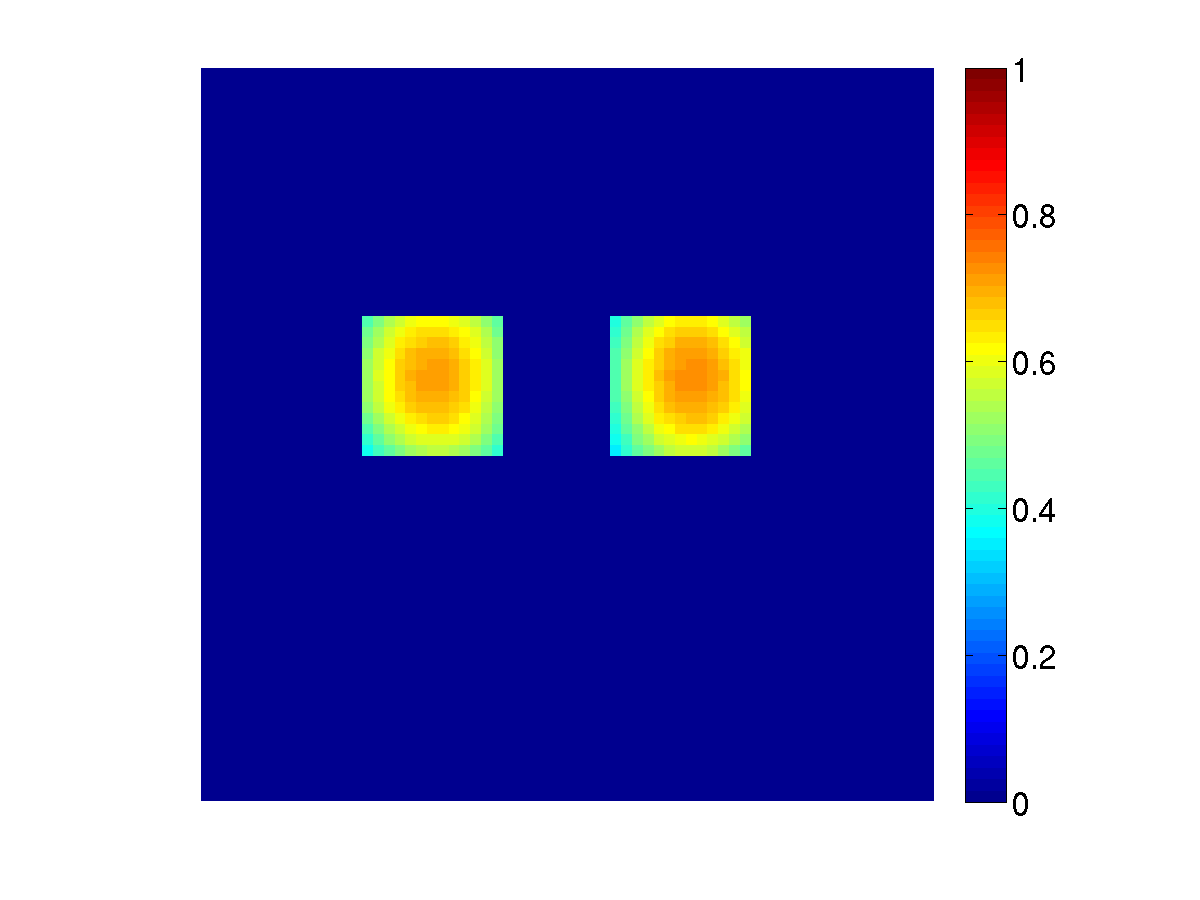} & \includegraphics[width=.25\textwidth]{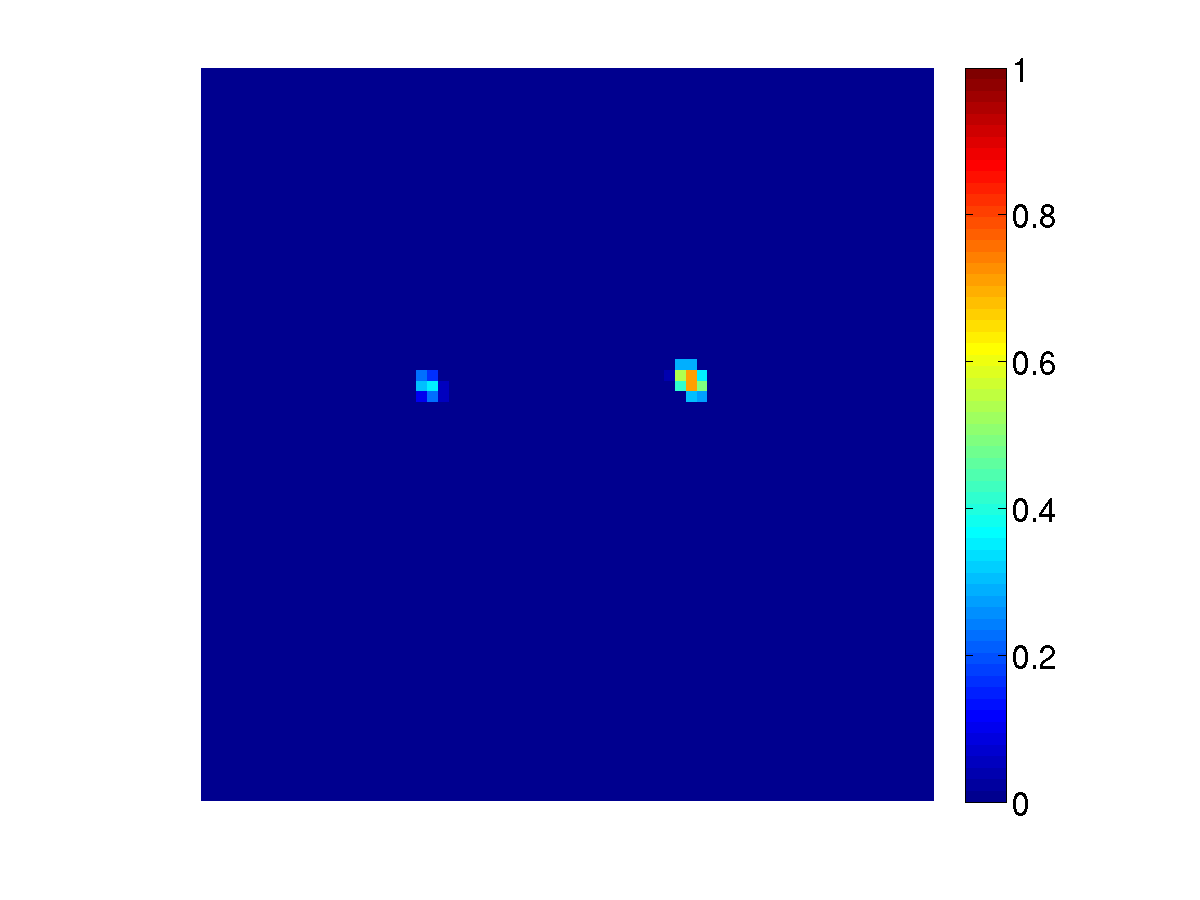}\\
     \includegraphics[width=.25\textwidth]{exam3d_rec_true_10.png} & \includegraphics[width=.25\textwidth]{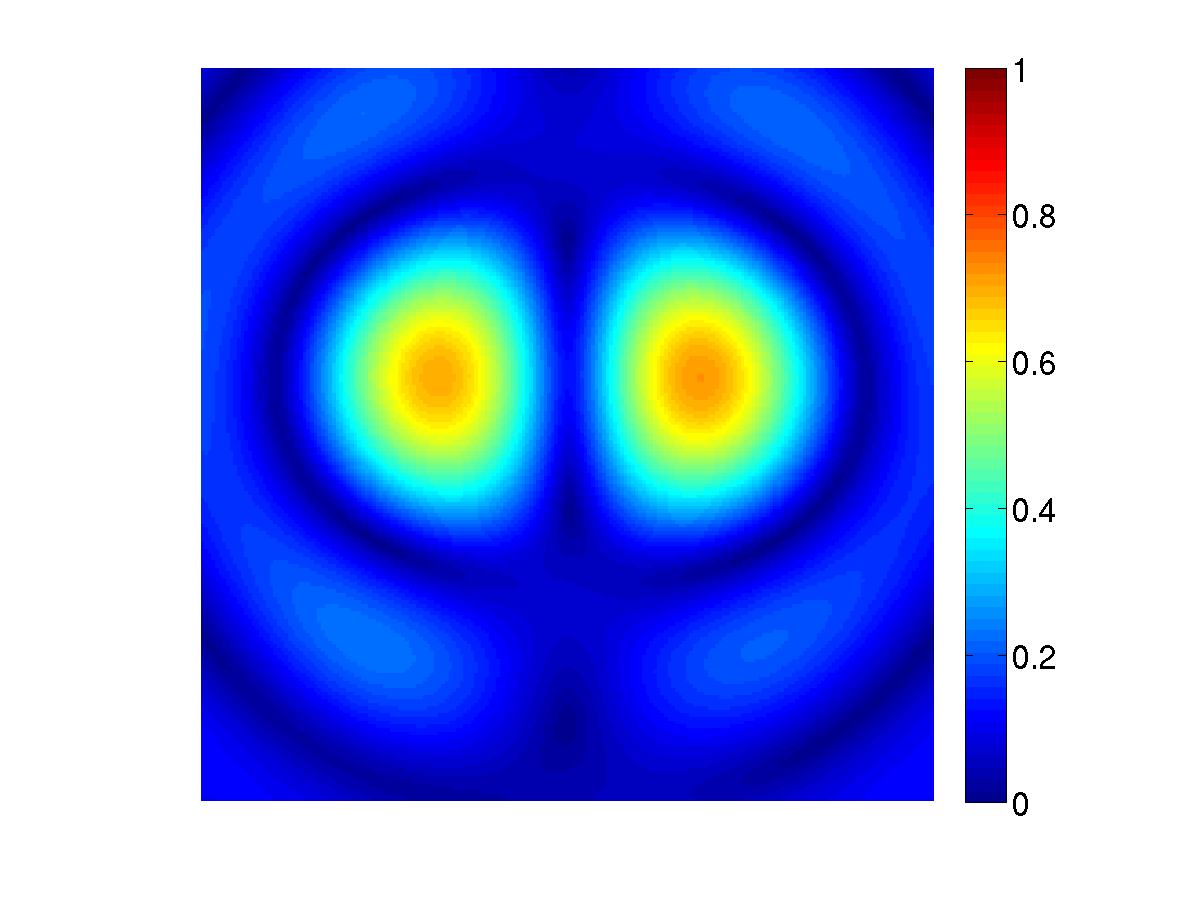}
    &\includegraphics[width=.25\textwidth]{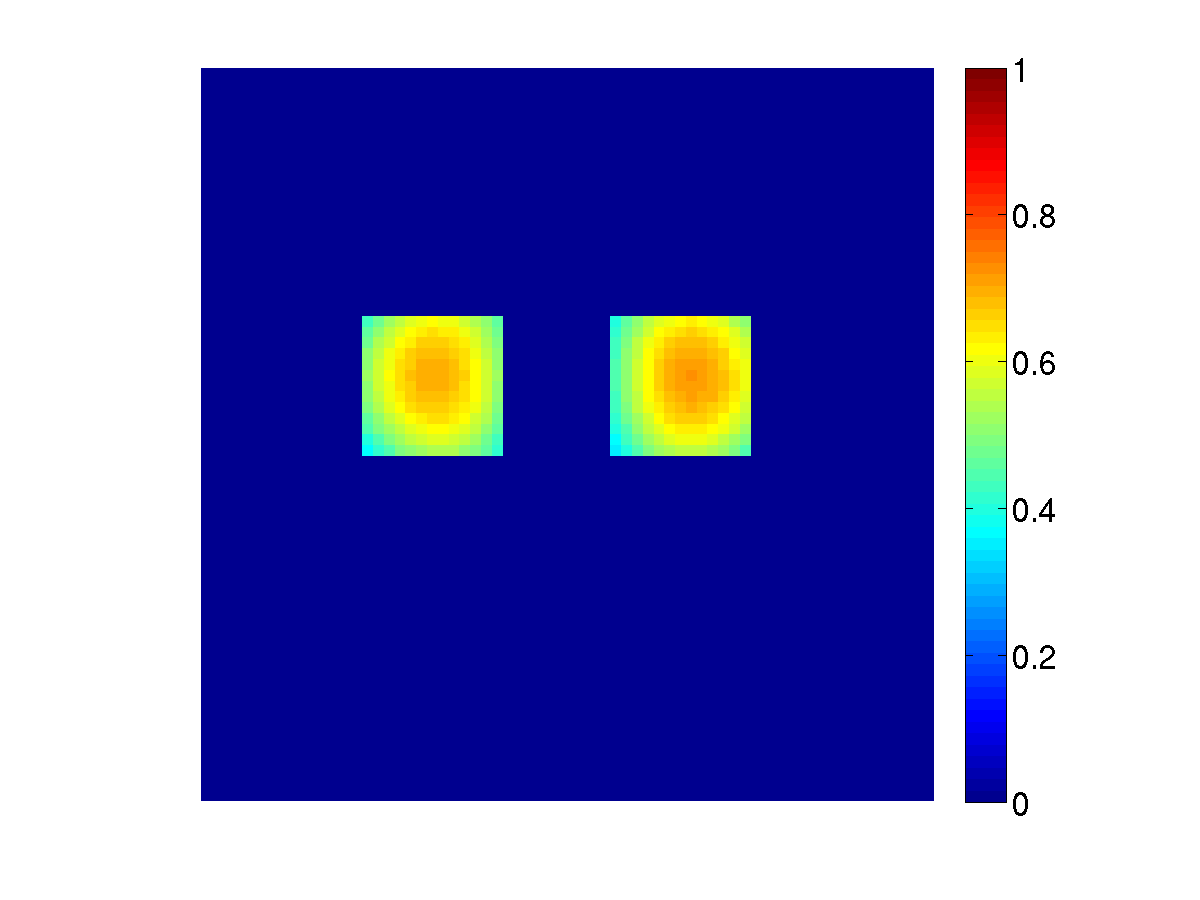}& \includegraphics[width=.25\textwidth]{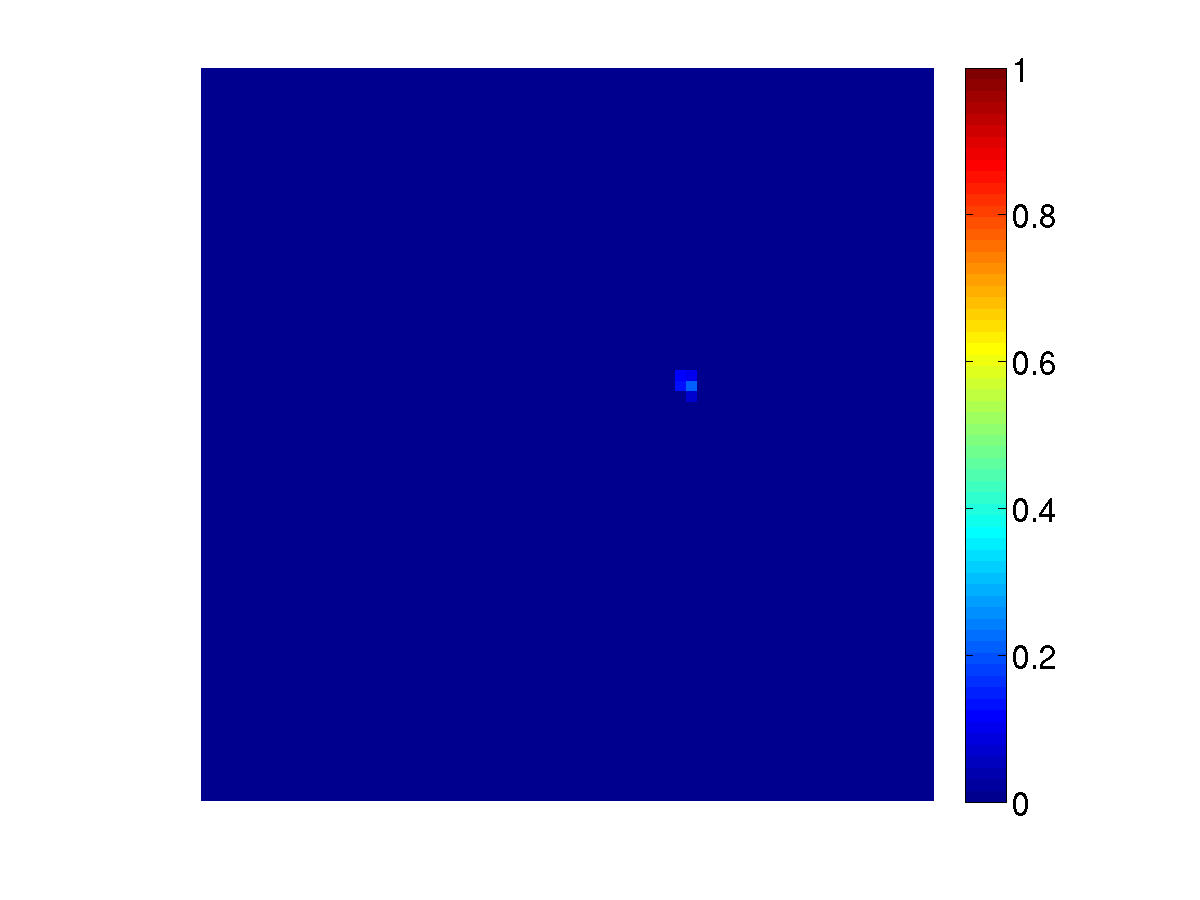}\\
    (a) true scatterer  & (b) index $\Phi$ & (c) index $\Phi|_D$ & (d) sparse recon.
  \end{tabular}
  \caption{Numerical results for Example \ref{exam:cube} with $20\%$ noise in the data: (a) true scatter,
  (b) index $\Phi$, (c) index $\Phi|_D$ (restriction to the subdomain $D$) and (d) sparse reconstruction. From the top to
  bottom: the cross sectional images at $x_2=1.07, 1.10, 1.13, 1.16$, $1.19$, and $1.22$ respectively.}
  \label{fig:cubn20}
\end{figure}

The numerical results for Example \ref{exam:cube} with exact data are shown in Fig.
\ref{fig:cubex}(b), where each row represents a cross-sectional image along the second
coordinate axis $x_2$. The scatterer support estimated by the index $\Phi$ agrees
reasonably with the exact one, and away from the boundary of the true scatterers, the
magnitude of $\Phi$ decreases quickly. However, the reconstructed profile is slightly
diffusive in comparison with the exact one, which is reminiscent of the decay property of
fundamental solutions. The nonsmooth mixed regularization \eqref{mixed} is carried out on
two cubic subregions (of width $0.36\lambda$), cf. Fig. \ref{fig:cubex}(c). Like before,
a significant improvement in the resolution is observed: the sparse estimate is much more
localized in comparison with the index $\Phi$, and also the magnitude is close to the
exact one; see Fig. \ref{fig:cubex}(d). The presence of $20\%$ data noise does not worsen
much the index $\Phi$ and the sparse reconstruction, cf. Fig. \ref{fig:cubn20}. Hence the
reconstruction algorithm is highly tolerant with respect to data noise.

Lastly, we briefly comment on the computational efficiency of the overall procedure. The
first step with the index involves only computing inner products and is embarrassingly
cheap and easily parallelized. The accuracy of the support detection is quite
satisfactory, and thus a large portion of the sampling domain $\widetilde{\Omega}$ can be
pruned from inversion, i.e., $|D|\ll|\widetilde{\Omega}|$. Hence, the enhancement via
mixed regularization is also rather efficient.

\section{Concluding remarks}
We have presented a novel two-stage inverse scattering method for the inverse medium
scattering problem of recovering the refractive index from near-field scattered data. The
efficiency and accuracy of the method stem from accurate support detection by the
sampling strategy and group sparsity-promoting of the mixed regularization technique. The
former is computationally very efficient, and reduces greatly the computational domain
for the more expensive inversion via nonsmooth mixed regularization, while the latter
achieves an enhanced resolution with the magnitudes and sizes comparable with the exact
ones. The numerical results for two- and three-dimensional examples clearly confirm these
observations.

These promising experimental results raise a number of interesting questions for further
studies. First, the potentials of mixed regularization have been clearly demonstrated. It
is of great interest to shed theoretical insights into the model as well as to design
efficient acceleration strategies, which for three-dimensional problems remains very
challenging. Some partial theoretical results can be found in \cite{ItoJinTakeuchi:2011}.
Also of much practical relevance is an automated choice of regularization parameters.
Second, the reconstructions were obtained with the linearized model, which represents
only an approximation to the genuine nonlinear IMSP model. It would be interesting to
justify the excellent performance of the linearization procedure. Third, the robustness
of the approach to noise is outstanding when compared with more conventional inverse
scattering algorithms, especially noting the limited data for inversion. The mechanism of
the robustness is not yet clear.

\section*{Acknowledgements}
The work of BJ is supported by Award No. KUS-C1-016-04, made by King Abdullah University
of Science and Technology (KAUST), and that of JZ is substantially supported by Hong Kong
RGC grants (projects 405110 and 404611).

\appendix
\section{Numerical method for forward scattering}\label{app:int}

We denote by $\mathbb{J}$ the index set of grid points of a uniformly distributed mesh
with a mesh size $h>0$ and consider square cells
\begin{equation*}
  B_j=B_{j_1,j_2}=(x^1_{j_1},x^2_{j_2})+[-\tfrac{h}{2},\tfrac{h}{2}]\times [\tfrac{h}{2},\tfrac{h}{2}]
\end{equation*}
for every tuple $j=(j_1,j_2)$ belonging to the index set $\mathbb{J}$. Assume that the
domain $\cup_{j \in \mathbb{J}} \,B_j$ contains the scatterer support $\Omega$. We use
the mid-point quadrature rule to evaluate the operator $K$, and hence the integral
\eqref{eqn:indcur} is approximated by
\begin{equation*}
  I_k - \eta_{k}\, \sum_{j \in \mathbb{J}} G_{k,j} I_jh^2 = \eta_k \, u^{inc}(x_k)
\end{equation*}
where $I_k=I(x_k)$ and $\eta_k=\eta(x_k)$, and the off-diagonal entries $G_{k,j}$ and the
diagonal entries $G_{k,k}$ are given by $G_{k,j}=G(x_k,x_j)$ and
\begin{equation*}
  G_{k,k}=\frac{1}{h^2} \int_{(-\tfrac{h}{2},\tfrac{h}{2})^2}  G(x,0)dx,
\end{equation*}
respectively. The diagonal entries can be accurately computed by tensor-product Gaussian quadrature rules. 
The resulting system can be solved using standard numerical solvers, e.g., Gaussian
elimination, if the cardinality of the index set $\mathbb{J}$ is medium, and iterative
solvers like GMRES. The extension of the procedure to 3D problems is straightforward.

\section{Semi-smooth Newton method}\label{app:ssn}
In this part, we derive a semi-smooth Newton method for minimizing \eqref{mixed}. The
optimality condition of the variational problem reads
\begin{equation*}
   \left\{\begin{aligned}
     K^*K\eta+\alpha\lambda-\beta\Delta \eta - K^\ast u^s  &= 0,\\
     \lambda -\frac{\lambda+c\eta}{\max(1,|\lambda+c\eta|)}&=0,
   \end{aligned}\right.
\end{equation*}
where $\lambda$ is the Lagrange multiplier (dual variable). The second line, the
complementarity function, equivalently expresses the inclusion $\lambda\in
\partial\|\eta\|_{L^1}$, {which can be checked directly by pointwise
inspection. Thereby, we effectively transforms the inclusion \eqref{fixed} into a
numerically amenable nonlinear system}.
It follows directly from the complementarity relation
\begin{equation}\label{eqn:complem}
  \lambda =\frac{\lambda+c\eta}{\max(1,|\lambda+c\eta|)}
\end{equation}
that on the active set $\mathcal{A}=\{x\in D: |\lambda+c\eta|(x)\leq1\}$, $\eta$ vanishes
identically. Otherwise, both the dual variable $\lambda$ and the primal variable $\eta$
need to be solved. We shall solve the system by a semi-smooth Newton method
\cite{ItoKunisch:2008}. First observe that the Newton step (with the increments for
$\lambda$ and $\eta$ denoted by $\delta\lambda$ and $\delta\eta$, respectively) applied
to the following reformulation of equation \eqref{eqn:complem} (on the set
$\mathcal{I}=D\setminus\mathcal{A}$)
\begin{equation*}
   \lambda|\lambda+c\eta| - \lambda + c\eta=0
\end{equation*}
is given by
\begin{equation*}
|\lambda+c\eta|\delta\lambda+\lambda\frac{\lambda+c\eta}{|\lambda+c\eta|}[\delta\lambda+c\delta\eta]-(\delta\lambda+c\delta
\eta)+\lambda|\lambda+c\eta|-(\lambda+c\eta)=0,
\end{equation*}
or equivalently with the notation $\lambda^+=\lambda+\delta\lambda$ and
$\eta^+=\eta+\delta\eta$, we have
\begin{equation*}
   \lambda^+|\lambda+c\eta|+\lambda\frac{\lambda+c\eta}{|\lambda+c\eta|}[\lambda^++c\eta^+]=\lambda|\lambda+c\eta|+[\lambda^++c\eta^+].
\end{equation*}
Next we apply the idea of damping and regularization to the equation and thus get
\begin{equation*}
\lambda^+|\lambda+c\eta|+\theta[\lambda^++c\eta^+]\frac{\lambda+c\eta}{|\lambda+c\eta|}\frac{\lambda}{\max(|\lambda|,1)}
=[\lambda^++c\eta^+]+\theta|\lambda+c\eta|\frac{\lambda}{\max(|\lambda|,1)}.
\end{equation*}
Here, the purpose of the regularization step $\frac{\lambda}{\max(|\lambda|,1)}$ is to
automatically constrain the dual variable $\lambda$ to $[-1,1]$. The damping factor
$\theta$ is automatically selected to achieve the stability. To this end, we let
$d=|\lambda+c\eta|$, $\widetilde{\eta}=d-1$, $a=\tfrac{\lambda}{\max(|\lambda|,1)}$, and
$b=\tfrac{\lambda+c\eta}{|\lambda+c\eta|}$. We arrive at
\begin{equation*}
\lambda^+(\widetilde{\eta}+1)+\theta[\lambda^++c\eta^+]ab=
[\lambda^++c\eta^+]+\theta ad.
\end{equation*}
Thus we have
\begin{equation*}
\lambda^+=\frac{1}{\widetilde{\eta}+\theta ab}[1-\theta ab]c\eta^++\frac{\theta
d}{\widetilde{\eta}+\theta ab}a
\end{equation*}
To arrive at a simple iteration scheme, we set $\tfrac{\theta d}{\widetilde{\eta}+\theta
ab}=1$, i.e., $\theta=\tfrac{d-1}{d-ab}\leq1$. Consequently, we obtain a simple
iteration
\begin{equation*}
   \lambda^+=\frac{1-ab}{d-1}c\eta^++\frac{\lambda}{\max(|\lambda|,1)},
\end{equation*}
where we have used the relation $\frac{1-\theta ab}{\widetilde\eta+\theta
ab}=\frac{1-ab}{d-1}$. Substituting this into the first equation gives
\begin{equation}\label{eqn:ssneq}
   K^\ast K\eta^+ + \alpha c\frac{1-ab}{d-1} \eta^+ -\beta \Delta \eta^+ = K^\ast u^s-\alpha\frac{\lambda}{\max(|\lambda|,1)}.
\end{equation}
We note that one only needs to solve equation \eqref{eqn:ssneq} on the inactive set
$\mathcal{I}$, since on the active set $\mathcal{A}$, there always holds $\eta^+=0$. This
has an enormous computational consequence: the size of the linear system in
\eqref{eqn:ssneq} can be very small if $|\mathcal{I}|$ is small, i.e., the solution is
sparse. This last relation shows also clearly the sparsity of the solution, and this
provides a crispy estimate of the background. Upon obtaining the solution $\eta^+$, one
can update $\lambda^+$ on the sets $\mathcal{I}$ and $\mathcal{A}$ according to the
second and the first equation, respectively. Lastly, we would like to remark on the
consistency of the scheme: if the sequence generated by the semi-smooth Newton method
converges, then the limit satisfies the complementarity relation \eqref{eqn:complem} as
desired.

\bibliographystyle{abbrv}
\bibliography{scatter}

\end{document}